\documentclass[11pt]{article}

\usepackage{amsthm}
\usepackage{hyperref}
\input{pstricks.tex}

\newtheorem{thm}{Theorem}
\newtheorem{lem}{Lemma}
\newtheorem{cor} {Corollary}

\newtheorem{pro} {Proposition}
\newtheorem{df} {Definition}

\newtheorem{obs}{Observation}
\newtheorem{cons}{Construction}

\newcommand {\F}{{\mathcal{F}}}
\newcommand {\ol}{\overline{M-{\mathcal{F}}}}

\newcommand {\stk}{K^{str}_0\left(Q\right)}
\newcommand {\suk}{K^{str}_1\left(Q\right)}

\begin{document}
\title{Generalizations of Agol's inequality and nonexistence of tight
laminations}
\author{Thilo Kuessner}
\date{}
\maketitle
%\tableofcontents

%%%% **** The text of the paper starts here **** %%%%

\section{Results}

Agol's inequality (\cite{ag}, Theorem 2.1.) is the following:\\
\\
{\bf Agol's inequality:} {\em If $M$ is a hyperbolic 
3-manifold containing an incompressible, properly embedded surface $F$, 
then $$Vol\left(M\right)\ge -2 V_3 \chi\left(Guts\left(\overline{M - F}\right)\right),$$
where $V_3$ is the volume of a regular ideal tetrahedron in hyperbolic 3-space.}\\

In \cite{ag2}, this inequality has been improved to 
$$Vol\left(M\right)
\ge Vol\left(Guts\left(\overline{M-F}\right)\right)
\ge -V_{oct}\chi\left(Guts\left(\overline{M-F}\right)\right),$$
where $V_{oct}$ is the volume of a regular ideal octahedron in hyperbolic 3-space.\\

In this paper we will, building on ideas from \cite{ag}, prove a general inequality for the (transversal)
Gromov norm $\parallel M\parallel_{\mathcal{F}}$ and the normal Gromov norm
$\parallel M\parallel^{norm}_{\mathcal{F}}$ of laminations.\\

To state the result in its general form we first need two definitions.\\
\\
{\bf Definition (Pared acylindrical)}: {\em Let $Q$ be a manifold with 
a given decomposition $$\partial Q=\partial_0Q\cup\partial_1Q.$$
The pair $\left(Q,\partial_1Q\right)$ is called a pared
acylindrical manifold, if any continuous mapping of pairs
$f:\left({\bf S}^1\times \left[0,1\right],{\bf S}^1\times\left\{0,1\right\}\right)
\rightarrow \left(Q,\partial_1Q\right)$, which is $\pi_1$-injective as a map of pairs, must be homotopic, as a map of pairs 
$$\left({\bf S}^1\times \left[0,1\right],{\bf S}^1\times\left\{0,1\right\}\right)\rightarrow \left(Q,\partial_1Q\right),$$
into $\partial Q.$}\\
\\
{\bf Definition (Essential Decomposition)}: {\em 
Let $\left(N,\partial N\right)$ be a pair of topological spaces such that $N=Q\cup R$ for two subspaces
$Q,R$. Let $$\partial_0Q=Q\cap R, \partial_1 Q=Q\cap\partial N,\partial_1 R=
R\cap\partial N, \partial Q=\partial_0Q\cup\partial_1Q, \partial R=\partial_0Q \cup\partial_1R.$$
We say that the
decomposition $N=Q\cup R$ is an essential decomposition of $\left(N,\partial N\right)$
if the inclusions 
$$\partial_1Q\rightarrow Q\rightarrow N, \partial_1R\rightarrow R\rightarrow N,
\partial N\rightarrow N,
\partial_0Q\rightarrow Q,
\partial_0Q\rightarrow R$$ 
are each $\pi_1$-injective (for each path-component).}\\
\\
\begin{thm}\label{Thm1}  
Let $M$ be a compact, orientable, connected n-manifold and $\mathcal{F}$ a lamination (of codimension one)
of
$M$. 

Assume that $N:=\overline{M - {\mathcal{F}}}$
has a decomposition $N=Q\cup R$ into orientable n-manifolds (with boundary) $Q,R$ such that the 
%$\partial Q=\partial_0Q\cup\partial_1Q$ with
following assumptions are satisfied for $\partial_0Q=Q\cap R,
\partial_1Q=
Q\cap\partial N, \partial_1R=R\cap \partial N$:\\
i) each path-component of
$\partial_0Q$ has
amenable fundamental group, \\
ii) $\left(Q,\partial_1Q\right)$ is pared acylindrical, $\partial_1Q$ is acylindrical\\
%$\partial_1Q\rightarrow Q$ is $\pi_1$-injective,\\
iii) $Q, \partial N,\partial_1Q,\partial_1R,\partial_0Q$ are aspherical,\\
iv) the decomposition $N=Q\cup R$
is an essential decomposition of $\left(N,\partial N\right)$.

Then
$$\parallel M,\partial M\parallel^{norm}_{\mathcal{F}}
\ge\frac{1}{n+1}\parallel \partial Q\parallel.$$\end{thm}
%$$\parallel M
%\parallel_{\mathcal{F}}\ge
%\frac{1}{\left[\frac{n}{2}\right]+1}\parallel \partial Q\parallel.$$\end{thm}

In the case of 3-manifolds $M$ carrying an essential lamination $\mathcal{F}$, considering $Q=Guts\left(\overline{M - {\mathcal
{F}}}\right)$ yields then as a special case:\\

\begin{thm}\label{Thm2} Let $M$ be a compact 3-manifold with (possibly empty) boundary consisting of incompressible tori, and let
$\mathcal{F}$ be an essential lamination of 
$M$. Then $$\parallel M,\partial M\parallel_{\mathcal{F}}^{norm}\ge
-\chi\left(Guts\left(\overline{M - {\mathcal{F}}}\right)\right).$$
%, \parallel M\parallel_{\mathcal{F}} \ge -2\chi\left(Guts\left(\overline{M - {\mathcal{F}}}\right)\right).$$
More generally, if $P$ is a polyhedron with $f$ faces, then
 $$ \parallel M,\partial M\parallel_{{\mathcal{F}},P}^{norm}\ge
-\frac{2}{f-2}\chi\left(Guts\left(\overline{M - {\mathcal
{F}}}\right)\right).$$\end{thm}

\noindent The following corollary applies, for example, to all hyperbolic manifolds $M$
obtained by Dehn-filling the complement of the figure-eight knot in ${\bf S}^3$. (It is known that each of these $M$ contains tight laminations. By the following corollary, all these tight laminations have empty guts.)\\

{\bf Corollary 4}: {\em If $M$ is a finite-volume hyperbolic 3-manifold with 
$Vol\left(M\right)< 2V_3 =2.02...$, then $M$ carries no essential lamination $\mathcal{F}$
with $\parallel M\parallel_{{\mathcal{F}},P}^{norm}=\parallel M\parallel_P$ 
for all polyhedra P, and nonempty guts. In particular, there is no tight essential
lamination with nonempty guts.}\\

It was observed by Calegari-Dunfield in \cite{cd} that a 
generalization of Agol's inequality to the case of tight laminations, 
together with the results in \cite{cd} about tight laminations with empty guts, would imply the following corollary.\\

{\bf Corollary 5} (\cite{cd}, Conjecture 9.7.): {\em The Weeks manifold admits no tight 
lamination $\mathcal{F}$.}\\

Putting this together with the main result of a recent paper by Tao Li (\cite{li}), one can even improve this result as follows.\\

{\bf Corollary 6}: {\em The Weeks manifold admits no transversely orientable essential lamination.}\\

Finally, we also have an application of  \hyperref[Thm1]{Theorem \ref*{Thm1}} to higher-dimensional manifolds.\\

{\bf Corollary 7}:
{\em  Let $M$ be a compact Riemannian $n$-manifold of negative sectional curvature and finite volume. Let $F\subset M$ be a geodesic $n-1$-dimensional
hypersurface of finite volume. Then $\parallel
F\parallel\le \frac{n+1}{2}\parallel M\parallel$.}\\

\noindent
The basic idea of \hyperref[Thm1]{Theorem \ref*{Thm1}}, say for simplicity in the special situation of \hyperref[Cor7]{Corollary \ref*{Cor7}}, is the following: a simplex which contributes to a normalized fundamental
cycle of $M$ should intersect $\partial Q=2F$ in at most $n+1$ codimension one simplices. This is of course not true in general: simplices can wrap around $M$ many times and
intersect $F$ arbitrarily often, and even a homotopy rel.\ vertices will not change this. As an obvious examle, look at the following situation: let $\gamma$ be a closed geodesic transverse to $F$,
and for some large $N$ let $\sigma$ be a straight simplex contained in a small neighborhood of $\gamma^N$. Then $\sigma$ intersects $F$ $N$ times and, since $\sigma$ is already straight, this number of intersections can of 
course not be reduced by straightening. This shows that some more involved straightening must take place, and that
the acylindricity of $F$ is an essential condition. The way to use acylindricity will be 
to find a normalization such that many subsets of simplices are mapped to cylinders, which degenerate and thus can be removed without changing the homology class.

We remark that many technical points, in particular the use of multicomplexes, can be omitted if (in the setting of \hyperref[Thm2]{Theorem \ref*{Thm2}}) one does not consider incompressible surfaces or essential laminations,
but just geodesic surfaces in hyperbolic manifolds. In this case, all essential 
parts of the proof of \hyperref[Thm1]{Theorem \ref*{Thm1}} enter without 
the notational complications caused by the use of multicomplexes. 
Therefore we have given a fairly detailed outline of the proof for this special case in 
Section \ref{sec:motivating}.
This should help to motivate the general proof in Section \ref{sec:proof}. (We mention that \hyperref[Thm1]{Theorem \ref*{Thm1}} is not true without assuming amenability of $\pi_1\partial_0Q$. This indicates that the proof of multicomplexes in the proof of \hyperref[Thm1]{Theorem \ref*{Thm1}} seems unavoidable.)

Acknowledgements: It is probably obvious that this paper is strongly influenced by Agol's preprint
\cite{ag}. 
Moreover, 
the argument that a generalization of Agol's inequality would imply 
\hyperref[Cor5]{Corollary \ref*{Cor5}} is due to \cite{cd}.
\section{Preliminaries}
\subsection{Laminations}
Let $M$ be an n-manifold, possibly with boundary. In this paper 
all manifolds will be smooth and orientable. (Hence they are triangulable by Whitehead's theorem and
possess a locally finite fundamental class.)
A (codimension 1)
lamination ${\mathcal{F}}$ of $M$
is a foliation of a closed subset $\F$ of $M$, i.e., a decomposition
of a closed subset $\F\subset M$ into immersed codimension 1 submanifolds (leaves) 
so that $M$ is covered by charts $\phi_j:{\bf R}^{n-1}\times{\bf R}\rightarrow 
M$, the intersection of any leaf with the image of any chart $\phi_j$ 
being a union of plaques of the form $\phi_j\left({\bf R}^{n-1}\times\left\{*\right\}\right)$. (We
will denote by $\F$ both the lamination and the laminated subset of $M$, i.e.\ the union of leaves.)
If $M$ has boundary, we will always assume without further mentioning that 
$\mathcal{F}$ is {\em either transverse to $\partial M$} (that is, every leaf is transverse to $\mathcal{F}$) 
{\em or tangential to $\partial M$} (that is, $\partial M$ is a leaf of $\mathcal{F}$). If neither of these two conditions were true, then the transverse and normal Gromov norm would be infinite, therefore all lower bounds will be trivially true.

To construct the leaf space $T$ of $\F$,
one considers the pull-back lamination $\widetilde{\F}$
on the universal covering $\widetilde{M}$.
The space of leaves $T$ is defined as
the quotient of $\widetilde{M}$ under the following equivalence
relation $\sim$. Two points $x,y\in \widetilde{M}$ are equivalent if either they belong to 
the same leaf of $\widetilde{\F}$, or they belong to the same
connected component of the metric completion $\overline{\widetilde{M} - \widetilde{\F}}$ (for the path metric inherited by $\widetilde{M} - \widetilde{\F}$ from an arbitrary Riemannian metric on $\widetilde{M}$).\\
\\
{\em Laminations of 3-manifolds.}
A {\bf lamination} $\F$
of a {\bf 3-manifold} $M$
is called {\bf essential}
if no leaf is a sphere or a torus bounding a solid torus, $\overline{M - \F}$
is irreducible, and $\partial\left(\overline{M - \F}\right)$ is incompressible
and end-incompressible in $\overline{M - \F}$, where again the metric completion
$\overline{M- \F}$ of $M-\F$ is taken w.r.t.\ the path metric inherited from any Riemannian metric on $M$,
see \cite{goe}, ch.1. (Note that $\overline{M- \F}$ is immersed in $M$, the leaves of
$\mathcal{F}$ in the image of the immersion are called boundary leaves.)

Examples
of essential laminations are taut foliations or compact, incompressible, boundary-incompressible surfaces
in compact 3-manifolds.
(We always consider laminations without isolated leaves. If a lamination has isolated 
leaves, then it can be converted into a lamination without isolated leaves by replacing each two-sided isolated leaf $S_i$ with the
trivially foliated product $S_i\times\left[0,1\right]$, resp.\ each one-sided 
isolated leaf with the canonically foliated normal $I$-bundle, without changing the topological type of $M$.)

If $\mathcal{\F}$ is an essential lamination, then the
leaf space $T$ is an order tree, with segments 
corresponding to directed, transverse, efficient arcs. (An order tree $T$ is a set $T$ with a collection of linearly ordered subsets, called segments, such that the axioms of \cite{goe}, Def.\ 6.9., are satisfied.)
Moreover, $T$ is an ${\bf R}$-order tree, that is, it is a countable union of segments
and each segment is order isomorphic to a closed interval in ${\bf R}$. $T$ can be topologized by the order topology on segments (and declaring that a set is closed if the intersection with each segment is closed). 
For this topology, $\pi_0T$ and $\pi_1T$ are trivial (see, for example, \cite{rss}, Chapter 5, and its references).

The order tree $T$ comes 
with a fixed-point 
free action of $\pi_1M$. Fenley (\cite{fe}) 
has exhibited hyperbolic 3-manifolds whose fundamental groups do not admit any fixed-point free
action on ${\bf R}$-order
trees. Thus there are hyperbolic 3-manifolds not carrying any essential lamination.  

If $M$ is hyperbolic and $\mathcal{F}$ an essential lamination, then $\overline{M - {\mathcal{F}}}$
has a characteristic submanifold which is the maximal submanifold 
that can be decomposed into $I$-bundles and solid tori, respecting boundary patterns (see \cite{js}, \cite{joh} for precise definitions). The complement of this characteristic submanifold
is denoted by $Guts\left({\mathcal{F}}\right)$. It admits a hyperbolic 
metric with geodesic boundary 
and cusps. (Be aware that some authors, like \cite{cd}, include the solid tori into the guts.)
If ${\mathcal{F}}=F$ is a properly embedded, incompressible, boundary-incompressible surface, then
{\bf Agol's inequality} states that $Vol\left(M\right)\ge -2 V_3\chi\left(Guts\left(F\right)\right)$. This implies, for example, that a hyperbolic manifold of volume $< 2V_3$ can
not contain any geodesic surface of finite area. Recently, this inequality has been improved to $Vol\left(M\right)
\ge Vol\left(Guts\left(F\right)\right)
\ge -V_{oct}\chi\left(Guts\left(F\right)\right)$ in \cite{ag2}, using estimates coming from
Perelman's work on the Ricci flow.\\

Assume that $\mathcal{F}$ is a codimension one lamination
of an n-manifold $M$ such that its leaf space $T$ is an ${\bf R}$-order tree. (For example this is the case if $n=3$ and $\mathcal{F}$ is essential.)
An essential lamination is called {\bf tight} if $T$
is Hausdorff. It is called unbranched if  
$T$
is homeomorphic to $
{\bf R}$. It is said to have two-sided branching (\cite{cal}, Definition 2.5.2) 
if there are leaves $\lambda,\lambda_1,\lambda_2,\mu,\mu_1,\mu_2$ such that the corresponding points in the
$T$ satisfy $\lambda<\lambda_1,\lambda<\lambda_2,\mu>\mu_1,\mu>\mu_2$ but $\lambda_1,\lambda_2$ are incomparable and $\mu_1,\mu_2$ are incomparable.
%if it is not unbranched but there is some point $x\in T$ such that one connected component 
%of $T-\left\{x\right\}$ is order-isomorphic to a connected subset of ${\bf R}$. 
It is said to have one-sided branching if it is neither unbranched nor has two-sided branching.\\

If $M$
is a hyperbolic 3-manifold and carries a tight lamination with empty guts, then
Calegari and Dunfield have shown (\cite{cd}, Theorem 3.2.) that $\pi_1M$ acts effectively on the circle, i.e., there is an injective homomorphism
$\pi_1M\rightarrow Homeo\left({\bf S}^1\right)$. This
implies that the Weeks manifold (the closed hyperbolic manifold of smallest
volume) can not carry a tight lamination with empty guts (\cite{cd}, Corollary 9.4.). The aim of this paper
is to find obstructions to the existence of laminations with nonempty guts.

\subsection{Simplicial volume and refinements}

Let $M$ be a compact, orientable, connected n-manifold, possibly with boundary. Its top integer (singular) homology
group $H_n\left(M,\partial M;{\bf Z}\right)$ is cyclic. The image of a generator under the change-of-coefficients homomorphism $H_n\left(M,\partial M;{\bf Z}\right)\rightarrow
H_n\left(M,\partial M;{\bf R}\right)$ is called a fundamental class and is 
denoted $\left[M,\partial M\right]$. If $M$ is not connected, we define $\left[M,\partial M\right]$ to be the formal sum of the fundamental classes of its connected components.

The simplicial volume $\parallel M,\partial M\parallel$ is defined as $\parallel M,
\partial M\parallel=inf\left\{\sum_{i=1}^r\mid a_i\mid\right\}$ 
where the infimum is taken over all singular chains $\sum_{i=1}^r a_i\sigma_i$
(with real coefficients) representing the fundamental 
class in $H_n\left(M,\partial M;{\bf R}\right)$. 

If $M - \partial M$ carries a complete hyperbolic metric of finite volume $Vol\left(M\right)$,
then $\parallel M,\partial M\parallel=\frac{1}{V_n}Vol\left(M\right)$ with $V_n=
sup\left\{Vol\left(\Delta\right):\Delta\subset{\bf H}^n\ \mbox{geodesic simplex}\right\}$ (see \cite{gro},\cite{thu},\cite{bp}, \cite{fr}).

More generally, let $P$ be any polyhedron. Then the invariant $\parallel M,\partial M\parallel_P$
is defined in \cite{ag} as follows: denoting by $C_*\left(M,\partial M;P;{\bf R}\right)$ the complex of $P$-chains with real coefficients, and by $H_*\left(
M,\partial M;P;{\bf R}\right)$ its homology, there is a canonical chain homomorphism $\psi:C_*\left(
M,\partial M;P;{\bf R}\right)\rightarrow C_*\left(M,\partial M;{\bf R}\right)$, given by some 
triangulations of $P$ which is to be chosen such that all possible cancellations of boundary faces are preserved. 
$\parallel M,\partial M\parallel_P$ is defined as the infimum of $\sum_{i=1}^r\mid a_i\mid$
over all $P$-chains $\sum_{i=1}^r a_iP_i$ such that $\psi\left(\sum_{i=1}^r a_iP_i\right)$ represents
the fundamental class $\left[M,\partial M\right]$. Let $V_P:=sup\left\{Vol\left(\Delta\right)\right\}$,
where the supremum is taken over all straight $P$-polyhedra $\Delta\subset{\bf H}^3$. 
%The following 
\hyperref[Prop1]{Proposition \ref*{Prop1}} is Lemma 4.1. in \cite{ag}. (The proof in \cite{ag} is quite short, and it does not give details for the cusped case. However, the
proof in the cusped case can be completed using the arguments in sections 5 and 6
of Francaviglia's paper \cite{fr}.)

\begin{pro}\label{Prop1} If $M - \partial M$ admits a hyperbolic metric of finite volume $Vol\left(M\right)$, then
$$\parallel M,\partial M\parallel_P=\frac{1}{V_P}Vol\left(M\right).$$\end{pro}

Let $M$ be a manifold and $\mathcal{F}$ a {\bf codimension one lamination} of $M$.
Let $\Delta^n$ be the standard simplex in ${\bf R}^{n+1}$, and $\sigma:\Delta^n\rightarrow M$ some continuous
singular simplex.
The lamination $\mathcal{F}$ induces an equivalence relation on $\Delta^n$ by:
$x\sim y\Longleftrightarrow \sigma\left(x\right)$ and $\sigma\left(y\right)$ belong to the same connected component
of $L\cap\sigma\left(\Delta^n\right)$ for some leaf $L$ of $\mathcal{F}$.
We say that a singular simplex $\sigma:\Delta^n\rightarrow M$ is laminated if
the equivalence relation $\sim$
is induced by a lamination ${\mathcal{F}}\mid_{\sigma}$
of $\Delta^n$.
We call a lamination $\mathcal{F}$ of $\Delta^n$ affine if there is an affine mapping $f:\Delta^n\rightarrow{\bf R}                                                                                                                  $ such that $x,y\in\Delta^n$ belong to the same leaf if and only if $f\left(x\right)=f\left(y\right)$. We 
say that a lamination $\mathcal{G}$ of $\Delta^n$ is conjugate to an affine lamination if there is a 
simplicial homeomorphism $H:\Delta^n                   
\rightarrow\Delta^n$ such that $H^*\mathcal{G}$ is an affine lamination.\\   
We say that a singular $n$-simplex $\sigma:\Delta^n\rightarrow M$, $n\ge 2$, is {\bf transverse} to 
$\mathcal{F}$ if it is laminated and  it is 
either contained in a leaf,
or ${\mathcal{F}}\mid_\sigma$ is conjugate to 
an affine lamination $\mathcal{G}$ of $\Delta^n$.\\
For $n=1$, we say that a singular 1-simplex $\sigma:\Delta^1\rightarrow M$ is transverse to $\mathcal{F}$ if it is either contained in a leaf,                                                                                        or for each lamination chart $\phi:U\rightarrow {\bf R}^{m-1}\times{\bf R}^1$ (with m-th coordinate map           $\phi_m:U\rightarrow{\bf R}^1$) one has that $\phi_m\circ\sigma\mid_{\sigma^{-1}\left(U\right)}:                  \sigma^{-1}\left(U\right)\rightarrow{\bf R}^1$ is locally surjective at                                           all points of $int\left(\Delta^1\right)$, i.e.\ for all                                                            $p\in int\left(\Delta^1\right)\cap \sigma^{-1}\left(U\right)$, the                                                 image of $\phi_m\circ\sigma\mid_{\sigma^{-1}\left(U\right)}$ contains a neighborhood of $\phi_m\circ\sigma\left(p\right)$.\\ 
We say that
the simplex $\sigma:\Delta^n\rightarrow M$ is {\bf normal} to ${\mathcal{F} }$ if, for each leaf $F$,
$\sigma^{-1}\left(F\right)$ consists of normal disks, i.e.\ disks meeting each edge of $\Delta^n$ at most once.
(If $F=\partial M$ is a leaf of ${\mathcal{F}}$ we also allow that $\sigma^{-1}\left(F\right)$
can be a face of $\Delta^n$). In particular, any transverse simplex is normal.\\
In the special case of 
{\bf foliations} $\mathcal{F}$
one has that
the transversality of a singular simplex $\sigma$ is implied by (hence equivalent to) the normality of $\sigma$,
as can be shown along the lines of \cite{k3}, section 1.3.\\
More generally, let $P$ be any polyhedron. Then we say that a singular polyhedron $\sigma:P\rightarrow M$ is normal to $\mathcal{F}$ if, 
for each leaf $F$,
$\sigma^{-1}\left(F\right)$ consists of normal disks, i.e.\ disks meeting each edge of $P$ at most once
(or being equal to a face of $P$, if $F$ is a boundary leaf).
%each edge is either contained in a leaf or
%intersects each leaf of $\mathcal{F}$ at most once.

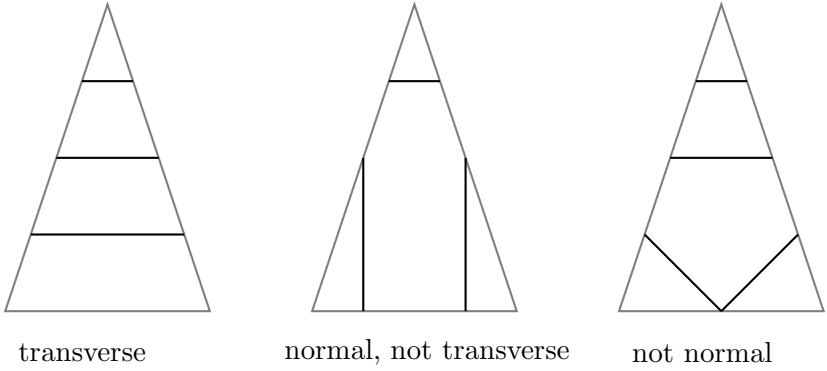
\psset{unit=0.1\hsize}
$$\pspicture(0,-1)(10,3)
\pspolygon[linecolor=gray](0,0)(2,0)(1,3)(0,0)
\psline(0.5,1.5)(1.5,1.5)
\psline(0.25,0.75)(1.75,0.75)
\psline(0.75,2.25)(1.25,2.25)
\uput[0](0,-0.4){transverse}

\pspolygon[linecolor=gray](3,0)(5,0)(4,3)(3,0)
\psline(3.5,0)(3.5,1.5)
\psline(3.75,2.25)(4.25,2.25)
\psline(4.5,0)(4.5,1.5)
\uput[0](2.6,-0.4){normal, not transverse}

\pspolygon[linecolor=gray](6,0)(8,0)(7,3)(6,0)
\psline(6.5,1.5)(7.5,1.5)
\psline(6.75,2.25)(7.25,2.25)
\psline(6.25,0.75)(7,0)
\psline(7,0)(7.75,0.75)
\uput[0](6,-0.4){not normal}

\endpspicture$$

\begin{df} Let $M$ be a compact, oriented, connected n-manifold, possibly with boundary, and ${\mathcal{F}}$ a foliation or lamination on $M$. Let $\Delta^n$ be the standard simplex and $P$ any polyhedron. Then $$\parallel M,\partial M\parallel_{{\mathcal{F}}}:=inf\left\{\sum_{i=1}^r\mid a_i\mid: \psi\left(\sum_{i=1}^r a_i\sigma_i\right) \mbox{ represents }\left[M,\partial M\right], \sigma_i:\Delta^n\rightarrow M \mbox{ transverse to }{\mathcal{F}}\right\}$$
and 
$$\parallel M,\partial M\parallel^{norm}_{{\mathcal{F}},P}:=inf\left\{
\sum_{i=1}^r\mid a_i\mid: \psi\left(\sum_{i=1}^r a_i\sigma_i\right) \mbox{ represents
 }\left[M,\partial M\right], \sigma_i:P\rightarrow M \mbox{ normal to }{
\mathcal{F}}\right\}.$$
In particular, we define $\parallel M,\partial M\parallel_{\mathcal{F}}^{norm}
=\parallel M,\partial M
\parallel_{{\mathcal{F}},\Delta^n}^{norm}$.\end{df}
All norms are finite, under the assumption that $\mathcal{F}$ is transverse
or tangential to $\partial M$.\\
There is an obvious inequality $$\parallel M,\partial M\parallel\le
\parallel M,\partial M\parallel^{norm}_{{\mathcal{F}}}\le
\parallel M,\partial M\parallel_{{\mathcal{F}}}.$$
In the case of foliations, equality $\parallel M,\partial M\parallel^{norm}_{{\mathcal{F}}}=\parallel M,\partial M\parallel_{{\mathcal{F}}}$ holds.\\
(We remark that all definitions extend in an obvious way to disconnected manifolds by summing over the connected components.)\\

\hyperref[Prop2]{Proposition \ref*{Prop2}} and \hyperref[Lemma1]{Lemma \ref*{Lemma1}} are
a straightforward generalisation of \cite{cal}, Theorem 2.5.9, and of arguments in \cite{ag}.
\begin{pro}\label{Prop2} Let $M$ be a compact, oriented 3-manifold.\\
a) If ${\mathcal{F}}$ is an
essential lamination which is either unbranched or has 
one-sided branching such that the induced lamination of $\partial M$ is unbranched, then
$$\parallel M,\partial M\parallel^{norm}_{
{\mathcal{F}},P}=\parallel M,\partial M\parallel_P$$ for each polyhedron $P$.\\
b) If ${\mathcal{F}}$ is a tight essential lamination, then
$$\parallel M,\partial M\parallel^{norm}_{
{\mathcal{F}},P}=\parallel M,\partial M\parallel_P$$ for each polyhedron $P$.\end{pro}
\begin{proof} Since ${\mathcal{F}}$ is an essential lamination, we know from \cite{goe}, Theorem 6.1., that the leaves are $\pi_1$-injective, 
the universal covering $\widetilde{M}$ is homeomorphic
to ${\bf R}^3$ and that the leaves of the pull-back lamination are planes, in particular aspherical. 
%Moreover, it follows with $\widetilde{M}={\bf R}^3$ and $\widetilde{F}\simeq{\bf R}^2$ from
%Alexander duality that
%$\widetilde{M}-\widetilde{F}$ has
%two connected components for each leaf $F$. 
Therefore Proposition 2 is a special case of \hyperref[Lemma1]{Lemma \ref*{Lemma1}}.\end{proof}

\begin{lem}\label{Lemma1} Let $M$ be a compact, oriented, 
aspherical manifold, and $\F$ a lamination of codimension 
one. \\
Assume that the leaves are $\pi_1$-injective and aspherical, and that the leaf space $T$ is an ${\bf R}$-order tree.\\
a) If the leaf space $T$ is either ${\bf R}$ or branches 
in only one direction, such that the induced lamination of $\partial M$ has leaf space ${\bf R}$, then $$\parallel M,\partial M
\parallel^{norm}_{{\mathcal{F}},P}=\parallel M,\partial M\parallel_P$$ for each polyhedron $P$.\\
b) If the leaf space is a Hausdorff tree, 
%and $\widetilde{M}-\widetilde{F}$ has
%two connected components for each leaf $F$, 
then $$\parallel M,\partial M
\parallel^{norm}_{{\mathcal{F}},P}=\parallel M,\partial M\parallel_P$$ for each polyhedron $P$.\end{lem}

\begin{proof}
To prove the wanted equalities, it suffices in each case to show that any (relative)
cycle can be homotoped to a cycle consisting of normal polyhedra.\\
We denote by $\widetilde{\mathcal{F}}$ the pull-back lamination of $\widetilde{M}$ and $p:\widetilde{M}\rightarrow T=\widetilde{M}/\widetilde{\mathcal{F}}$ the projection to the leaf space.\\
a)
First we consider the case that $P$=simplex (\cite{cal}, Section 4.1) and $\mathcal{F}$ unbranched. 
For this case, we can repeat the argument in \cite{cal}, Lemma 2.2.8. Namely, let us be given a 
(relative) cycle $\sum_{i=1}^r a_i\sigma_i$, lift it to a $\pi_1M$-equivariant (relative) cycle on $\widetilde{M}$ and then perform an (equivariant) straightening, by induction on the dimension of subsimplices of the lifts $\widetilde{\sigma_i}$ as follows: for each edge $\tilde{e}$
of any 
lift $\widetilde{\sigma_i}$, its projection $p\left(\tilde{e}\right)$
to the leaf space $T$ is homotopic to a unique straight
arc $str\left(p\left(\tilde{e}\right)\right)$
in $T\simeq{\bf R}$. It is easy to see (covering the arc by foliation charts and then extending the
lifted arc stepwise) that $str\left(p\left(\tilde{e}\right)\right)$ 
can be lifted to an arc $str\left(\tilde{e}\right)$ with 
the same endpoints as $\tilde{e}$, and that 
the homotopy between $str\left(p\left(\tilde{e}\right)\right)$
and $p\left(\tilde{e}\right)$ can be lifted to a homotopy between $str\left(\tilde{e}\right)$ and $\tilde{e}$.
$str\left(\tilde{e}\right)$ is transverse to $\mathcal{F}$, because its projection is a straight arc in $T$. 
These homotopies of edges can be extended to a homotopy of the whole (relative) cycle. Thus we have straightened the 1-skeleton of the given (relative) cycle.

Now let us be given a 2-simplex $\tilde{f}:
\Delta^2\rightarrow\widetilde{M}$ with transverse edges. There is 
an obvious straightening $str\left(p\left(\tilde{f}\right)\right)$
of $p\left(\tilde{f}\right):\Delta^2\rightarrow T$ as 
follows: if, for $t\in T$, $\left(p \tilde{f}\right)^{-1}\left(t\right)$ 
has two preimages $x_1,x_2$ on edges of $\Delta^2$ (which are necessarily unique), then 
$str\left(p\left(\tilde{f}\right)\right)$ maps the line which connects $x_1$ and $x_2$ in
$\Delta^2$ constantly to $t$. It is clear that this defines a continuous map $str\left(p\left(\tilde{f}\right)\right):\Delta^2\rightarrow T$.

Since leaves $\widetilde{F}$
of $\widetilde{\mathcal{F}}$ are connected ($\pi_0\widetilde{F}=0$), $str\left(p\left(\tilde{f}\right)\right)$ can be lifted to a map $str\left(\tilde{f}\right):\Delta^2\rightarrow\widetilde{M}$ with $p\left(str\left(\tilde{f}\right)\right)=
str\left(p\left(\tilde{f}\right)\right)$. $str\left(\tilde{f}\right)$ is transverse to $\mathcal{F}$, because its projection is a straight simplex in $T$.

There is an obvious homotopy between $p\left(\tilde{f}\right)$ and
$str\left(p\left(\tilde{f}\right)\right)$. For each $t\in T$, the restriction of the homotopy to $\left(p \tilde{f}\right)^{-1}\left(t\right)$
can be lifted to a homotopy in $\widetilde{M}$, because $\pi_1\widetilde{M}=0$.
Since $\pi_2\widetilde{M}=0$, these homotopies for various $t\in T$ fit together continuously to give a homotopy between $\tilde{f}$ and $str\left(\tilde{f}\right)$.

These homotopies of 2-simplices leave the (already transverse) boundaries pointwise fixed, thus they
can be extended to a homotopy of the whole (relative) cycle. 
Hence we have straightened the 2-skeleton of the given (relative) cycle.

Assume that we have already straightened the $k$-skeleton, for some $k\in{\bf N}$.
The analogous procedure, using $\pi_{k-1}\widetilde{F}=0$ for all leaves, and $\pi_k\widetilde{M}=0, \pi_{k+1}\widetilde{M}=0$, allows to straighten the $\left(k+1\right)$-skeleton of the (relative) cycle. This finishes the proof in
the case that ${\mathcal{F}}$ is unbranched.

The generalization to the case that $\mathcal{F}$ has one-sided branching such that the induced lamination of $\partial M$ is unbranched works as in \cite{cal}, Theorem 2.6.6.

We remark that in the case $P=simplex$ we get not only a normal cycle, but even a transverse cycle.

Now we consider the case of arbitrary polyhedra $P$. Let $\sum_{i=1}^r a_i\sigma_i$ be a P-cycle. It can be 
subtriangulated to a simplicial cycle $\sum_{i=1}^r a_i\sum_{j=1}^s\tau_{i,j}$. Again the argument in \cite
{cal}, Lemma 2.2.8 (resp.\ its version for manifolds with boundary), shows that this simplicial cycle can be homotoped such that each $\tau_{i,j}$ is transverse (and such that boundary cancellations are preserved).
But transversality of each $\tau_{i,j}$ implies by definition that $\sigma_i=\sum_{j=1}^s\tau_{i,j}$ is normal (though in general not transverse) to $\mathcal{F}$.\\
b) By assumption $\widetilde{M}/\widetilde{\mathcal{F}}$ is a Hausdorff tree. We observe that its branching points are 
the projections of complementary regions. Indeed, let $F$ 
be a leaf of $\F$, then $\widetilde{F}$ is a submanifold of 
the contractible manifold $\widetilde{M}$.
By asphericity and $\pi_1$-injectivity of $F$, $\widetilde{F}$ must be contractible.
By Alexander duality it follows that
$\widetilde{M}-\widetilde{F}$ has 
two connected components. Therefore the complement of the point $p\left({\widetilde{F}}\right)$ in the leaf space
has (at most) two connected components, thus $p\left(\widetilde{F}\right)$
can not be a branch point.\\
Again, to define a straightening of $P$-chains it suffices to
define a canonical straightening of singular polyhedra $P$ such that straightenings of common boundary faces will agree. Let $\tilde{v}_0,\ldots,\tilde{v}_n$ be the vertices of the image of $P$.
For each pair $\left\{\tilde{v}_i,\tilde{v}_j\right\}$ 
there exists at most one edge $\tilde{e}_{ij}$ with
vertices $\tilde{v}_i,\tilde{v}_j$
in the image of $P$. Since the leaf space is a tree, we have a unique 
straight arc $str\left(p\left(\tilde{e}_{ij}\right)\right)$ 
connecting the points $p\left(\tilde{v}_i\right)$ and $p\left(\tilde{v}_j\right)$ 
in the leaf space. As in a), one can lift this straight arc
$str\left(p\left(\tilde{e}_{ij}\right)\right)$ to an arc $str\left(\tilde{e}_{ij}\right)$ 
in $\widetilde{M}$, connecting $\tilde{v}_i$ and $\tilde{v}_j$, which
is transverse to $\mathcal{F}$. 
We define this arc $str\left(\tilde{e}_{ij}\right)$ to be
the straightening of $\tilde{e}_{ij}$.
As in a), we have homotopies of 1-simplices, which extend to a homotopy of the whole (relative) cycle.
Thus we have straightened the 1-skeleton. 

Now let us be given the 3 vertices $\tilde{v}_0,\tilde{v}_1,\tilde{v}_2$ of a 
2-simplex $\tilde{f}$
with straight edges. If the projections $p\left(\tilde{v}_0\right),p\left(\tilde{v}_1\right),p\left(\tilde{v}_2\right)$
belong to a subtree isomorphic to a connected subset of
${\bf R}$, then we can straighten $\tilde{f}$
as in a). If not, we have that the projection of the 1-skeleton of this simplex has exactly one branch point, which corresponds to a complementary region. (The projection may of course meet many branch points of the tree, but the image of the projection, considered as a subtree, can have at most one branch point. In general, a subtree with $n$ vertices can have at most $n-2$ branch points.)
The preimage of the complement of this complementary region
consists of three connected subsets of the 2-simplex ("corners around the vertices").
We can straighten each of these subsets and do not need to care about the complementary region corresponding to the branch point. Thus we have straightened the 2-skeleton.\\
Assume that we have already straightened the $k$-skeleton, for
some $k\in{\bf N}$.
Let us be given the $k+2$ vertices $\tilde{v}_0,\tilde{v}_1,\ldots,\tilde{v}_{k+1}$ of a $\left(k+1\right)$-simplex 
with straight faces. Then we have (at most $k$) branch points in the projection of
the simplex, which correspond to complementary regions. Again we can straighten 
the parts of the simplex which do not belong to these 
complementary regions as in a), since they are projected to linearly ordered subsets of the tree. Thus we have straightened the $\left(k+1\right)$-skeleton.

Since, by the recursive construction,
we have defined straightenings of simplices with common faces by first 
defining (the same) straightenings of their common faces,
the straightening of a 
(relative)
cycle will be again a (relative) cycle, in the same (relative) homology class.
\end{proof}

Remark: For $\parallel M\parallel_{\mathcal{F}}$ instead of
$\parallel M\parallel^{norm}_{\mathcal{F}}$, equality b) is in general wrong, and equality a) is unknown (but presumably wrong). \\
If $\mathcal{F}$ is essential but not tight, one may still try to homotope cycles to be transverse, 
by possibly changing the lamination. In the special case 
that the cycle is coming from a triangulation, this has been done in \cite{bri}
and \cite{gab} by Brittenham resp.\ Gabai. It is not obvious how to generalize their arguments to cycles with overlapping simplices.

\section{Retracting chains to codimension zero submanifolds}

\subsection{Definitions}
 
The results of this section are essentially all due to Gromov, but we follow mainly our exposition in \cite{k2}.
We start with some recollections about multicomplexes (cf.\ \cite{gro}, Section 3, or
\cite{k2}, Section 1).\\
\\
A multicomplex $K$ is a topological space $\mid K\mid$ with a decomposition into simplices,
where each $n$-simplex is attached to the $n-1$-skeleton $K_{n-1}$
by a simplicial homeomorphism $f:\partial \Delta^n\rightarrow K_{n-1}$.
(In particular, each $n$-simplex has $n+1$ distinct vertices.) 

As opposed to simplicial complexes, in a multicomplex
there may be $n$-simplices with
the same $n-1$-skeleton.

We call a multicomplex minimally complete if the following holds:
whenever $\sigma:\Delta^n\rightarrow \mid K\mid$ is a singular n-simplex, such that
$\partial_0\sigma,\ldots,\partial_n\sigma$ are distinct simplices of $K$, then $\sigma$ is homotopic relative $\partial\Delta^n$ to
a {\em unique} simplex in $K$.

We call a minimally complete
multicomplex $K$ aspherical if all simplices $\sigma\not =\tau$ in $K$
satisfy $\sigma_1\not=\tau_1$. That means, simplices are uniquely determined by their 1-skeleton.        \\
Orientations of multicomplexes are defined as
usual in the simplicial theory. For a simplex $\sigma$, $\overline{\sigma}$ will denote the simplex with the
opposite orientation.\\
A submulticomplex $L$ of a multicomplex $K$ consists of a subset of the set of simplices closed under face maps.
$\left(K,L\right)$ is a pair of multicomplexes if $K$ is a multicomplex and $L$ is a submulticomplex of $K$.\\
A group $G$ acts simplicially on a pair of multicomplexes $\left(K,L\right)$ if it acts on the set of simplices of $K$, mapping simplices in $L$ to simplices in $L$, such that the action commutes with all face maps. For $g\in G$ and $\sigma$ a simplex in $K$, we denote by
$g\sigma$ the simplex obtained by this action.\\
%Notational remark: for a connected multicomplex $K$ we write
%$\pi_1K$ for $\pi_1\left(\mid K\mid,p\right)$ with some $p\in\mid K\mid$. The choice of $p$ will not be of importance for our arguments.

\subsection{Construction of $K\left(X\right)$}\label{sec:construction}

We recall the construction from \cite{k2}, section 1.3 (originally due to \cite{gro}, page 45-46).
 
For a topological space $X$, we denote by $S_*\left(X\right)$ the simplicial set of all singular simplices in $X$ and $\mid S_*\left(X\right)\mid$ its geometric realization.

For a topological space $X$, a 
multicomplex $\widehat{K}\left(X\right)\subset \mid S_*\left(X\right)\mid$ is 
constructed as follows. The 0-skeleton $\widehat{K}_0\left(X\right)$ equals $S_0\left(X\right)$. The 1-skeleton $\widehat{K}_1\left(X\right)$
contains one element
%\footnotemark\footnotetext{We mention a trivial observation that will play a role in ?? Let $X$ be a manifold of dimension $\ge 2$, x\i X and $\left\{e_1,\ldots,e_n\right\}$ a finite set of singular 1-simplices with $\partial_1ei_i=x$ for $i=1,\ldots,n$ and $\partial_0e_i\not=\partial_0e_j$ for $i\not=j$. Then, in the construction of $K\left(X\right)$, we can choose the 1-simplices $\overline{e}_i\in K_1\left(X\right)$, which are homotopic rel.\ $\left\{0,1\right\}$ to $e_i$, such that their images in $X$ intersect only in $x$ and are disjoint otherwise. In particular, the union of their images in $X$ is contractible.} 
in each 
homotopy class (rel.\ $\left\{0,1\right\}$) of singular 1-simplices $f:\left[0,1\right]\rightarrow X$ 
with $f\left(0\right)\not=f\left(1\right)$.
For $n\ge 2$,
assuming by recursion
that the n-1-skeleton is defined, the n-skeleton $\widehat{K}_n\left(X\right)$ contains {\em one} 
singular n-simplex in each homotopy class 
(rel.\ boundary) of singular n-simplices $f:\Delta^n\rightarrow X$ with $\partial f\in 
\widehat{K}_{n-1}\left(X\right)$. We can choose simplices in $\widehat{K}\left(X\right)$ 
such that $\sigma\in \widehat{K}\left(X\right)\Leftrightarrow \overline{\sigma}\in \widehat{K}\left(X\right)$, where $\overline{\sigma}$ denotes the simplex with the opposite orientation. We will henceforth assume that $\widehat{K}\left(X\right)$ is constructed according to this condition.

According to \cite{gro}, $\mid \widehat{K}\left(X\right)\mid $ is weakly homotopy equivalent to $X$.

The multicomplex $K\left(X\right)$ is defined as the quotient $$K\left(X\right):=\widehat{K}\left(X\right)/\sim$$
where simplices in $\widehat{K}\left(X\right)$ are identified if and only if they have the same 1-skeleton. Let $p$ be the canonical projection $p:\widehat{K}\left(X\right)\rightarrow K\left(X\right).$

$K\left(X\right)$ is minimally complete and aspherical.

If $X^\prime\subset X$ is a subspace, then 
%we will simplices in $K\left(X\right)$ to have image in $X^\prime$ whenever this is possible.
% is $\pi_1$-injective (for each connected component), then
%$K\left(X\right)$ (resp.\ $\widehat{K}\left(X\right)$) can be constructed
%such that $K\left(X^\prime\right)$ naturally embeds into $K\left(X\right)$. If $X^\prime\subset X$ is 
%not $\pi_1$-injective, 
we have (not necessarily injective) simplicial mappings $\hat{j}:\widehat{K}\left(X^\prime\right)
\rightarrow\widehat{K}\left(X\right)$ and  
$j:K\left(X^\prime\right)\rightarrow K\left(X\right)$. 

If $\pi_1X^\prime\rightarrow\pi_1 X$ is injective (for each path-connected component of $X^\prime$), then
$j$ is injective (\cite{k2}, Section 1.3) 
and we can (and will) consider
$K\left(X^\prime\right)$ as a submulticomplex of
$K\left(X\right)$.
%We will, by abuse of notation, denote the images of $\tilde{j}$ resp.\ $j$ by $\widetilde{K}\left(X^\prime\right)$ resp.\ $K
%\left(X^\prime\right)$. 
(Since simplices in $\widehat{K}\left(X^\prime\right)$ have image in $X^\prime$, 
this means that we assume to have constructed $\widehat{K}\left(X\right)$ such that
simplices in $\widehat{K}\left(X\right)$ have image in $X^\prime$ whenever this is possible.) If moreover $\pi_nX^\prime\rightarrow \pi_nX$ is injective for all $n\ge 2$ (e.g.\ if $X^\prime$ is aspherical), then also 
$\hat{j}$ is injective and $\widehat{K}\left(X^\prime\right)$ can be considered as 
a submulticomplex of $\widehat{K}\left(X\right)$.

%We will denote $C_*^{simp}\left(K\left(X\right), K\left(X^\prime\right)\right)$ (resp.\ $C_*^{simp}\left(\widehat{K}\left(X\right), \widehat{K}
%\left(X^\prime\right)\right)$ the
%relative chain complex of simplices in $K\left(X\right)$ (resp.\ $\widehat{K}\left(X\right)$), 
%modulo simplices in the image of $K\left(X^\prime\right)$ (resp.\ $\widehat{K}\left(X^\prime\right)$).

In particular, if {\bf $X$ and $X^\prime$ are aspherical and $\pi_1X^\prime\rightarrow \pi_1X$ is injective}, 
then there is an inclusion
$$i_*:C_*^{simp}\left(K\left(X\right),K\left(X^\prime\right)\right)=
C_*^{simp}\left(\widehat{K}\left(X\right),\widehat{K}\left(X^\prime\right)\right)
\rightarrow
C_*^{sing}\left(X,X^\prime\right)$$
into the relative singular chain complex of $\left(X,X^\prime\right)$.\\
% (This is true regardless whether $X^\prime$ is aspherical or not, since we only consider images of simplices in $X$.)\\ 

%In section 5 we will consider manifolds $Q$ with a decomposition $\partial
%Q=\partial_0Q\cup \partial_1Q$. In this case we assume
%to have chosen the simplices in $K\left(Q\right)$ such that, whenever a simplex $\sigma\in K\left(Q\right)$ is homotopic rel.\ boundary 
%to a singular simplex with image in $\partial Q$ resp.\ $\partial_0Q$ resp.\ $\partial_1Q$, then $\sigma$ has image
%in $\partial Q$ resp.\ $\partial_0Q$ resp.\ $\partial_1Q$. 
%Such a choice is indeed possible. Namely, 
%if a simplex $\sigma\in S_n\left(Q\right)$ is homotopic rel.\ boundary to simplices $\sigma_0\in S_n\left(\partial_0Q\right)$
%and $\sigma_1\in S_n\left(\partial_1Q\right)$,
%then we get a homotopy rel.\ boundary between $\sigma_0$ and $\sigma_1$, giving a map
%${\bf D}^{n+1}\rightarrow Q$ such that
%$\partial {\bf D}^{n+1}$
%is divided into two halves, one mapped by $\sigma_0$ to $\partial_0Q$, the other mapped by $\sigma_1$
%$to $\partial_1Q$. Hence there must be an n-simplex in the given homotopy class
%rel.\ boundary mapped to $\partial_0Q\cap\partial_1Q$. This shows that we can
%choose $\sigma_0=\sigma_1\in S_n\left(\partial_0Q\right)\cap S_n\left(
%\partial_1Q\right)$ in the homotopy class rel.\ boundary of $\sigma$.\\

\noindent
{\bf Infinite and locally finite chains.} In this paper we will also work with infinite chains, and in particular
with locally finite chains on noncompact manifolds, as introduced in \cite{gro}, section 0.2. \\
For a topological space $X$, a formal sum $\sum_{i\in I}a_i\sigma_i$ of singular k-simplices with real coefficients
(with a possibly infinite index set $I$, and the convention $a_i\not=0$ for $i\in I$)
is an infinite singular k-chain. It is
said to be a locally finite chain if each point of $X$ is contained in 
the image of at most finitely many $\sigma_i$. Infinite resp.\ locally finite k-chains form 
real vector spaces $C_k^{inf}\left(X\right)$ resp.\ $C_k^{lf}\left(X\right)$.
The boundary operator maps locally finite k-chains
to locally finite k-1-chains, hence, for a pair of spaces $\left(X,X^\prime
\right)$ the homology $H_*^{lf}\left(X,X^\prime\right)$ of the complex of locally finite chains can be defined.\\
For a noncompact, orientable $n$-manifold $X$ with (possibly noncompact)
boundary $\partial X$, one has a fundamental class $\left[X,\partial X\right]
\in H_n^{lf}\left(X,\partial X\right)$. We will say that an infinite 
chain $\sum_{i\in I}a_i\sigma_i$ represents $\left[X,\partial X\right]$ 
if it is homologous to a locally finite chain representing $\left[X,\partial 
X\right]\in H_n^{lf}\left(X,\partial X\right)$.\\
For a simplicial complex $K$, we denote by $C_k^{simp,inf}\left(K\right)$ the ${\bf R}$-vector space
of (possibly infinite) formal sums $\sum_{i\in I}a_i\sigma_i$ with 
$a_i\in{\bf R}$ and $\sigma_i$ k-simplices in $K$.
%is a locally finite simplicial chain if each $\sigma_i$ is a simplex of $K$ and 
%each vertex of $K$ is adjacent to only finitely many $\sigma_i$. 
%We can define the homology of the complex of locally finite simplicial chains. \\
If $\pi_nX^\prime\rightarrow \pi_nX$ is injective for $n\ge 1$, we have again the obvious inclusion $i_*:
C_*^{simp,inf}\left(\widehat{K}\left(X\right),\widehat{K}\left(X^\prime\right)\right)
\rightarrow
C_*^{inf}\left(X,X^\prime\right)$.\\
The following observation is of course a
well-known application of the {\bf homotopy extension property}, but we will use it so often that we state it here for reference.
\begin{obs}\label{Obs1}: Let $X$ be a topological space and
$\sigma_0:\Delta^n\rightarrow  X$ a singular simplex. Let $H:\partial\Delta^n\times I\rightarrow X$ be a homotopy
with $H\left(x,0\right)=\sigma_0\left(x\right)$ for all $x\in\partial \Delta^n$. Then there exists a homotopy $\overline {H}:\Delta^n\times I\rightarrow X$ with $\overline{H}\mid_{\partial\Delta^n\times I}=H$ and
$\overline{H}\mid_{\Delta^n\times\left\{0\right\}}=\sigma_0$.

If $X^\prime\subset X$ is a subspace and the images of $\sigma_0$ and $H$ belong to $X^\prime$, then we can choose $\overline{H}$ such that its image belongs to $X^\prime$.
\end{obs}

\begin{lem}\label{Lemma2} Let $\left(X,X^\prime\right)$ be a pair of topological spaces. Assume $\pi_nX^\prime\rightarrow\pi_nX$ is injective for each path-component of $X^\prime$ and each $n\ge 1$.\\
%Assume that each path-connec
%ted component of $X$ and $X^\prime$ contains infinitely many points or is empty. 
a) Let $\sum_{i\in I}a_i\tau_i\in C_n^{inf}\left(X,X^\prime\right)$ be a (possibly infinite) singular n-chain.
Assume that $I$ is countable, and that each path-component of $X$ and each non-empty path-component
of $X^\prime$ contain uncountably many points.
Then $\sum_{i\in I}a_i\tau_i\in C_n^{inf}\left(X,X^\prime\right)$ 
is homotopic to a (possibly infinite) simplicial chain $\sum_ia_i\tau_i^\prime\in
C_n^{simp,inf}\left(\widehat{K}\left(X\right),\widehat{K}\left(X^\prime\right)\right)$. 
In particular, $$\sum_ia_i\tau_i^\prime\in
C_n^{simp,inf}\left(\widehat{K}\left(X\right),\widehat{K}\left(X^\prime\right)\right)\subset C_*^{inf}\left(X,X^\prime\right)$$ is homologous to $\sum_{i\in I}a_i\tau_i$. \\
b) Let $\sigma_0\in \widehat{K}\left(X\right)$ and $H:\Delta^n\times I\rightarrow X$ a homotopy with $H\left(.,0\right)=\sigma_0$. Consider a minimal triangulation $\Delta^n\times I= \Delta_0\cup
\ldots\Delta_n$ of $\Delta^n\times I$ into n+1 
n+1-simplices. Assume that $H\left(\partial\Delta^n\times I\right)$ consists of simplices in $\widehat{K}\left(X\right)$. Then $H$ is homotopic (rel.\ $\Delta^n\times\left\{0\right\}\cup\partial\Delta^n\times I$) to a map $\overline{H}:\Delta^n\times I\rightarrow X$ such that $\overline{H}\mid_{\Delta_i}\in
\widehat{K}\left(X\right)$, in particular $\sigma_1:=\overline{H}\left(.,1\right)\in\widehat{K}\left(X\right)$.\\
\end{lem}
\begin{proof} 

a) From the assumptions it follows that there exists a homotopy of the 0-skeleton such that each vertex is moved
into a distinct point of $X$, and such that vertices in $X^\prime$ remain in $X^\prime$ during the homotopy.
By Observation 1, this homotopy can by induction be extended to a homotopy of the whole chain. 

Now we prove the claim by induction on $k$ ($0\le k <n$). We assume that the $k$-skeleton of $\sum_{i\in I}a_i\tau_i$
consists of simplices in $\widehat{K}\left(X\right)$ and we want to homotope 
$\sum_{i\in I}a_i\tau_i$
such that the homotoped k+1-skeleton consists of simplices in $\widehat{K}\left(X\right)$. 

By construction, each singular k+1-simplex $\sigma$ in $X$ with boundary a simplex in 
$\widehat{K}\left(X\right)$ is homotopic (rel.\ boundary)
to a unique k+1-simplex in $\widehat{K}\left(X\right)$. Since the homotopy keeps the boundary fixed, the homotopies of different k+1-simplices are compatible.
By  \hyperref[Obs1]{Observation \ref*{Obs1}}, the homotopy of the k+1-skeleton can by induction be extended to a homotopy of the whole chain.

If the image of the k+1-simplex $\sigma$ is contained in $X^\prime$, then it is homotopic 
rel.\ boundary to a simplex in $\widehat{K}\left(X^\prime\right)$, 
for a homotopy with image in $X^\prime$. Thus we can realise 
the homotopy such that all simplices with image in $X^\prime$ are homotoped inside $X^\prime$.\\
b) follows by the same argument as a), succesively applied to $\Delta_0,\ldots,\Delta_n$.

\end{proof}

We remark that there exists a canonical
simplicial map
$$p: C_*^{simp,inf}\left(\widehat{K}\left(X\right),\widehat{K}\left(X^\prime\right)\right)\rightarrow C_*^{simp,inf}\left(K\left(X\right),K\left(X^\prime\right)\right).$$
$p$ is defined by induction. It is defined to be the identity on the 1-skeleton. If it is defined on the n-1-skeleton, for $n\ge 2$, then, for an n-simplex $\tau$, $p\left(\tau\right)\in K\left(X\right)$ is the unique simplex with
$\partial_i p\left(\tau\right)=p\left(\partial_i\tau\right)$ for $i=0,\ldots,n$.

\subsection{Action of $G=\Pi\left(A\right)$}\label{sec:action}
We repeat the definitions from \cite{k2}, section 1.5.\ (originally due to \cite{gro}), as they will be frequently used in the remainder of the paper. 

Let $\left(P,A\right)$ be a pair
of minimally complete multicomplexes. \\
We define its nontrivial-loops space $\Omega^* A$ as the
set of homotopy classes (rel.\ $\left\{0,1\right\}$)
of continuous maps $\gamma:\left[0,1\right]\rightarrow \mid A\mid$ with $\gamma\left(0\right)=\gamma\left(1\right)$ and not homotopic (rel.\ $\left\{0,1\right\}$) to a
constant map.

We define
$$
\Pi\left(A\right):=\left\{\begin{array}{c}\left\{\gamma_1,\ldots,\gamma_n
\right\}
:n\in N,
\gamma_1,\ldots,\gamma_n\in A_1\cup \Omega^* A\\
\gamma_i\left(0\right)\not=\gamma_j\left(0\right),\gamma_i\left(1\right)\not=\gamma_j\left(1\right)\mbox{\ for\ }i\not= j,\\
\left\{\gamma_1\left(
0\right),\ldots,\gamma_n\left(0\right)\right\}=\left\{
\gamma_1\left(1\right),\ldots,\gamma_n\left(1\right)\right\}.
\end{array}\right\}.$$
If $\gamma,\gamma^\prime$ are elements of $A_1$ with $\gamma^\prime\not=\overline{\gamma}$ and $\gamma\left(0\right)=\gamma^\prime\left(1\right)$,
we denote\footnotemark\footnotetext[1]{We follow the usual convention to define the concatenation of paths
by $\gamma*\gamma^\prime\left(t\right)=\gamma\left(2t\right)$ if $t\le\frac{1}{2}$
and $\gamma*\gamma^\prime\left(t\right)=\gamma^\prime\left(2t-1\right)$
if $t\ge \frac{1}{2}$. Unfortunately this implies that, in order to let $\Pi\left(
A\right)$ act on $P$, we will have the multiplication in $\Pi\left(A\right)$
such that, for example, $\left\{\gamma\right\}\left\{\gamma^\prime\right\}=\left\{
\gamma^\prime*\gamma\right\}$. We hope that this does not lead to confusion.}
 $\gamma*\gamma^\prime\in A_1$ to be the unique edge of $A$ in 
the homotopy class of the concatenation. If $\gamma\in A_1$ and 
$\gamma^\prime\in\Omega^*A$ (or vice versa), 
with $\gamma\left(1\right)\not=\gamma\left(0\right)=\gamma^\prime\left(1\right)=\gamma^\prime\left(0\right)$, we also
denote $\gamma*\gamma^\prime\in A_1$ the unique edge in the homotopy class of the concatenation. 
If $\gamma,\gamma^\prime\in\Omega^*A$ with $\gamma\left(1\right)
=\gamma\left(0\right)=\gamma^\prime\left(1\right)=\gamma^\prime\left(0\right)$, we denote $\gamma*\gamma^\prime\in\Omega^*A$ the concatenation of homotopy classes of loops.\\
\\
This can be used to define a multiplication on
$\Pi\left(A\right)$ as follows:\\
\\
given $\left\{\gamma_1,\ldots,\gamma_m\right\}$
and $\left\{\gamma_1^\prime,\ldots,\gamma_n^\prime\right\}$,
we chose a reindexing of the unordered sets $\left\{\gamma_1,\ldots,\gamma_m\right\}$
and $\left\{\gamma_1^\prime,\ldots,\gamma_n^\prime\right\}$
such that we have: $\gamma_j\left(1\right)
=\gamma_j^\prime\left(0\right)$ for $1\le j\le i$ and $\gamma_j\left(1\right)
\not=\gamma_k^\prime\left(0\right)$ for $j\ge i+1, k\ge i+1$.
(Since we are assuming that all $\gamma_j\left(1\right)$ are pairwise distinct, and also all
$\gamma_j^\prime\left(0\right)$ are pairwise distinct, such a reindexing exists for some $i\ge0$, and
it is unique up to permuting the indices $\le i$ and permuting separately the
indices of $\gamma_j$'s and $\gamma_k^\prime$'s with $j\ge i+1, k\ge i+1$.) \\
Moreover we permute the indices $\left\{1,\ldots,i\right\}$ such that there exists some $h$ with $0\le h\le i$ satisfying the following conditions: \\
- for $1\le
j\le h$ we have either $\gamma_j^\prime\not=\overline{\gamma_j}\in A_1$ or $\gamma_j^\prime
\not=\gamma_j^{-1}\in\Omega^*A$,\\
- for $h<
j\le i$ we have either $\gamma_j^\prime=\overline{\gamma_j}\in A_1$ 
or $\gamma_j^\prime=\gamma_j^{-1}\in\Omega^*A$.\\

With this fixed reindexing we define
$$\left\{\gamma_1,\ldots,\gamma_m\right\}
\left\{\gamma_1^\prime,\ldots,\gamma_n^\prime\right\}:=
\left\{\gamma_1^\prime*\gamma_1,\ldots,\gamma_h^\prime*\gamma_h,\gamma_{i+1},\ldots,
\gamma_m,\gamma_{i+1}^\prime,
\ldots,\gamma_n^\prime\right\}.$$
(Note that we have omitted all $\gamma_j^\prime *\gamma_j$ with $j>h$. The choice of $\gamma_j^\prime *\gamma_j$ rather than $\gamma_j*\gamma_j^\prime$ is just because we want to define a left action on $\left(P,A\right)$.)\\
We have shown in \cite{k2} (footnote to Section 1.5.1) that the product belongs to $\Pi\left(A\right)$.
Moreover, the so-defined multiplication is independent of the chosen reindexing. It is clearly associative.
A neutral element is given by the empty set. The inverse to $\left\{\gamma_1,\ldots,\gamma_n\right\}$ is given by $\left\{\gamma_1^\prime,\ldots,\gamma_n^\prime\right\}$ with $\gamma_i^\prime=\overline{\gamma_i}$ if 
$\gamma_i\in A_1$ resp.\ $\gamma_i^\prime=\gamma_i^{-1}$ if $\gamma_i\in\Omega^*A$. 
(Indeed, in this case $h=0$, thus $\left\{\gamma_1,\ldots,\gamma_n\right\}
\left\{\gamma_1^\prime,\ldots,\gamma_n^\prime\right\}$ is the empty set.)
Thus we have defined a group law on $\Pi\left(A\right)$.\\
\\
We remark that there is an
inclusion $\Pi\left(A\right)\subset map_0\left(A_0,\left[\left[0,1\right],\mid A\mid\right]_{\mid P\mid}\right)$,
where $\left[\left[0,1\right],\mid A\mid\right]_{\mid P\mid}$ is the set of homotopy
classes (in $\mid P\mid$) rel.\ $\left\{0,1\right\}$ of maps from $\left[0,1\right]$ to $\mid A\mid$,
and $map_0\left(A_0,\left[\left[0,1\right],\mid A\mid\right]_{\mid P\mid}\right)$
is the set of maps $f:A_0\rightarrow \left[\left[0,1\right],\mid A\mid\right]_{\mid P\mid}$ with \\
- $f\left(y\right)\left(0\right)=y$ for all $y\in A_0$ and \\
- $f\left(.\right)\left(1\right):A_0\rightarrow A_0$ is a bijection.\\
This inclusion is given by sending $\left\{\gamma_1,\ldots,\gamma_n\right\}$
to the map $f$ defined by $f\left(\gamma_i\left(0\right)\right)=\left[\gamma_i\right]$ for $i=1,\ldots,n$, and $f\left(y\right)=\left[c_y\right]$
(the constant path) for $y\not\in \left\{\gamma_1\left(0\right),\ldots,\gamma_n\left(0\right)\right\}$.\\
The inclusion is a homomorphism with respect to the group law defined on
$map_0\left(A_0,\left[\left[0,1\right],\mid A\mid\right]_{\mid P\mid}\right)$
by $\left[gf\left(y\right)\right]:=\left[f\left(y\right)\right]*\left[g\left(f\left(
y\right)\left(1\right)\right)\right]$.\\

{\bf Action of $\Pi\left(A\right)$ on $P$:}\\
From now on we assume: {\em $P$ is aspherical.}
We define an action of $map_0\left(A_0,\left[\left[0,1\right],\mid A\mid\right]_{\mid P\mid}\right)$ on $P$.
This gives, in particular, an action of $\Pi\left(A\right)$ on $P$.\\
Let $g\in map_0\left(A_0,\left[\left(0,1\right),\mid A\mid\right]_{\mid P\mid}\right)$.
Define $gy=g\left(y\right)\left(1\right)$ for $y\in A_0$ and $gx=x$ for $x\in P_0-A_0$. This defines the action on the 0-skeleton of $P$.\\
We extend this to an action on the 1-skeleton of $P$:
Recall that, by minimal completeness of $P$, $1$-simplices $\sigma$ are in 1-1-correspondence with homotopy classes (rel.\ $\left\{0,1\right\}$)
of (nonclosed) singular 1-simplices in $\mid P\mid$ with vertices in $P_0$.
Using this correspondence,
define $g\sigma:=\left[\overline{g\left(\sigma\left(0\right)\right)}\right]*\left[\sigma\right]*\left[g\left(\sigma\left(1\right)\right)\right]$, where $*$ denotes concatenation of (homotopy classes of) paths.

In \cite{k2}, Section 1.5.1, we proved that this defines an action on $P_1$ and that there is an extension of ths action to an action on $P$. (The extension is unique because $P$ is aspherical.)

We remark, because this will be one of the assumptions to apply Lemma 7, that the action of any 
element $g\in\Pi\left(A\right)$ is homotopic to the identity. The homotopy between 
the action of the identity and the action of $\left\{\gamma_1,\ldots,\gamma_r\right\}$  
given by the action of
$\left\{\gamma_1^t,\ldots,\gamma_r^t\right\}, 0\le t\le 1$, 
with $\gamma_i^t\left(s\right)=\gamma_i\left(st\right)$.

The next Lemma follows directly from the construction, but we will use it so often that we want to
explicitly state it.

\begin{lem}\label{Lemma3} Let $\left(P,A\right)$ be a pair of aspherical, minimally complete multicomplexes, with the action of $G=\Pi\left(A\right)$. If $\sigma\in P$ is a simplex, all of whose vertices are not in $A$, then 
$g\sigma=\sigma$ for all $g\in G$.\end{lem}

For a topological space and a subset $P\subset S_*\left(X\right)$,
closed under face maps,
the (antisymmetric)
bounded cohomology $H_b^*\left(P\right)$
and its pseudonorm are defined literally like for multicomplexes 
in \cite{gro}, Section 3.2. The following well-known fact will
be needed for applications of  \hyperref[Lemma7]{Lemma \ref*{Lemma7}} (to the setting of  \hyperref[Thm1]{Theorem \ref*{Thm1}})
with $P=K^{str}\left(\partial Q\right),G=\Pi\left(K\left(\partial_0Q\right)\right)$.

\begin{lem}\label{Lemma4} a) Let $\left(P,A\right)$ be a pair of minimally complete multicomplexes.
If each connected component of $\mid A\mid$ has amenable fundamental group, 
then $\Pi\left(A\right)$ is amenable.\\
b) Let $X$ be a topological space, $P\subset S_*\left(X\right)$ a subset closed under face maps, 
and $G$ an amenable group acting on $P$.
Then the canonical homomorphism $$id\otimes 1:
C_*^{simp}\left(P\right)\rightarrow C_*^{simp}\left(P\right)\otimes_{{\bf Z}G}{\bf Z}$$
induces an isometric monomorphism in bounded cohomology.\end{lem}
\begin{proof} a) The proof is an obvious adaptation of the proof of \cite{k2}, Lemma 4.\\ 
b) This is proved by averaging bounded cochains, see \cite{gro}.\end{proof}

\subsection{Retraction to central simplices}\label{sec:central}

\begin{lem}\label{Lemma5}
Let $\left(N,\partial N\right)$ be a pair of topological spaces with $N=Q\cup R$ for two subspaces
$Q,R$. Let $$\partial_0Q=Q\cap R, \partial_1 Q=Q\cap\partial N,\partial_1 R=
R\cap\partial N, \partial Q=\partial_0Q\cup\partial_1Q, \partial R=\partial_0Q \cup\partial_1R.$$

Assume that $\partial_1Q\rightarrow Q\rightarrow N, \partial_1R\rightarrow R\rightarrow N,
\partial N\rightarrow N,
\partial_0Q\rightarrow Q,
\partial_0Q\rightarrow R$ are $\pi_1$-injective, and that $\partial N, \partial_1Q, \partial_1 R, \partial_0Q$ are aspherical (and thus the corresponding $K\left(.\right)$ can be considered as submulticomplexes of $K\left(N\right)$.)

Consider the simplicial action of $G=\Pi\left(K\left(\partial_0Q\right)\right)$ on $K\left(N\right)$.

Then there is a chain homomorphism $$r:C_*^{simp,inf}\left(K\left(N\right)\right)\otimes_{{\bf Z}G}{\bf Z}\rightarrow C_*^{simp,inf}\left(K\left(Q\right)
\right)\otimes_{{\bf Z}G}{\bf Z}
$$
in degrees $*\ge2$, mapping
$C_*^{simp,inf}\left(GK\left(\partial N\right)\right)\otimes_{{\bf Z}G}{\bf Z}$ to $C_*^{simp,inf}\left(GK\left(\partial_1Q\right)\right)\otimes_{{\bf Z}G}{\bf Z}$
such that\\
- if $\sigma$ is a simplex in $K\left(N\right)$, then $r\left(\sigma\otimes 1 \right)=\kappa\otimes 1$, where either
$\kappa$ is a simplex in $K\left(Q\right)$
 or $\kappa=0$,\\
- if $\sigma$ is a simplex in $K\left(Q\right)$, then $r\left(\sigma\otimes 1 \right)=\sigma\otimes 1$,\\
- if $\sigma$ is a simplex in $K\left(R\right)$, then $r\left(\sigma\otimes 1 \right)=0$.\\
\end{lem}

\begin{proof}
This is \cite{k2}, Proposition 6. (We have replaced the assumption $ker\left(\pi_1\partial_0Q\rightarrow
\pi_1Q\right)=ker\left(
\pi_1\partial_0Q\rightarrow\pi_1R\right)$ from \cite{k2} by the stronger assumption of $\pi_1$-injectivity, since 
this will be true in all our applications and 
we have no need for the more general assumption.)
The Conclusion is stated in \cite{k2} for locally finite chains, but of course $r$ extends linearly to infinite chains. 
% b) is \cite{k2}, Observation 4.
\end{proof}

{\bf Remark}: {\em If some edge of $\sigma$ is
contained in $K\left(\partial_0Q\right)=K\left(Q\right)\cap K\left(R\right)$, then $$\sigma\otimes 1=0\in C_*^{simp,inf}\left(K\left(N\right)
\right)\otimes_{{\bf Z}G}{\bf Z},$$
see \cite{k2}, Section 1.5.2. (The proof is essentially the same as that of Lemma 15 below.)
In particular, if $\sigma$ is contained in both $K\left(Q\right)$ and $K\left(R\right)$, then $r\left(\sigma\otimes 1\right)=r\left(0\right)=0$.}\\ 
\\
{\bf Fundamental cycles in $K\left(N\right)$ and $K\left(Q\right)$.} 
Let $N$ be a (possibly noncompact) connected,
orientable n-manifold with (possibly noncompact) boundary $\partial N$. Then $H_n^{lf}\left(N,\partial N\right)\simeq{\bf Z}$ by Whitehead's theorem and a generator is called $\left[N,\partial N\right]$. (It is only defined up to sign, but this will not concern our arguments.)
Recall that an infinite chain is said to represent
$\left[N,\partial N\right]$ if it is homologous to a locally finite chain representing $\left[N,\partial N\right]$.

If $\partial N\rightarrow N$ is $\pi_1$-injective and $\partial N$ is aspherical, then
$C_*^{simp,inf}\left(\widehat{K}\left(N\right),\widehat{K}\left(\partial N\right)\right)
\subset
C_*^{sing,inf}\left(N, \partial N\right)$, see Section 3.2. Thus it makes sense to say that some chain $z\in 
C_*^{simp,inf}\left(\widehat{K}\left(N\right),\widehat{K}\left( \partial N\right)\right)$ represents the
fundamental class $\left[N, \partial N\right]$. 

If $\partial_1Q\rightarrow Q$ is $\pi_1$-injective and $Q$ and 
$\partial_1Q$ are aspherical, and if $G:=\Pi\left(K\left(\partial_0Q\right)\right)$, 
then $C_*^{simp,inf}\left(
GK\left(\partial_1 Q\right)\right)=C_*^{simp,inf}\left(G\widehat{K}\left(\partial_1 Q\right)\right)\subset C_*^{sing,inf}\left(\partial Q\right)
%\subset C_*^{sing,inf}\left(\partial Q\right)
$, because $G$ maps simplices in $im\left(K\left(\partial Q\right)\rightarrow K\left(Q\right)\right)$ to simplices in $im\left(K\left(\partial Q\right)\rightarrow K\left(Q\right)\right)$.
%and 
%$C_*^{simp,inf}\left(
%GK\left(\partial_1 Q\right)\right)\subset  C_*^{sing,inf}\left(\partial Q\right)$ 
Thus it makes sense to say that some chain $z\in
C_*^{simp,inf}\left(K\left(Q\right),GK\left( \partial_1Q\right)\right)$ represents the
fundamental class $\left[Q, \partial Q\right]$.

The projection $p:\widehat{K}\left(N\right)\rightarrow K\left(N\right)$ is defined at the end of \ref{sec:construction}.

\begin{lem}\label{Lemma6} Let $N^{n\ge2}$ be an orientable n-manifold with boundary, and $Q,R\subset N$ 
orientable n-manifolds
with 
boundary, such that $N=Q\cup R$ satisfies the assumptions of Lemma 5 and
$\partial_0Q,\partial_1Q,\partial_1R$ are n-1-dimensional submanifolds (with boundary) of $\partial Q$ resp.\ $\partial R$.
Assume, in addition, that $Q$ is aspherical.\\

If $\sum_i a_i\sigma_i\in C_n^{simp,inf}\left(\widehat{K}\left(N\right),\widehat{K}\left(\partial N\right)\right)$ represents
$\left[N,\partial N\right]$, then 
$$ \sum_i a_ir\left(p\left(\sigma_i\right)\right)\otimes 1\in
C_n^{simp,inf}\left(K\left(Q\right),GK\left(\partial_1 Q\right)\right)\otimes_{{\bf Z}G}{\bf Z}$$
represents
\footnotemark\footnotetext[2]{This means that it
represents the image of
$h\otimes 1$ under the canonical homomorphism
$H_n^{sing,inf}\left(Q,\partial Q\right)\otimes_{{\bf Z}G}{\bf Z}\rightarrow H_n\left(
C_*^{sing,inf}\left(Q,\partial Q\right)\otimes_{{\bf Z}G}{\bf Z}\right)$,
where $h\in H_n^{simp,inf}\left(K\left(Q\right),
GK\left(\partial_1 Q\right)\right)$ represents $\left[Q,\partial Q\right]
\in H_n^{sing}\left(Q,\partial Q\right)$} 
$\left[Q,\partial Q\right]\otimes 1$

and 
$$\partial \sum_i a_ir\left(p\left(\sigma_i\right)\right)\otimes 1\in
C_n^{simp,inf}\left(GK\left(\partial Q\right)\right)\otimes_{{\bf Z}G}{\bf Z}$$
represents
\footnotemark\footnotetext[3]{This means that it
represents the image of
$ h\otimes 1$ under the canonical homomorphism
$H_n^{simp,inf}\left(GK\left(\partial_1 Q\right)\right)\otimes_{{\bf Z}G}{\bf Z}\rightarrow H_n\left(
C_*^{simp,inf}\left(GK\left(
\partial_1 Q\right)\right)\otimes_{{\bf Z}G}{\bf Z}\right)$, where
$h\in H_n^{simp,inf}
\left(GK\left(\partial_1 Q\right)\right)$ represents $\left[\partial Q\right]\in H_n^{sing}\left(\partial Q\right).$} 
$\left[\partial Q\right]\otimes 1$.
\end{lem}

\begin{proof} Since $p$ and $r$ are chain maps, it suffices to 
check the claim for some chosen representative of $\left[N,\partial N\right]$. So let $z\in
C_*^{simp,inf}\left(\widehat{K}\left(N\right),\widehat{K}\left(\partial N\right)\right)$ be a representative of
$\left[N,\partial N\right]$ chosen such that $$p\left(z\right)=z_Q+z_R$$ where
$z_Q$ represents $\left[Q,\partial Q\right]$ and $z_R$ represents $\left[R,\partial R\right]$ and such that 
$$\partial z_Q=w_1+w_2,\partial z_R=-w_2+w_3$$ with $w_1\in C_{n-1}^{simp,inf}\left(K\left(\partial_1Q\right)\right), 
w_2\in C_{n-1}^{simp,inf}\left(K\left(\partial_0Q\right)\right), w_3 \in C_{n-1}^{simp,inf}\left(K\left(\partial_1R\right)\right)$ representing
$\left[\partial_1Q\right], \left[\partial_0Q\right], \left[\partial_1R\right]$, respectively.

From  \hyperref[Lemma5]{Lemma \ref*{Lemma5}}:$$r\left(p\left(z\right)\otimes 1\right)=z_Q\otimes 1,$$
which implies the first claim, and
$$\partial r\left(p\left(z\right)\otimes 1\right)=\partial z_Q\otimes 1=w_1\otimes 1+w_2\otimes 1.$$
Since $w_1+w_2$ represents $\left[\partial Q\right]$, this implies the second claim.\\
\\
(Remark: From the Remark after  \hyperref[Lemma5]{Lemma \ref*{Lemma5}} we have $w_2\otimes 1=0$. This implies $\partial r\left(p\left(z\right)\otimes 1\right)=\partial z_Q\otimes 1=w_1\otimes 1$, that is, $\partial r\left(p\left(z\right)\otimes 1\right)$ represents at the same time $\left[\partial Q\right]\otimes 1$ and $\left[\partial_1Q\right]\otimes 1$.)

%Let $N^\prime:=\partial N\cup\partial_0Q$. Then we have the canonical morphism 
%$j_*:H_n\left(N,\partial N\right)\rightarrow H_n\left(N,N^\prime\right)$, which is the identity map on chains.  
%With respect to the decomposition $H_n\left(N,N^\prime\right)=i_{Q*}H_n\left(Q,\partial Q\right)\oplus 
%i_{R*}H_n\left(R,\partial R\right)$ we have 
%$$j_*\left[N,N^\prime\right]=\left(i_Q\right)_*\left[Q,Q^\prime\right]
%+\left(i_R\right)_*\left[R,R^\prime\right].$$
%By Lemma 5, we have $r\left(i_Q\otimes id\right)=id, r\left(i_R\otimes id\right)=0$, which implies a). 

%It is well-known that the boundary operator $\partial:Z_n\left(Q,\partial Q\right)\rightarrow
%Z_{n-1}\left(\partial Q\right)$ (mapping relative cycles to genuine cycles in $\partial Q$)
%induces the connecting morphism (of the long exact homology sequence) in
%homology and maps the relative fundamental class
%$\left[Q,\partial Q\right]\in H_n\left(Q,\partial Q\right)$ to
%the fundamental class $\left[\partial Q\right]\in H_{n-1}\left(\partial Q\right)$.
%Hence the boundary
%operator on $C_*^{simp}\left(K\left(Q\right)\right)$ maps (any) relative cycle representing
%$\left[Q,\partial Q\right]$ to
%a cycle representing $\left[\partial Q\right]$.
%Since $\partial$ commutes with $i_Q$,
%After tensoring with ${\bf Z}$, we obtain claim b). 
\end{proof}

\subsection{Using amenability}\label{sec:amenable}

\hyperref[Lemma7]{Lemma \ref*{Lemma7}} is well-known in slightly different formulations and we reprove it here
only for completeness.
We will apply\footnotemark\footnotetext[4]{If a group $G$ acts simplicially on a multicomplex $M$, then $C_*\left(M\right)\otimes_{{\bf Z}G}{\bf Z}$
are abelian groups with well-defined boundary operator $\partial_*\otimes 1$, even though $M/G$ may not be a multicomplex, like for the action
of $G=\Pi_X\left(X\right)$ on $K\left(X\right)$, for a topological space $X$. 

We remark that $C_*\left(M\right)\otimes_{{\bf Z}G}{\bf Z}\simeq C_*\left(M\right)\otimes_{{\bf R}G}{\bf R}$ is just the quotient chain complex for the $G$-action. In particular, even though $C_*\left(M\right)$ is an ${\bf R}G$-module, it does not make any difference whether we tensor over ${\bf Z}G$ or ${\bf R}G$.}
Lemma 7 in the proof of Theorem 1
with $X=\partial Q, G=q_*\left(\Pi\left(K\left(\partial_0Q\right)\right)\right)$ and $ K=GK^{str}\left(\partial_1 Q\right)$.
(The following lemma has of course also a relative version, but we will not need that for our argument.)

\begin{lem}\label{Lemma7}: Let $X$ be a closed, orientable manifold and $K\subset S_*\left(X\right)$ closed under face maps.
Assume that \\
- there is an amenable group $G$ acting on $K$, such that the action of each $g\in G$ on $\mid K\mid$
is homotopic to the identity\\
- there is a fundamental cycle
$z\in C_*^{simp}\left(K\right)$
such that $z\otimes 1$ is homologous to a cycle
$h=\sum_{j=1}^s b_j \tau_j\otimes 1\in C_*^{simp}
\left(K\right)\otimes_{{\bf Z}G}{\bf Z}$. \\
Then 
$$\parallel X\parallel\le \sum_{j=1}^s \mid b_j\mid.$$\end{lem}

\begin{proof}
If $\parallel X\parallel=0$, there is nothing to prove. \\
Thus we 
may assume $\parallel X\parallel\not=0$, which implies (\cite{gro}, p.17) that there is
$\beta\in H_b^n\left(X\right)$, a bounded cohomology 
class dual to 
$\left[X\right]\in H_n\left(X\right)$,
with $\parallel\beta\parallel=
\frac{1}{\parallel X\parallel}$. 

Let $p:C_*^{simp}\left(K
\right)\rightarrow C_*^{simp}\left(K
\right)\otimes_{{\bf Z}G}{\bf Z}$ be the homomorphism defined by
$p\left(\sigma\right)=\sigma\otimes 1$.
Since $G$ is amenable we have, by the proof of Lemma 4b) in \cite{gro},
an 'averaging homomorphism' $Av:H_b^*
\left(K\right)\rightarrow 
H_b^*\left(C_*\left(K\right)
\otimes_{{\bf Z}G}{\bf Z}\right)$ such that $Av$ is left-inverse to $p^*$
and $Av$ is an isometry.
Hence we have $$\parallel Av\left(\beta\right)\parallel=\parallel\beta\parallel=\frac{1}{\parallel X\parallel}.$$
Moreover, denoting by $\left[\sum_{j=1}^sb_j\tau_j\otimes 1\right]$ the homology class of 
$\sum_{j=1}^sb_j\tau_j\otimes 1$, we have obviously
$$\mid Av\left(\beta\right)\left[\sum_{j=1}^sb_j\tau_j\otimes 1\right] \mid \le
\parallel Av\left(\beta\right)\parallel \sum_{j=1}^s\mid b_j\mid$$
and therefore
$$\parallel X\parallel=\frac{1}{\parallel Av\left(\beta\right)\parallel}\le\frac{
\sum_{j=1}^s\mid b_j\mid}{\mid Av\left(\beta\right)\left[\sum_{j=1}^sb_j\tau_j\otimes 1\right]\mid}.$$
It remains to prove $Av\left(\beta\right)\left[\sum_{j=1}^sb_j\tau_j\otimes 1\right]= 1$.\\
For this we have to look at the definition of $Av$, which is as follows:
Let $\gamma\in C_b^*\left(K\right)$ be a bounded cochain. By amenability there exists a bi-invariant mean $av: B\left(G\right)\rightarrow{\bf R}$ on the bounded functions on $G$ with $inf_{g\in G} \delta
\left(g\right)\le av\left(\delta\right)\le sup_{g\in G} 
\delta\left(g\right)$ for all $\delta\in B\left(G\right)$.
Then, given any $p\left(\sigma\right)\in C_*\left(K\right)\otimes_{{\bf Z}G}{\bf Z}$ one can fix an identification between $G$ and
$G\sigma$, the set of all $\sigma^\prime$ with $p\left(\sigma^\prime\right)=p\left(\sigma\right)$, and thus consider the restriction of $\gamma$ to $G\sigma$ as a
bounded cochain on $G$. Define $Av\left(\gamma\right)\left(p\left(\sigma\right)\right)$ to be the average 
$av$ of this bounded cochain on $G\simeq G\sigma$. (This definition is independent of all choices, see \cite{iva}.)\\
Now, if $z=\sum_{j=1}^sb_j\tau_j$ is a fundamental cycle, 
then we have $\beta\left(z\right)=1$. \\
If $g\in G$ is arbitrary, then left multiplication with $g$ is a chain map
on $C_*^{simp}\left(K\right)$, as well as on $C_*^{sing}\left(X\right)$. Since the
action of $g$ on $\mid K\mid$ is
homotopic to the identity, it induces the identity on the image of $C_*^{simp}\left(K\right)\rightarrow C_*^{sing}\left(X\right)$.
%This chain map has an inverse, given by left multiplication with $g^{-1}$. Both chain maps are preserving the
%integer lattice in real homology. Since $H_n\left(X\right)$ is one-dimensional for $n=dim\left(X\right)$, 
Thus, for each cycle $z\in
C_*^{simp}\left(K\right)$ representing $\left[X\right]\in H_*^{sing}\left(X\right)$, the cycle $gz\in
C_*^{simp}\left(K\right)$ must also 
represent $\left[X\right]$. \\
If $gz$ represents $\left[X\right]$, 
then $\beta\left(gz\right)=\beta\left(\left[X\right]\right)=1$. 
%If $gz$ represents $-\left[X\right]$, then,
%since the bounded cohomology of the multicomplex $K\left(X\right)$ is defined using {\em antisymmetric}
%bounded cochains, we have $\beta\left(gz\right)=-1$. Moreover, since the action of every element $g$ is 
%homotopic to the identity, $\beta\left(gz\right)$ must not depend on $g$, for each cycle $z$.
In conclusion, we have $\beta\left(p\left(z^\prime\right)\right)=1$ for each $z^\prime$
with $p\left(z^\prime\right)=p\left(z\right)$.
%or $\beta\left(p\left(z^\prime\right)\right)=-1$ for each $z^\prime$
%with $p\left(z^\prime\right)=p\left(z\right)$. 
By definition of 
$Av$, this implies $Av\left(\beta\right)\left(p\left(z\right)\right)= 1$ for 
each fundamental cycle $z$.\\
In particular, $Av\left(\beta\right)\left[\sum_{j=1}^sb_j\tau_j\otimes 1\right]=1$, finishing the proof of the lemma.
\end{proof}

Remark: In the proof of  \hyperref[Thm1]{Theorem \ref*{Thm1}}, we will work with $C_*^{simp}\left(K\right)\otimes_{{\bf Z}G}{\bf Z}$
rather than $C_*^{simp}\left(K\right)$. This is analogous to
Agol's construction of "crushing the cusps to points"
in \cite{ag}. However 
$C_*^{simp}\left(K\left(Q\right)\right)\otimes_{{\bf Z}\Pi\left(\partial_0Q
\right)}{\bf Z}\not= C_*^{simp}\left(K\left(Q/\partial_0Q\right)\right)$,
thus one can not simplify our arguments by working directly with $Q/\partial_0Q$.

\section{Disjoint planes in a simplex}

In this section, we will discuss the possibilities how a simplex can be cut by planes 
without producing parallel arcs in the boundary. (More precisely, 
we pose the additional condition that the components of the complement can be coloured by black and white such that all vertices belong to black components, and we actually 
want to avoid only parallel arcs in the boundary of white components.)
For example, for the 3-simplex, it
will follow that there is essentially only the possibility in Case 1, pictured below, meanwhile in Case 2
each triangle has a parallel arc with another triangle, regardless how the quadrangle is triangulated. 
%(Of course, in case 1 one may have less than 4 triangles.)

\psset{unit=0.1\hsize}
$$\pspicture(0,-1)(10,3.8)
\pspolygon[linecolor=gray](0,0)(2,0)(1,4)(0,0)
\pspolygon[linecolor=gray](2,0)(1,4)(3,0.5)(2,0)
\psline(0.75,3)(1.25,3)
\psline(1.25,3)(1.3,3.5)
\psline(0.5,0)(0.375,1.5)
\psline[linestyle=dashed](0.375,1.5)(0.9,0.15)
\psline(1.5,0)(1.75,1)
\psline(1.75,1)(2.25,0.15)
\psline(2,2.25)(2.5,0.3)
\psline[linestyle=dashed](2,2.25)(2.7,0.45)
\psline[linecolor=gray,linestyle=dashed](0,0)(3,0.5)
\uput[0](0,-0.4){Case 1}

\pspolygon[linecolor=gray](5,0)(7,0)(6,4)(5,0)
\pspolygon[linecolor=gray](7,0)(6,4)(8,0.5)(7,0)
\psline[linestyle=dashed](6.25,0)(6.5,0.25)
\psline[linestyle=dashed](6.5,0.25)(6.6,3)
\psline(6.25,0)(5.5,2)
\psline[linestyle=dashed](5.5,2)(6.6,3)
\psline[linecolor=gray,linestyle=dashed](5,0)(8,0.5)
\psline(6.75,1)(7.25,0.15)
\psline(7,2.25)(7.5,0.3)
\psline[linestyle=dashed](7,2.25)(7.7,0.45)
\psline(6.75,1)(6.75,0)
\uput[0](5.6,-0.4){Case 2}

\endpspicture$$

Let $\Delta^n\subset{\bf R}^{n+1}$ be the standard simplex\footnotemark\footnotetext[5]{As usual, $v_i$ is the vertex with all coordinates, except the i-th, equal to zero, and $\partial_i\Delta^n$ denotes the subsimplex spanned by all vertices except $v_i$. We will occasionally identify singular 1-simplices $\sigma:\Delta^1\rightarrow M$ with paths $e:\left[0,1\right]\rightarrow M$ by the rule $e\left(t\right)=\sigma\left(t,1-t\right)$. In particular, $e\left(0\right)=\sigma\left(v_0\right)=\partial_1\sigma$ and $e\left(1\right)=\sigma\left(v_1\right)=\partial_0\sigma$.} with vertices $v_0,\ldots,v_n$. It is contained in the plane $E=\left\{\left(x_1,\ldots,x_{n+1}\right)\in{\bf R}^{n+1}:x_1+\ldots+x_{n+1}=1\right\}$.

In this section we will be interested in n-1-dimensional affine planes $P\subset
E$ whose intersection with $\Delta^n$ either contains no vertex, consists of exactly one vertex, or consists of a face of $\Delta^n$.
For such planes we define their type as follows.

\begin{df}\label{Def2}
Let $P\subset E$ be an n-1-dimensional affine plane such that $P\cap\Delta^n$ either contains no vertex,
consists of exactly one vertex, or consists of a face of $\Delta^n$.\\
If $P\cap\Delta^n=\partial_0\Delta^n$, then we say that $P$ is of type $\left\{0\right\}$.\\
If $P\cap\Delta^n=\partial_j\Delta^n$ with $j\ge 1$, then we say that $P$ is of type $\left\{01\ldots \hat{j}\ldots n\right\}$.\\
If $P\cap\left\{v_0,\ldots,v_n\right\}=\left\{v_0\right\}$, then we say that $P$ is of type $\left\{0\right\}$.\\
If $P\cap\left\{v_0,\ldots,v_n\right\}=\emptyset$ or $P\cap\left\{v_0,\ldots,v_n\right\}=\left\{v_j\right\}$ with $j\ge 1$, then we say that $P$ is of type $\left\{0 a_1\ldots a_k\right\}$ with $a_1,\ldots,a_k\in\left\{1,\ldots,n\right\}$ if: 

$v_i$ belongs to the same connected component of $\Delta^n - \left(P\cap\Delta^n\right)$ as 
$v_0$ 

if and only if $i\in\left\{a_1,\ldots,a_k\right\}$.
\end{df}

\begin{obs}\label{Obs2} Let $P_1,P_2$ be two planes of type $\left\{0a_1\ldots a_k\right\}$ resp.\ $\left\{0b_1\ldots b_l\right\}$ and let $Q_1=P_1\cap\Delta^n\not=\emptyset, Q_2=P_2\cap\Delta^n\not=\emptyset$. 
Then $Q_1\cap Q_2=\emptyset$ implies that 
either $\left\{a_1,\ldots,a_k\right\}= \left\{b_1,\ldots,b_l\right\}$ or exactly one of the following conditions holds:\\
- $\left\{a_1,\ldots,a_k\right\}\subset \left\{b_1,\ldots,b_l\right\}$,\\ 
- $\left\{b_1,\ldots,b_l\right\}\subset
\left\{a_1,\ldots,a_k\right\}$,\\
- $\left\{a_1,\ldots,a_k\right\}\cup \left\{b_1,\ldots,b_l\right\} =\left\{1,\ldots,n\right\}$.\end{obs}
%(equivalently
%$\left\{1,\ldots,n\right\}-\left\{b_1,\ldots,b_l\right\}\subset
%\left\{a_1,\ldots,a_k\right\}$).\end{obs}
\begin{proof} $\Delta^n-Q_1$ consists of two connected components, $C_1$ 
and $C_2$. W.l.o.g.\ assume that $v_0\in C_1$. $\Delta^n-Q_2$ 
consists of two connected components, $D_1$ and $D_2$. W.l.o.g.\ assume that $v_0\in D_1$. In particular, $C_1\cap D_1\not=\emptyset$.

Since $Q_1\cap Q_2=\emptyset$, it follows that $Q_2$ is contained in one of $C_1$ or
$C_2$, and
$Q_1$ is contained in one of $D_1$ or $D_2$.\\
Case 1: $Q_1\subset D_1$. Then either we have $C_1\subset D_1$, which implies $\left\{a_1,\ldots,a_k\right\}\subset \left\{b_1,\ldots,b_l\right\}$, or we have $C_2\subset D_1$, which implies $\left\{1,\ldots,n\right\}-
\left\{a_1,\ldots,a_k\right\}\subset \left\{b_1,\ldots,b_l\right\}$, hence 
$\left\{a_1,\ldots,a_k\right\}\cup \left\{b_1,\ldots,b_l\right\} =\left\{1,\ldots,n\right\}$.\\
Case 2: $Q_1\subset D_2$. This implies $Q_2\subset C_1$ and after interchanging $Q_1$ and $Q_2$ we are in Case 1.\end{proof}

{\bf Notational remark}: {\em 'arc'} will mean the intersection of 
an n-1-dimensional affine plane $P\subset E$ (such that $P\cap\Delta^n\not=\emptyset$ 
either contains no vertex, consists
of exactly one vertex or consists of a face) with a 2-dimensional subsimplex $\tau^2\subset\Delta^n$. If an arc consists of only one vertex, we call it a {\em degenerate arc}.

\begin{df}\label{Def3}{\bf (Parallel arcs)} Let $P_1,P_2\subset E$ be n-1-dimensional affine planes. Let $\tau^2$
be a 2-dimensional subsimplex of $\Delta^n$ with vertices $v_r,v_s,v_t$.
We say that disjoint arcs $e_1,e_2$ obtained 
as intersections of $P_1$ resp.\ $P_2$ with (the same) $\tau^2$ are parallel arcs if one of the following holds:\\
- both are nondegenerate and any two of $\left\{v_r,v_s,v_t\right\}$ belong to the same connected component of $\tau^2-e_1$ if and only if they belong to the same
connected component of $\tau^2 - e_2$,\\
- one, say $e_1$ is nondegenerate, the other, say with vertices $v_s,v_t$ is contained in a face, and $v_r$ belongs to another connected component of $\tau^2-e_1$ as both $v_s$ and $v_t$,\\ 
- one, say $e_1$, is nondegenerate, the other is degenerate, say equal to $v_r$, and both $v_s,v_t$ belong to another
connected component of $\tau^2-e_1$ as $v_r$,\\
- both are degenerate and equal,\\
- both are contained in a face and equal,\\
- one is degenerate, the other is contained in a face.\end{df}

\begin{lem}\label{Lemma8} Let $\Delta^n\subset{\bf R}^{n+1}$
be the standard simplex. Let $P_1,P_2\subset E$ be n-1-dimensional affine planes with $Q_i=P_i\cap\Delta^n\not=\emptyset$ for
$i=1,2$.\\
Let $P_1$ be of type $
\left\{0a_1\ldots a_k\right\}$ with $1\le k\le n-2$ and $P_2$ of type $\left\{0b_1\ldots b_l\right\}$ with $l$ arbitrary.\\
Then either $Q_1\cap Q_2\not=\emptyset$, or $Q_1$ and $Q_2$ have a parallel arc.\end{lem}

\begin{proof}
Assume that $Q_1\cap Q_2=\emptyset$.\\
By  \hyperref[Obs2]{Observation \ref*{Obs2}}, there are 4 possible cases if $Q_1\cap Q_2=\emptyset$. \\
Case 1: $\left\{0a_1\ldots a_k\right\}=\left\{0b_1\ldots b_l\right\}$. Then we clearly have parallel arcs.\\
Case 2: $\left\{0a_1\ldots a_k\right\}$ is a proper subset of $\left\{0b_1\ldots b_l\right\}$, i.e.\ 
$1\le k<l\le n-1$ and $a_1=b_1,\ldots,a_k=b_k$. There is at least one index, say $i$, not contained in 
$\left\{0b_1\ldots b_l\right\}$. Consider the 2-dimensional subsimplex $\tau^2\subset
\Delta^n$ with vertices $v_0,v_{a_1},v_i$. It intersects $P_1$ and $P_2$ in parallel arcs, because $P_1$ and $P_2$ both separate $v_0$ and $v_{a_k}$ from $v_i$. \\
Case 3: $\left\{0b_1\ldots b_l\right\}$ is a proper subset of $\left\{0a_1\ldots a_k\right\}$, i.e.\ 
$ 0\le l <k\le n-2$ and $a_1=b_1,\ldots,a_l=b_l$. There are two indices $i,j$ not contained in 
$\left\{0a_1\ldots a_k\right\}$. Consider the 2-dimensional subsimplex $\tau^2\subset
\Delta^n$ with vertices $v_0,v_i,v_j$. It intersects $P_1$ and $P_2$ in parallel arcs, because $P_1$ and $P_2$ both separate $v_0$ from $v_i$ and $v_j$.\\
Case 4: 
$\left\{a_1,\ldots,a_k\right\}\cup\left\{b_1,\ldots,b_l\right\}=\left\{1,\ldots,n\right\}$.
By $k\le n-2$, there are two indices $i,j$ with
$i,j\not\in\left\{0a_1\ldots a_k\right\}$. Hence $i,j\in \left\{b_1,\ldots,b_l\right\}$. Moreover, there exists
an index $h$ such that $h\in\left\{a_1,\ldots,a_k\right\}$ but
$h\not\in\left\{b_1,\ldots,b_l\right\}$. (If not, we would have $\left\{a_1,\ldots,a_k\right\}\subset 
\left\{b_1,\ldots,b_l\right\}$, hence $\left\{1,\ldots,n\right\}=
\left\{a_1,\ldots,a_k\right\}\cup 
\left\{b_1,\ldots,b_l\right\}\subset\left\{b_1,\ldots,b_l\right\}$,
contradicting $Q_2\not=\emptyset$.) Consider the 2-dimensional subsimplex $\tau^2\subset
\Delta^n$ with vertices $v_i,v_j,v_h$. It intersects $P_1$ and $P_2$ in parallel arcs, because both $P_1$ and $P_2$ separate $v_i$ and $v_j$ from $v_h$..

\end{proof} 

\begin{df}\label{Def4} {\bf (Canonical colouring of complementary regions)}\\
Let $P_1,P_2,\ldots\subset E$ be a (possibly infinite)
set of n-1-dimensional affine planes with
$Q_i:=P_i\cap \Delta^n\not=\emptyset$ and $Q_i\cap Q_j=\emptyset$ for all $i\not=j$. Assume that each $Q_i$ either
contains no vertices or consists of exactly one vertex.

A colouring of\\
- the connected components
of $\Delta^n-\cup_i Q_i$ by colours black and white, and \\
- of all $Q_i$ by black,\\
is called a canonical colouring (associated to $P_1,P_2,\ldots$) if:\\
- all vertices of $\Delta^n$ are coloured black,\\
- each $Q_i$ is incident to at least one white component. \end{df}
%It is clear that a canonical colouring is unique when it exists.

\begin{df}\label{Def5} {\bf (White-parallel arcs)} Let $\left\{P_i: i\in I\right\}$ be a set of
of $n-1$-dimensional affine planes $P_i\subset E$, with $Q_i:=P_i\cap\Delta^n\not=\emptyset$ for $i\in I$.
Assume that $Q_i\cap Q_j=\emptyset$ for all $i\not=j\in I$, and that we
have a canonical colouring associated to $\left\{P_i: i\in I\right\}$.
We say that arcs $e_i,e_j$ obtained
as intersections of $P_i,P_j$ ($i,j\in I$) with some 2-dimensional subsimplex
$\tau^2$ of $\Delta^n$ are white-parallel arcs if they are parallel arcs
and, moreover, belong to the boundary of the closure of the same white component.\end{df}

We mention two consequences of  \hyperref[Lemma8]{Lemma \ref*{Lemma8}}. These will not be needed for the proof of  \hyperref[Lemma10]{Lemma \ref*{Lemma10}}, but they will be necessary for the proof of  \hyperref[Thm1]{Theorem \ref*{Thm1}}.

\begin{cor}\label{Cor1} Let $\Delta^n\subset{\bf R}^{n+1}$ 
be the standard simplex. Let $P_1,\ldots,P_m\subset E$ be a finite set of
n-1-dimensional
affine planes and let $Q_i=P_i\cap\Delta^n$ for 
$i=1,\ldots,m$. 

Assume that $Q_i\cap Q_j=\emptyset$ for all $i\not=j$, and that
we have an associated 
canonical colouring, such that
$Q_i$ and $Q_j$ 
do not have a white-parallel arc for $i\not=j$.

Then either $m=0$, or\\
$m=n+1$ and $P_1$ is of type $\left\{0\right\}$, $P_{n+1}$ is of type $\left\{0\ 1\ldots n-1\right\}$,
and $P_i$ is of type
$\left\{01\ldots\widehat{i-1}\ldots n\right\}$ for $i=2,\ldots,n$.\end{cor}

\begin{proof}
If the conclusion were not true, there would exist a plane $P_1$ of type $\left\{0a_1\ldots a_k\right\}$ with $1\le k\le n-2$. Let $W$ be the white component of the canonical colouring, which is incident to $P_1$. Because, for a canonical colouring, no vertex belongs to a white component, there must be at least one more plane $P_2$ incident to $W$. Since $Q_1\cap Q_2=\emptyset$, from Lemma 8
we get that $Q_1$ and $Q_2$ have a parallel arc. Because $Q_1$ and $Q_2$ are incident to $W$, the arc is white-parallel.\end{proof}

\begin{cor}\label{Cor2} Let $\Delta^n\subset{\bf R}^{n+1}$                                                                   be 
the standard simplex. Let $P_1,P_2,\ldots\subset E$ be a (possibly infinite) set of
n-1-dimensional                                                                                                    affine 
planes and let $Q_i=P_i\cap\Delta^n$ for             
$i=1,2,\ldots$.  
Assume that we have an associated canonical colouring.
\\
Let $P_i$ be of type 
$                        \left\{0a_1^i\ldots a^i_{c\left(i\right)}\right\}$, for $i=1,2,\ldots$.
Then \\
- either $c\left(1\right)\in\left\{0,n-1\right\}$,\\
- or whenever, for some $i\in\left\{2,3,\ldots\right\}$,
$P_1$ and $ P_i$ bound a white component of $\Delta^n-\cup_jQ_j$, then they must have a white-parallel arc.\end{cor}

\begin{proof} Assume that $c\left(1\right)\not \in\left\{0,n-1\right\}$. 
Let $W$ be the white component bounded by $P_1$.
%In particular, 
%there is no connected component of $\Delta^n-\cup_i Q_i$ bounded by n+1 planes of type $\left\{0\right\},
%\left\{023\ldots n\right\},\left\{013\ldots n\right\},\ldots,\left\{012\ldots n-2,n\right\},\left\{012\ldots n-2,n-1\right\}$. 
%Thus each white component $W$
%is bounded by at least one plane $P_1$ of type $\left\{0a_1^1\ldots a^1_{c\left(1\right)}\right\}$ with $1\le c\left(1\right)\le n-2$. 
$W$ is bounded by a finite number of planes, thus we can
apply  \hyperref[Cor1]{Corollary \ref*{Cor1}}, and conclude that $P_1$ has a 
white-parallel arc with each other plane adjacent to $W$.  \end{proof}

\begin{df}\label{Def6} Let $P\subset E$ be an n-1-dimensional affine plane, and $T$ a triangulation of the 
polytope $Q:=P\cap\Delta^n$. We say that $T$ is minimal, if all vertices of $T$ are vertices of $Q$. 
We say that an edge of some simplex in $T$ is an exterior edge if it is an edge of $Q$.\end{df}

\begin{obs}\label{Obs3} Let $P\subset E$ be an n-1-dimensional affine plane, and $T$ a triangulation of the
polytope $Q:=P\cap\Delta^n$. If $T$ is minimal, then each edge of $Q$ is an (exterior)
edge of (exactly one) simplex in $T$.\end{obs}
\begin{proof} By minimality, the triangulation does not introduce new vertices. Thus every edge of $Q$ is an edge of some simplex.\end{proof}

\begin{obs}\label{Obs4} Let $P\subset E$ be an n-1-dimensional affine plane with $Q:=P\cap \Delta^n\not=\emptyset$.
Assume that $P$ is of type
$\left\{0a_1\ldots a_k\right\}$.\\
a) Each vertex of $Q$ arises as the intersection of $P$ with an edge $e$ of $\Delta^n$. The vertices of $e$ are $v_i$ and $v_j$ with $
i\in\left\{0,a_1,\ldots,a_k\right\}$ and
$j\not\in \left\{0,a_1,\ldots,a_k\right\}$. (We will denote such a vertex by $\left(v_iv_j\right)$.)\\
b) Two vertices
$\left(v_{i_1}v_{j_1}\right)$ and $\left(v_{i_2}v_{j_2}\right)$ of $Q$ are connected by an edge of $Q$ (i.e.\ an exterior edge of any triangulation)
if either $i_1=i_2$ or $j_1=j_2$. \end{obs}
\begin{proof} a) holds because $e$ has to connect vertices in distinct components of $\Delta^n-Q$. b) holds because the edge of $Q$ has to belong to some 2-dimensional subsimplex of $\Delta^n$, with vertices either $v_{i_1},v_{j_1},v_{j_2}$ or $v_{i_1},v_{i_2},v_{j_1}$.\end{proof}

Remark: if, for an affine hyperplane $P\subset E$, $Q=P\cap \Delta^n$ consists of exactly one vertex, then we will consider the minimal triangulation
of $Q$ to consist of one (degenerate) n-1-simplex. This convention helps to avoid needless case distinctions.

\begin{lem}\label{Lemma9} 
Let $\left\{P_i\subset E:i\in I\right\}$ be a set of n-1-dimensional affine
planes and let $Q_i:=P_i\cap\Delta^n$ for $i\in I$.
Assume that $Q_i\cap Q_j=\emptyset$ for all $i\not=j$ and that
we have an associated
canonical colouring.
Assume that we have fixed, for each $i\in I$, a minimal triangulation $Q_i=\cup_a\tau_{ia}$ of $Q_i$. 

If $P_1$ is of type $\left\{0a_1^1\ldots a_{c\left(1\right)}^1\right\}$ with $1\le c\left(1\right)\le n-2$, then for each simplex $\tau_{1a}\subset Q_1$
there exists some $j\in I$ and some simplex $\tau_{jb}\subset Q_j$ (of the fixed triangulation of $Q_j$)
such that $\tau_{ia}$ and $\tau_{jb}$
have a white-parallel arc.\end{lem}

\begin{proof} 
%Let $P_i$ be of type
%$\left\{0a_1^i\ldots a_{c\left(i\right)}^i\right\}$. 
Let $w_1,\ldots,w_n$ be the $n$ vertices of the n-1-simplex $\tau_{1k}$.  
By Observation 4a), each 
$w_l$ arises as intersection of $Q_1$ with some edge $\left(v_{r_l}v_{s_l}\right)$ of $\Delta^n$, and the 
vertices $v_{r_l},v_{s_l}$ 
satisfy $r_l\in \left\{0,a_1^1,\ldots, a_{c\left(1\right)}^1\right\}$ and
$s_l\not\in \left\{0,a_1^1,\ldots, a_{c\left(1\right)}^1\right\}$. 

For the canonical colouring, there must be a white component 
$W$ bounded by $P_1$. We distinguish the cases whether $W$ and $v_0$ belong to the same connected component
of $\Delta^n-Q_1$ or not.\\
Case 1: $W$ and $v_0$ belong to the same connected component
of $\Delta^n-Q_1$. 

Since $c\left(1\right)\le n-2$, there exist at most $n-1$ possible values for $r_l$.
Hence there exists $l\not=m\in\left\{1,\ldots,n\right\}$
such that $v_{r_l}=v_{r_m}$. 

Let $e$ be the edge of $\tau_{1k}\subset Q_1$ connecting $w_{l}$ and $w_{m}$. By Observation 4b), $e$ is an exterior edge.
Consider the 2-dimensional subsimplex $\tau^2\subset\Delta^n$ with vertices $v_{r_l},v_{s_l},v_{s_m}$. We have that $P_1$ intersects $\tau^2$ in $e$, i.e.\ in an
arc separating $v_{r_l}$ from the other two vertices of $\tau^2$. 

Note that $r_l\in \left\{0,a_1^1,\ldots, a_{c\left(1\right)}^1\right\}$, hence $v_{r_l}$
belongs to the same component of $\Delta^n-Q_1$ as $v_0$. In particular,
$v_{r_l}$
belongs to the same component of $\Delta^n-Q_1$ as $W$. On the other hand, since the colouring is canonical, all vertices are coloured black and
$v_{r_l}$ can not belong to the white component $W$. Thus 
there must be some plane $P_j$ such that $Q_j$ bounds $W$ and separates $v_{r_l}$ from $Q_1$. (The possiblity $P_j\cap\Delta^n=\left\{v_{r_l}\right\}$ is allowed.)
In particular, some (possibly degenerate)
exterior edge $f$ of $Q_j$ separates $v_{r_l}$ from $v_{s_l},v_{s_m}$. Thus $e$ and $f$ are white-parallel arcs.
By  \hyperref[Obs3]{Observation \ref*{Obs3}}, $f$ is an edge of some $\tau_{jl}$.\\
\\
Case 2:
$W$ and $v_0$ don't belong to the same connected component
of $\Delta^n-Q_1$.

Since $n-c\left(1\right)\le n-1$, there exist some $l\not=m\in\left\{1,\ldots,n\right\}$
such that $v_{s_l}=v_{s_m}$.

Let $e$ be the edge of $\tau_{1k}\subset Q_1$ connecting $w_{l}$ and $w_{m}$. $e$ is an 
exterior edge by  \hyperref[Obs4]{Observation \ref*{Obs4}}b).
Consider the 2-dimensional subsimplex $\tau^2\subset\Delta^n$ with vertices $v_{r_l},v_{r_m},v_{s_l}$.
$P_1$ intersects $\tau^2$ in $e$, i.e.\ in an
arc separating $v_{s_l}$ from the other two vertices of $\tau^2$. 

We have that $s_l\not\in \left\{0,a_1^1,\ldots, a_{c\left(1\right)}^1\right\}$, hence $v_{s_l}$
does not belong to the same component of $\Delta^n-Q_1$ as $v_0$. This implies that
$v_{s_l}$
belongs to the same component of $\Delta^n-Q_1$ as $W$. On the other hand, since the colouring is canonical, $v_{s_l}$ 
can not belong to the white component $W$ and there must be 
some plane $P_j$ such that $Q_j$ bounds $W$ and separates $v_{s_l}$ from $Q_1$. 
In particular, some exterior edge $f$ of $Q_j$ separates $v_{s_l}$ from $v_{r_l},v_{r_m}$. Thus $e$ and $f$ are white-parallel arcs.
By  \hyperref[Obs3]{Observation \ref*{Obs3}}, $f$ is an edge of some $\tau_{jl}$.\end{proof}

\begin{lem}\label{Lemma10}
Let $\left\{P_i:i\in I\right\}$ be a set of n-1-dimensional affine
planes with $Q_i:=P_i\cap\Delta^n\not=\emptyset$ for $i\in I$. Let $P_i$ be of type $\left\{
0a_1^{\left(i\right)}\ldots a^{\left(i\right)}_{k_i}\right\}$ for $i\in I$.
Assume that
$Q_i\cap Q_j=\emptyset$ for $i\not=j\in I$,
and that we have an associated canonical colouring.
Assume that for each $Q_i$ one has fixed a minimal triangulation $Q_i=\cup_{k=1}^{t\left(i\right)}\tau_{ik}$.

For each $i\in I$, let $$D_i=\sharp\left\{\tau_{ik}\subset Q_i:\mbox{\ there\ is\ no\ }\tau_{jl}\subset Q_j\mbox{\ such\ that\ }\tau_{ik},\tau_{jl}\mbox{\ have\ a\ white-parallel\ arc}\right\}.$$
Then $$\sum_{i\in I} D_i=0\ \ or\ \ \sum_{i\in I} D_i = n+1.$$
%In particular, the minimal triangulations of $Q_1,Q_2,\ldots$ contain altogether at most $n+1$ simplices
%without white-parallel arcs.
\end{lem}

 \begin{proof}
First we remark that the number of planes may be infinite, but we may of course remove pairs
of planes $P_i,P_j$ whenever they are of the same type and bound the same white component. This removal (of $P_i,P_j$ and the common white component) 
does not affect $\sum_{i\in I} D_i$. Since there are only finitely many different types of planes,
we may w.l.o.g.\ assume that we start with a finite number $P_1,\ldots,P_m$ of planes. (It may happen 
that after this removal no planes and no white components remain. In this case 
$\sum_{i\in I} D_{i\in I}=0$.) So we assume now that we have a finite number 
of planes $P_1,\ldots,P_m$, and no two planes of the same type bound a white region.

The first case to consider is that all planes are of type $\left\{0a_1\ldots a_k\right\}$ with 
$k=0$ or $k=n-1$. Since all vertices are coloured black, this means that $m=n+1$ and (upon renumbering) 
$P_1$ is of type $\left\{0\right\}$, $P_{n+1}$ is of type $\left\{0\ 1\ldots n-1\right\}$,
and $P_i$ is of type
$\left\{01\ldots\widehat{i-1}\ldots n\right\}$ for $i=2,\ldots,n$. 
Hence $D_1=\ldots=D_{n+1}=1$ and $\sum_{i=1}^{n+1} D_i=n+1$.

Now we assume that there exists $P_i$, w.l.o.g.\ $P_1$,
of type 
$\left\{0a_1^{\left(1\right)}\ldots a_{k_1}^{\left(1\right)}\right\}$ with $1\le c\left(1\right)\le n-2$. Let $W$ be the 
white component bounded by $P_1$ and let w.l.o.g.\ $P_2,\ldots,P_l$ be the other planes bounding $W$. 
Then  \hyperref[Lemma9]{Lemma \ref*{Lemma9}} says that each simplex in the chosen triangulation of $Q_1$ has a parallel arc with some simplex in the chosen triangulation of each of $Q_2,\ldots,Q_l$. 
In particular, $D_1=0$. For $j\in\left\{2,\ldots,l\right\}$,
if $1\le c\left(j\right)\le n-2$, the same argument shows that $D_j=0$. 
If $j\in\left\{2,\ldots,l\right\}$ and
$c\left(j\right)=0$ or $c\left(j\right)=n-1$, 
then $Q_j$ consists of only one simplex. By \hyperref[Cor2]{Corollary \ref*{Cor2}}, this simplex has a parallel arc with 
(some exterior edge of) $Q_1$ and thus
(by  \hyperref[Obs3]{Observation \ref*{Obs3}}) with (some) simplex of the chosen triangulation of $Q_1$. This shows $D_j=0$ 
also in this case. Altogether we conclude $\sum_{j=1}^l D_j=0$ and 
thus $\sum_{i=1}^m D_i=\sum_{i=l+1}^m D_i$. Hence we can remove\footnotemark\footnotetext[6]{To remove a white component means that this component together with the neighbouring black components will form {\bf one} new black component.} the white component $W$ and its
bounding planes $P_1,\ldots,P_l$ to obtain a smaller number of planes and
a new canonical colouring without changing $\sum_{i=1}^m D_i$.
Since we start with finitely many planes, we can repeat this reduction finitely many times and will end up either with an empty set of planes or with a set of planes of type $\left\{0a_1\ldots a_k\right\}$ with $k=0$ or $k=n-1$. Thus either $\sum_{i=1}^m D_i=0$ or $\sum_{i=1}^mD_i=n+1$.
\end{proof}

We have thus proved that, in presence of a canonical colouring, the number of n-1-simplices without white-parallel arcs in a minimal triangulation of the $Q_i$'s is $0$ or $n+1$.
We remark that in the proof of \hyperref[Thm1]{Theorem \ref*{Thm1}} we will actually count only those triangles which neither have a white-parallel arc nor a degenerate arc. Thus, in general, we may remain with 
even less than $n+1$ n-1-simplices.

\section{A straightening procedure}

In this section we will always work with the following set of assumptions.\\
\\
{\bf Assumption I}: {\em $Q$ is an aspherical n-dimensional
manifold with aspherical boundary $\partial Q$.
We have n-1-dimensional submanifolds $\partial_0Q,\partial_1Q\subset\partial Q$ such that
$\partial Q=\partial_0Q\cup\partial_1Q, \partial\partial_0Q=\partial\partial_1Q$ and
$\partial_1Q\not=\emptyset$ is aspherical.}\\
\\
The example that one should have in mind is a nonpositively curved manifold $Q$ with totally geodesic boundary $\partial_1Q$ and cusps corresponding to $\partial_0Q$.

In the case of nonpositively curved manifolds with totally geodesic boundary, there is a well-known 
straightening procedure (explained for closed hyperbolic manifolds in \cite{bp}, Lemma C.4.3.), which homotopes each relative cycle into a straight 
relative cycle. 

However, we will need a more subtle straightening procedure, which considers 
relative cycles with a certain 0-1-labeling of their edges and straightens the 1-labeled edges into certain distinguished 1-simplices.
This straightening procedure will be explained in Section \ref{sec:straighten}. Before, we  explain a construction which will
morally (although not literally) "reduce" the proof of  \hyperref[Thm1]{Theorem \ref*{Thm1}} to the case that $\partial_0Q\cap C$ is path-connected, for each path-component $C$ of
$\partial Q$. 

\subsection{Making $\partial_0Q\cap C$ connected}\label{sec:connected}

\begin{cons}\label{Cons1} Let Assumption I be satisfied. Then there exists a continuous map of triples $q:\left(Q,\partial Q,\partial_1 Q\right)\rightarrow \left(Q,\partial Q,\partial_1 Q\right)$ 
which is
(as a map of triples) homotopic to the identity and such that, for each path-component $C$ of $\partial Q$, the image $A:= q\left(\partial_0Q\cap C\right)$ is path-connected.

Moreover, for each path-component $F$ of $\partial_1Q$, the path-components of $\partial F\subset \partial_0Q\cap\partial_1Q$
can be numbered by $E_0^F,\ldots, E_s^F$ and one can choose points $x_{E_i^F}\in E_i^F$ 
%and open, contractible neighborhoods $U_i\subset E_i$ of $x_{E_i}$ for $i=0,\ldots,s$ 
such that $q\left(x_{E_i^F}\right)\equiv x_{E_0^F}$ for $i=0,\ldots,s$.\end{cons}

\begin{proof} 
For each path-component $F$ of $\partial_1Q$, number the path-components of $\partial F\subset \partial_0Q\cap\partial_1Q$ 
by $E_0^F,\ldots, E_s^F$, where $s$ depends on $F$. Choose one point $x_E^F\in E$ for each path-component $E\subset F$ of $\partial_0 Q\cap\partial_1 Q$. Whenever $E_0,E_i$ is a pair of path-components of $\partial_0Q\cap\partial_1 
Q$ adjacent to the same path-component $F$ of $\partial_1Q$, choose a 1-dimensional submanifold $l_{E_0^FE_i^F}\subset\partial_1Q$ with
$\partial l_{E_0^FE_i^F}=
\left\{x_{E_0^F}\right\}\cup\left\{x_{E_i^F}\right\}$. 
The $l_{E_0^FE_i^F}$ may be chosen succesively such that they are disjoint from each other (apart from the common vertex $x_{E_0^F}$) 
and disjoint from $\partial_0Q$ (apart from the vertices $x_{E_0^F}$ and $x_{E_i^F}$). \\
%Moreover we choose some neighborhood $U_i$ of $x_{E_i}$ in $E_i$ that is homeomorphic to an embedded $k$-ball, $k=dim\left(\partial_0Q\cap\partial_1Q\right)=dim\left(Q\right)-2$.

For each pair $\left\{E_0^F,E_i^F\right\}$ 
let $h:l_{E_0^FE_i^F}\rightarrow \left\{x_{E_0^F}\right\}$ be the constant map from $l_{E_0^FE_i^F}
%\cup U_i
$ to $x_{E_0^F}$. 
%$h$ is homotopic to the identity because $l_{E_0^FE_i^F}
%\cup U_i
%$ is contractible.
For each path-component $F$ of $\partial_1Q$, the union $$\bigcup_{i=1}^sl_{E_0^FE_i^F}$$ is an embedded wedge of arcs in
$\partial_1Q$, hence it is contractible. In particular, $h$ is homotopic to the identity. By the homotopy extension property exists $g:F\rightarrow 
F$ with $g\mid_{l_{E_0^FE_i^F}}=h\equiv x_{E_0}$ for all $l_{E_0^FE_i^F}$, and $g\sim id$ 
by a homotopy extending the homotopies between $h$ and $id$. 

Thus we defined $g$ on each path-component $F$ of $\partial_1Q$ with $F\cap\partial_0Q\not=\emptyset$. On path-components $F$ of 
$\partial_1Q$ with $F\cap\partial_0Q=\emptyset$ we define $g=id$.
Hence we have defined $g$ on all of $\partial_1Q$.

On path-components $C$ of $\partial_0Q$ with $C\cap\partial_1Q=\emptyset$, we define $f=id$.
Again by the homotopy extension property exists $f:\partial Q\rightarrow\partial Q$ with $f\mid_{\partial_1Q}=g$, $f\mid_C=id$ for path-components $C$ of $\partial_0Q$ with $C\cap\partial_1Q=\emptyset$,
and $f\sim id$ by a homotopy
extending the homotopy of $g$. (Of course, $f$ does not preserve those path-components of $\partial_0Q$ which intersect $\partial_1Q$.)

Once again
by the homotopy extension property exists $q:Q\rightarrow Q$ with $q\sim id$ such that $q$ extends 
$f$ and the homotopy between $q$ and $id$ extends the homotopy between $f$ and $id$.

Due to the stepwise construction, $q$ is a map of triples, homotopic to the identity by a homotopy of triples.
Moreover, $A:= q\left(\partial_0Q\cap C\right)$ is path-connected for each component $C$ of $\partial Q$. Indeed, any two points in $\partial_0Q\cap C$ 
can be connected by a sequence of paths which either have image in $\partial_0Q$ 
or belong to $\cup_{i=1}^sl_{E_0^FE_i^F}$ for some path-component
$F$ of $\partial_1Q\cap C$. The image of these paths under $q$, in both cases, is in $A$.
\end{proof}

{\bf Remark:} {\em $q$ induces a simplicial map $q:K\left(Q\right)\rightarrow K\left(Q\right)$ and a homomorphism $q_*:\Pi\left(K\left(\partial_0Q\right)\right)\rightarrow
\Pi\left(K\left(A\right)\right)$ defined by $q_*\left(\left\{\gamma_1,\ldots,\gamma_n\right\}\right):=\left\{q\left(\gamma_1\right),\ldots,q\left(\gamma_n\right)\right\}$ such that $$q_*\left(g\right)q\left(\sigma\right)=q\left(g\sigma\right)$$
holds for each $\sigma\in K\left(Q\right), g\in \Pi\left(K\left(\partial_0Q\right)\right)$.}

\begin{proof} Continuous maps $q:Q\rightarrow Q$ induce simplicial maps $q:K\left(Q\right)\rightarrow K\left(Q\right)$. (The simplicial map agrees with $q$ on the 0-skeleton, and it maps each 1-simplex $e\in K_1\left(Q\right)$ to the unique 1-simplex of $K_1\left(Q\right)$ that is in the homotopy class rel.\ $\left\{0,1\right\}$ of
$q\left(e\right)$.) 

Let $e\in K_1\left(Q\right)$. By construction $\left\{\gamma_1,\ldots,\gamma_n\right\}e=\left[\alpha*e*\overline{\beta}\right]$ for some $\alpha,\beta\in\left\{\gamma_1,\ldots,\gamma_n\right\}\cup\left\{c_{e\left(0\right)},c_{e\left(1\right)}\right\}$. Thus 
$$\left\{q\left(\gamma_1\right),\ldots,q\left(\gamma_n\right)\right\}q\left(e\right)=\left[q\left(\alpha\right)*q\left(e\right)*q\left(\overline{\beta}\right)\right\}=q\left(\left\{\gamma_1,\ldots,\gamma_n\right\}e\right).$$
This implies the claim for the 1-skeleton and thus, by asphericity of $K\left(Q\right)$, for all $\sigma\in K\left(Q\right)$.\end{proof}

\subsection{Definition of $K^{str}\left(Q\right)$}

Let $Q,\partial Q, \partial_1 Q, \partial_0 Q$ satisfy Assumption I.

Recall that we have defined
in Section \ref{sec:construction} an aspherical multicomplex $K\left(Q\right)\subset
 S_*\left(Q\right)$ with the property that (for aspherical $Q$) each singular simplex in $Q$, with boundary in $K\left(Q\right)$ and pairwise distinct vertices, is homotopic rel.\ boundary to a unique simplex 
in $K\left(Q\right)$. 

The aim of this subsection is to describe a selection
procedure yielding
a subset $K_*^{str}\left(Q\right)\subset S_*\left(Q\right)$.
The final purpose of the straightening procedure will be to produce a large number of 
(weakly) degenerate simplices, in the sense of the following definition.
\begin{df}\label{Def7} Let $Q$ be an compact manifold with boundary $\partial Q$.
We say that a simplex in $S_*\left(Q\right)$ is degenerate
if
one of its edges is a constant loop. We say that it is weakly degenerate if it is degenerate 
or its image is contained in $\partial Q$.
\end{df}

Notational remark:
%\footnotemark[3]\footnotetext[3]{For the proof of Theorem 1 it will, in view of Lemma 6, be sufficient to consider relative cycles in 
%$C_*\left(K\left(Q\right),GK\left(\partial_1Q\right)\right)\otimes_{{\bf Z}G}{\bf Z}$. 
%Since we prefer not to consider the G-action before Secton 5.4, but on the other hand $K\left(\partial Q\right)$ need not be a 
%submulticomplex of $K\left(Q\right)$, we will restrict ourselves to relative fundamental cycles $z\in C_*^{simp,inf}\left(K\left(Q\right)\right)$ with 
%the property that $\partial z=w_0+w_1, w_0\in C_*^{simp,inf}\left(K\left(\partial_0Q\right)\right), w_1\in C_*^{simp,inf}\left(
%K\left(\partial_1Q\right)\right)$.}
for subsets $K_*^{str}\left(Q\right)\subset
S_*\left(Q\right)$ we will denote $K_*^{str}\left(\partial_0
Q\right):=K_*^{str}\left(Q\right)\cap S_*\left(\partial_0 Q\right), K_*^{str}\left(\partial_1
Q\right):=K_*^{str}\left(Q\right)\cap S_*\left(\partial_1 Q\right),
K_*^{str}\left(\partial_0Q
%\cap \partial_1
Q
\right):=K_*^{str}\left(Q\right)\cap S_*\left(\partial_0Q
%\cap\partial_1 Q
\right).$

\begin{lem}\label{Lemma11} Let $Q,\partial Q, \partial_1 Q, \partial_0Q$ satisfy Assumption I. 
%Let $q:Q\rightarrow Q$ be given by  \hyperref[Cons1]{Construction \ref*{Cons1}} and let $A:=q\left(\partial_0 Q\right)$.
Let $K\left(Q\right)\subset  S_*\left(Q\right)$ be as defined in Section \ref{sec:construction}. 
Let $q:Q\rightarrow Q$  and $\left\{x_{E_i^F}\in \partial_0Q\cap\partial_1Q : 0\le i\le s\right\}$ be given by \hyperref[Cons1]{Construction \ref*{Cons1}}.

Then there exists a subset $K_*^{str}\left(Q\right)\subset
S_*\left(Q\right)$,
closed under face maps, 
such that:\\
i) If $C$ is a path-component of $\partial_0Q$ with $C\cap\partial_1Q=\emptyset$, then $K_0^{str}\left(Q\right)$ contains each point in $C$, \\
%and, for path-components $C$ of $\partial Q$ with $C\cap\partial_0Q\not=\emptyset$ it contains no point in $C\setminus\partial_0Q$,\\
ii) for a path-component $F$ of $\partial_1Q$ with $F\cap\partial_0Q=\emptyset$, there is exactly one point $x_F\in K_0^{str}\left(Q\right)\cap F$,\\
for a path-component $F$ of $\partial_1Q$ with $F\cap\partial_0Q\not=\emptyset$, we have
$K_0^{str}\left(Q\right)\cap F=\left\{x_{E_0^F},\ldots,x_{E_s^F}\right\}$,\\
% where $x_{E_0^F},\ldots,x_{E_s^F}$ are given by Construction 1,\\
iii) $K_0^{str}\left(Q\right)=K_0^{str}\left(\partial Q\right)$,\\
iv) $K_1^{str}\left(Q\right)$ consists of \\
- all 1-simplices $e\in K\left(Q\right)$ with $\partial e\in K_0^{str}\left(Q\right)$, and\\
- exactly one 1-simplex for each nontrivial homotopy class (rel.\ boundary) of loops $e$ with $\partial_0e=\partial_1e\in K_0^{str}\left(Q\right)$,\\
- the constant loop for the homotopy class of the constant loop at $x$, if $x\in K_0^{str}\left(Q\right)$,\\
v) for $n\ge 2$, if $\sigma\in S_n\left(Q\right)$ is an n-simplex with $\partial\sigma\in K^{str}_{n-1}
\left(Q\right)$,
then $\sigma$ is homotopic rel.\ boundary to a unique $\tau\in K_n^{str}\left(Q\right)$,\\
vi) if $\sigma\in K_n^{str}\left(Q\right)$ is homotopic rel.\ boundary to some $\tau\in K_n\left(Q\right)$, then $\sigma=\tau$,\\
vii) 
%for $n\ge 2$, 
if $\sigma\in K_n^{str}\left(Q\right)$ is homotopic rel.\ boundary to a simplex 
%$\tau\in K_n\left(\partial_0 Q
%\right)$ resp.\ 
$\tau\in S_n\left(\partial_1 Q\right)$, then
%\footnotemark[4]\footnotetext[4]{If $\tau$ is 
%at the same time homotopic rel.\ boundary to a simplex $\tau\in K_n\left(\partial_0 Q
%\right)$ and to a simplex $\tau\in K_n\left(\partial_1 Q\right)$, then it is actually homotopic to a simplex
%$\tau\in K_n\left(\partial_0 Q\right)\cap K_n^{str}\left(\partial_1Q\right)$. Indeed, the 
$\sigma\in K_n^{str}\left(\partial_1 Q\right)$; if $\sigma\in K_1^{str}\left(Q\right)$ is homotopic rel.\ boundary to a simplex $\tau\in S_1\left(\partial_0Q
%\cap \partial_1 Q
\right)$, then $\sigma\in K_1^{str}\left(\partial_0 Q
%\cap \partial_1 Q
\right)$,\\
%viii) for $n\ge 2$, if 
%$\sigma\in K_n^{str}\left(Q\right)$ is homotopic rel.\ boundary to a simplex 
%$\tau\in S_n\left(\partial_0 Q\right)$, then $\sigma\in K_n^{str}\left(\partial_0 Q\right)$,\\
viii) $K^{str}_*\left(Q\right)$
is aspherical, i.e.\ if $\sigma,\tau\in K^{str}_*\left(Q\right)$ have the same 1-skeleton, then $\sigma=\tau$.

\end{lem}

\begin{proof}
$K_*^{str}\left(Q\right)$ is
defined by induction on the dimension of simplices as follows.\\
\\
Definition of $\stk$:\\
%If $C$ is a path-component of $\partial Q$ with $C\cap \partial_0Q=\emptyset$, then
%choose exactly one point in $C$.\\
%Moreover choose each point in $\partial_0Q\cup q\left(\partial_0Q\right)$.\\
Choose $K_0^{str}\left(Q\right)$ such 
that conditions i),ii),iii) are satisfied.
Note that we have chosen a nonempty set of 0-simplices since we are 
assuming $\partial_1Q\not=\emptyset$. \\
%Moreover, by construction, 
%$K_0^{str}\left(Q\right)=K_0^{str}\left(\partial Q\right)$. \\
\\
Definition of $\suk$:\\
For an ordered pair $$\left(x,y\right)\in \stk
\times \stk$$ with $x\not=y$, there exists 
in each homotopy class (rel.\ boundary) of arcs $e$ with $$e\left(0\right)=x, e\left(1\right)=y$$ a unique simplex in $K_1\left(Q\right)$. Choose 
these 1-simplices to belong to $K_1^{str}\left(Q\right)$. 
(Uniqueness implies that vi)
) is true for $n=1$.)
%(By the construction in section 3.2. these simplices belong to $\partial Q$ resp.\ $\partial_1Q$ if they are homotopic into $\partial Q$ resp.\ $\partial_1Q$.)\\
Moreover, for pairs $$\left(x,x\right)\in \stk
\times \stk$$ choose one simplex in each homotopy class (rel.\ boundary) of loops $e$ with $$e\left(0\right)=e\left(1\right)=x.$$
For the homotopy class of the constant loop choose the constant loop. \\
Choose the 1-simplices in 
$\partial_0 Q$ and/or 
$\partial_1Q$  
whenever this is possible. 
(If a 1-simplex is homotopic into both $\partial_0Q$ and $\partial_1Q$, 
then it is necessarily homotopic into $\partial_0Q\cap\partial_1Q$. Indeed, 
a disk realizing a homotopy between 1-simplices in $\partial_0Q$ and $\partial_1Q$ can be made transversal to $\partial_0Q\cap\partial_1Q$ and then intersects $\partial_0Q\cap\partial_1Q$ in an arc resp.\ loop.) 
Hence vii) is satisfied for $n=1$.\\
%(Note that these choices do not contradict each other. Indeed, since 
%$S_0^{str}\left(\partial_0Q\right)\cap S_0^{str}\left(\partial_1Q\right)=\emptyset$, a 
%straight 1-simplex can not be homotopic rel.\ boundary
%into both $\partial_0Q$ and $\partial_1Q$. Moreover, if a 1-simplex,
%with vertices in $S_0^{str}\left(Q\right)\subset S_0\left(\partial Q\right)$, is homotopic to a degenerate 1-simplex, then this degenerate 1-simplex necessarily 
%has image in $\partial Q$. \\
\\
Definition of $K_n^{str}\left(Q\right)$ for $n\ge 2$, assuming that $K_{n-1}^{str}\left(Q\right)$ is defined:\\
For an $n+1$-tuple $\kappa_0,\ldots,\kappa_n$
of $n-1$-simplices in $K_{n-1}^{str}\left(Q\right)$, satisfying $$\partial_i\kappa_j=\partial_{j-1}\kappa_i$$ 
for all $i,j$, there are two possibilities:\\
- if no edge of any $\kappa_i$ is a loop, then, by asphericity of $Q$, there is 
a unique n-simplex
$$\sigma\in K_n\left(Q\right)$$
with $$\partial_i\sigma=\kappa_i$$ 
for $i=0,\ldots,n$. In this case set $\kappa:=\sigma$. Uniqueness implies that vi) is satisfied for $n$. (By the construction in Section \ref{sec:construction} 
%$\kappa\in K_n\left(\partial_0 Q\right)$ resp.\ 
$\kappa\in K_n\left(\partial_1Q\right)$ if $\kappa$ is homotopic rel.\ boundary into $\partial_1Q$.)\\
\\
- otherwise, choose 
an $n$-simplex $$\kappa\in S_n\left(Q\right)$$  with $$\partial_i\kappa=\kappa_i$$ 
for $i=0,\ldots,n$. By asphericity of $Q$, $\kappa$ exists
and is unique up to homotopy rel.\ boundary. Choose the simplices in 
%$\partial_0 Q$ and/or 
$\partial_1Q$ 
whenever this is possible.\\
%Hence vii) is satisfied for $n$.\\
\\
By construction, $K_*^{str}\left(Q\right)$ 
is closed under face maps and satisfies the conditions i)-vii). 
%Since $K\left(Q\right)$ is aspherical,
Condition viii) follows by induction on the dimension of subsimplices of $\sigma$ and $\tau$ from condition v). 
\end{proof}

The simplices in $K_*^{str}\left(Q\right)$
will be called the {\bf straight simplices}. \\
We remark that $K_*^{str}\left(Q\right)$ is 
not a multicomplex because simplices in $K_*^{str}\left(Q\right)$ need not have pairwise distinct vertices. (Note also that
simplices in $K\left(Q\right)$ belong to $K^{str}\left(Q\right)$ if and only if all their vertices belong to
$K_0^{str}\left(Q\right)$, by construction.)
%We emphasize that for $e=e_1*e_2*e_3\in
%\widehat{K_1^{str}}\left(Q\right)-K_1^{str}\left(Q\right)$ the decomposition $e_1*e_2*e_3$ can be recovered from $e$.\\

\begin{obs}\label{Obs5} Let $Q,\partial Q,\partial_1 Q,\partial_0Q$ satisfy Assumption I.  
%Let $q:Q\rightarrow Q$ be given by Construction 1 and let $A:=q\left(\partial_0 Q\right)$.
Let $K_*^{str}\left(Q\right)
\subset  S_*\left(Q\right)$ 
satisfy the conditions i)-viii) from  \hyperref[Lemma11]{Lemma \ref*{Lemma11}}.
Then \\
%a) the action of
%$$G:=\Pi\left(K\left(\partial_0Q\right)\right)$$ 
%on $K\left(Q\right)$ as defined in section 3.3.\ can be extended to 
%an action of $G$ on 
%$K_*^{str}\left(Q\right)$,\\
%b) 
$q:Q\rightarrow Q$ induces a simplicial map $q:K^{str}\left(Q\right)\rightarrow K^{str}\left(Q\right)$, compatible with the simplicial map $q:K\left(Q\right)\rightarrow K\left(Q\right)$ from Section \ref{sec:connected}.
\end{obs}
\begin{proof} 
%a) It is straightforward to check that $K_0^{str}\left(Q\right)$ is invariant under the $G$-action. Since $\sigma\in K\left(Q\right)$ belongs to $K^{str}\left(Q\right)$ if and only if all vertices of $\sigma$ belong to $K_0^{str}\left(Q\right)$, this implies that $K_*^{str}\left(Q\right)\cap K\left(Q\right)$ is invariant under the action of $G$.
%It remains to extend the action to all of $K_*^{str}\left(Q\right)$. (Note that a simplex in $K_*^{str}\left(Q\right)$ does {\em not} 
%belong to $K\left(Q\right)$ if and only if one of its edges is a loop.) By asphericity it suffices to 
%extend the action on the 1-skeleton, i.e.\ on loops. This is done in the obvious way.
%Let $e\in K_1^{str}\left(Q\right)$ be a loop with $e\left(0\right)=e\left(1\right)=x$ 
%and let $g=\left\{\gamma_1,\gamma_2,\ldots\right\}\in G$. If $x\not\in\left\{\gamma_1\left(0\right),\ldots,\gamma_m\left(0\right)\right\}$, then $ge=e$. Else, if $x=\gamma_i\left(0\right)$ for some $i$, then set 
%$$ge:=\left[\overline{\gamma}_i*e*\gamma_i\right],$$
%where as usual $\left[.\right]$ means the unique element of $K_1^{str}
%\left(Q\right)$ in the homotopy class rel.\ boundary. The same argument as in Section 3.3. shows that we have defined an action of $G$ on $K_*^{str}\left(Q\right)$.\\
%b) 
By construction, $q$ maps $K_0^{str}\left(Q\right)$ to itself. Indeed:\\ 
- if $C$ is a path-component of $\partial_0Q$ with $C\cap\partial_1Q=\emptyset$, then $q\left(v\right)=v$ for each $v\in C$,\\
- if $F$ is a path-component $F$ of $\partial_1Q$ with $F\cap\partial_0Q=\emptyset$, then $q\left(v\right)=v$ for {\em each} $v\in F$ (in particular for the unique $v\in F\cap K_0^{str}\left(Q\right)$),\\
- if $F$ is a path-component of $\partial_1Q$ with $F\cap\partial_0Q\not=\emptyset$, then we have
$K_0^{str}\left(Q\right)\cap F=\left\{x_{E_0^F},\ldots,x_{E_s^F}\right\}$, and $q\left(x_{E_i^F}\right)=x_{E_0^F}$ for $i=0,\ldots,s$ by Construction 1.

Hence $q$ induces a simplicial map on $K^{str}\left(Q\right)$. (The simplicial map 
agrees with $q$ on the 0-skeleton, and it maps each 1-simplex $e\in K_1^{str}\left(Q\right)$ to the unique 1-simplex of $K_1^{str}\left(Q\right)$ that is in the homotopy class rel.\ $\left\{0,1\right\}$ of
$q\left(e\right)$. Since $K^{str}\left(Q\right)$ is aspherical, this determines the simplicial map $q$ uniquely.) 
\end{proof}

\subsection{Definition of the straightening}\label{sec:straighten}

\begin{df}\label{Def8} Let $\left(Q,\partial_1 Q\right)$ be a pair
of topological spaces and 
let $z=\sum_{i\in I} a_i\tau_i\in C_n^{inf}\left(Q\right)$ a (possibly infinite)
singular chain.\\ 
a) A set of cancellations of $z$ is a symmetric set ${\mathcal{C}}
\subset S_{n-1}\left(Q\right)\times S_{n-1}\left(Q\right)$
with $\left(\eta_1,\eta_2\right)\in{\mathcal{C}}\Rightarrow \eta_1=\eta_2$ and $\eta_1=\partial_{k}\tau_{i_1},
\eta_2=\partial_{l}\tau_{i_2}$ for some $i_1,i_2\in I, k,l\in \left\{0,\ldots,n\right\}$.\\
b) If $z =\sum_{i\in I} a_i\tau_i\in C_n^{inf}\left(Q\right)$ and $\mathcal{C}$ is a set of cancellations for $z$, then the associated 
simplicial set $\Upsilon_{z,{\mathcal{C}}}$ is the simplicial set generated\footnotemark\footnotetext[7]{That is, the subset of $S_*^{sing}\left(Q\right)$ which contains the $\mid I\mid$ n-simplices $\Delta_i,
i\in I$, together with all simplices obtained by iterated applications of face and degeneracy operators, cf.\ \cite{may}, Example 1.5.}
by $\left\{\Delta_i: i\in I\right\}$,  subject to the identifications                                   
$\partial_{k}\Delta_{i_1}=\partial_{l}\Delta_{i_2}$ if and only if $\left(\partial_{k}\tau_{i_1},\partial_{l}\tau_{i_2}\right)\in{\mathcal{C}}$.\\
c) Let $z =\sum_{i\in I} a_i\tau_i\in C_n^{inf}\left(Q\right)$.  
Choose a minimal presentation for $\partial z$ (i.e.\ no further cancellation is possible). Let 
$$J=J_{\partial z}:=\left\{\begin{array}{c}\left(i,a\right)\in I\times\left\{0,\ldots,n\right\}:\\
\partial_a\tau_i\mbox{\ occurs\ with\ non-zero\ coefficient\ in\ the\ chosen\ presentation\ of\ }\partial z\end{array}\right\}.$$ 
Let $\mathcal{C}$ be a set of cancellations for $z$. 
Then the simplicial set $\partial\Upsilon_{z,{\mathcal{C}}}$ is defined as the set consisting of $\mid J\mid$ 
n-1-simplices $\Delta_{i,a},
\left(i,a\right)\in J$, together with all their iterated faces and degenerations, subject to the identifications
$\partial_a\partial_{a_1}\tau_{i_1}=\partial_a\partial_{a_2}\tau_{i_2}$ for all $a=0,\ldots,n-1$, 
whenever $\left(\partial_{a_1}\tau_{i_1},\partial_{a_2}\tau_{i_2}\right)\in{\mathcal{C}}$
and $ \left(i_1,a_1\right)\in J$ .\\
d) If $z =\sum_{i\in I} a_i\tau_i\in C_n^{inf}\left(Q\right)$ is a relative cycle, then a set of 
cancellations $\mathcal{C}$ is called sufficient if the formal sum
$\sum_{i\in I}\sum_{k=0}^n \left(-1\right)^k a_i\partial_k\tau_i$ can be reduced to a chain in 
$C_{n-1}^{inf}\left(\partial Q\right)$ by substracting (possibly infinitely many) multiples of $\left(\partial_{a_1}\tau_{i_1}-\partial_{a_2}\tau_{i_2}\right)$ with $\left(\partial_{a_1}\tau_{i_1},\partial_{a_2}\tau_{i_2}\right)\in{\mathcal{C}}$.
\end{df}
\begin{obs}\label{Obs6} Let $\left(Q,\partial_1 Q\right)$ be a pair of topological spaces.\\
a) If $z=\sum_{i\in I} a_i\tau_i\in C_n^{inf}\left(Q\right)$ is a singular chain, $\mathcal{C}$ a set of cancellations, and $\Upsilon:=\Upsilon_{z,{\mathcal{C}}}$ the associated simplicial set, 
then the geometric realisation $\mid\Upsilon\mid$ is obtained from $\mid I\mid $ copies of the standard n-simplex $\Delta_i,i\in I$, with identifications $\partial_{a_1}\Delta_{i_1}=\partial_{a_2}\Delta_{i_2}$ if and only if 
$\left(\partial_{a_1}\tau_{i_1},\partial_{a_2}\tau_{i_2}\right)\in{\mathcal{C}}$. 
Moreover, for a minimal presentation of $\partial z$ and $\partial\Upsilon:=\partial\Upsilon_{z,{\mathcal{C}}}$, $\mid\partial
\Upsilon\mid$ is the subspace of $\mid\Upsilon\mid$ containing all simplices $\partial_{a_1}\Delta_{i_1}$ with
$\left(i_1,a_1\right)\in J$. \\
b) There exists an associated continuous map $\tau:\mid\Upsilon\mid\rightarrow Q$ with $\tau\mid\Delta_i=\tau_i$ (upon the identification $\Delta_i=\Delta^n$). If $z$ is a relative cycle, i.e.\ 
$ \partial z\in C_{n-1}^{inf}\left(\partial_1 Q\right)$,
then $\tau$ maps $\mid\partial\Upsilon\mid$ to $\partial_1 Q$.\\
c) Let $z_1=\sum_{i\in I}a_i\tau_i,
z_2=\sum_{i\in I}a_i\sigma_i\in C_n^{inf}\left(Q,\partial_1 Q\right)$ be relative cycles and 
${\mathcal{C}}_1,{\mathcal{C}}_2$ sufficient sets of cancellations of $z_1$ resp.\ $z_2$.
Assume that $\left(\partial_{a_1}\tau_{i_1},\partial_{a_2}\tau_{i_2}\right)\in{\mathcal{C}}_1$ if and only if
$\left(\partial_{a_1}\sigma_{i_1},\partial_{a_2}\sigma_{i_2}\right)\in{\mathcal{C}}_2$, and that there exist minimal presentations of $\partial z_1,\partial z_2$ such that $J_{z_1}=J_{z_2}$. \\
If the associated continuous maps \\
$\tau,\sigma:\mid\Upsilon\mid\rightarrow Q$ are homotopic, \\
for a homotopy mapping $\mid\partial \Upsilon\mid$ to $\partial Q$, then $\sum_{i\in I}a_i\tau_i$
and $
\sum_{i\in I}a_i\sigma_i\in C_*^{inf}\left(Q,\partial Q\right)$ are relatively homologous.
\end{obs}
We emphasize that we do not assume that $\mathcal{C}$ is a complete list of cancellations, 
the simplicial map $\tau_*:C_*^{simp}\left(\Upsilon\right)\rightarrow C_*^{sing}\left(Q\right)$ need not be injective.\\
\\
%If a simplicial set $\widehat{S_*^{str}}\left(Q\right)\subset S_*\left(Q\right)$ of straight simplices is defined according
%to Lemma 11, we will denote the associated chain complex
%$$C_*^{str,inf}\left(Q\right):=C_*^{simp,inf}\left(\widehat{S_*^{str}}\left(Q\right)\right)\mbox{\ \ and\ \ }C_*^{str,inf}\left(Q,\partial Q\right):=C_*^{str,inf}\left(Q\right)/
%C_*^{str,inf}\left(\partial Q\right).$$
After having set up the necessary notations, we now start with the actual definition of the straightening.
We first mention that there is of course an analogue of the classical straightening (\cite{bp}, Lemma C.4.3.) in our setting.

\begin{obs}\label{Obs7}  Let $Q,\partial Q,\partial_1 Q,\partial_0Q$ satisfy Assumption I.
%, and let $G:=\Pi\left(K\left(\partial_0Q\right)\right)$.  
%Let $q:Q\rightarrow Q$ be given by  \hyperref[Cons1]{Construction \ref*{Cons1}} and let $A:=q\left(\partial_0 Q\right)$.
Let $K_*^{str}\left(Q\right)
\subset  S_*\left(Q\right)$
satisfy the conditions i)-viii) from  \hyperref[Lemma11]{Lemma \ref*{Lemma11}}.

Then there exists a 'canonical straightening' map
$$str_{can}:C_*^{simp,inf}\left(K\left(Q\right)\right)\rightarrow
C_*^{simp,inf}\left(K^{str}\left(Q\right)\right),$$
mapping $C_*^{simp,inf}\left(K\left(\partial_1 Q\right)\right)$
to $C_*^{simp,inf}\left(K^{str}\left(\partial_1 Q\right)\right)$, with the following properties:\\
i) $str_{can}$ is a chain map,\\
%ii) if $\tau\in K^{str}\left(Q\right)$, then $str_{can}\left(\tau\right)=\tau$,\\
ii) if $z=\sum_{i\in I} a_i\tau_i\in C_*^{simp,inf}\left(K\left(Q\right)\right)$ and $\sum_{i\in I}a_i\sigma_i:=\sum_{i\in I}
a_i str_{can}\left(\tau_i\right)$, then
the maps $$
\tau,\sigma:
%\left(
\mid\Upsilon\mid
%,\mid\partial\Upsilon\mid
%\right)
\rightarrow 
%\left(
Q
%,\partial Q\right)
$$ 
(defined by  \hyperref[Obs6]{Observation \ref*{Obs6}}b) after fixing 
a set of cancellations $\mathcal{C}$ and a minimal presentation of $\partial z$) are homotopic.

% as maps of pairs. 
Moreover,
if $z=\sum_{i\in I} a_i\tau_i$ is a relative cycle with 
$\partial z\in C_*^{simp,inf}
\left(K\left(\partial_1Q\right)\right)$, then the same is true for $\sum_{i\in I}a_i\sigma_i$ and 
$$
\tau,\sigma:
\left(
\mid\Upsilon\mid
,\mid\partial\Upsilon\mid
\right)\rightarrow 
\left(
Q
,\partial_1 Q\right)
$$
are homotopic as maps of pairs. 

In particular, 
$\sum_{i\in I} a_i str_{can}\left(\tau_i\right)$ is relatively homologous 
to $\sum_{i\in I} a_i\tau_i$,\\
%iii) if $v\in S_0\left(\partial_1 Q\right)$, then the homotopy $H\left(v\right)$ between $v$ and $str_{can}\left(v\right)$
%is either constant or $H\left(v\right)\in K_1\left(\partial_1 Q\right)$.
\end{obs}
\begin{proof}
We define $str_{can}$, and the homotopy to the identity,
by induction on the dimension of simplices. (During the construction we take care that $str_{can}$ and 
the homotopy preserve $K\left(\partial_1Q\right)$.)\\
% and $q\left(K\left(\partial_0Q\right)\right)$.) \\

{\em 0-simplices.}\\
If $C$ is a path-component of $\partial_0Q$ with $C \cap\partial_1Q=\emptyset$, then we define $str_{can}\left(v\right)=v$ for each 0-simplex $v$ in $C$. 
The homotopy $H\left(v\right)$ is for each $v$ given by the constant map.

If $C$ is a path-component of $\partial_0Q$ with $C \cap\partial_1Q
\not=\emptyset$, then there is at least one 
path-component $F$ of $\partial_1Q$ with $C\cap F\not=\emptyset$. By \hyperref[Cons1]{Construction \ref*{Cons1}} and condition ii) from \hyperref[Lemma11]{Lemma \ref*{Lemma11}}, for each such $F$, there 
is a straight 0-simplex $x_{E_i^F}\in C\cap F$. Choose one such straight 0-simplex (among the $x_{E_i^F}$'s)
for each path-component $C$ of $\partial_0Q$, denote it $x_C$,
and for each $v\in C$ we define 
$str_{can}\left(v\right):=x_C\in K_0^{str}\left(Q\right)\cap C$ and we choose the homotopy $H(v)$ to belong to $C$.

If $v\in\partial_1Q$,
then there is (at least) one straight 0-simplex in the same path-component $F$ of $\partial_1Q$, we choose $str_{can}\left(v\right)\in F\cap K_0^{str}\left(Q\right)$ and
%-\partial_0Q$ belongs to a path-component $C$ of $\partial Q$ with $C\cap \partial_0Q\not=\emptyset$,
%then we define $str_{can}\left(v\right)$ to be some 0-simplex in $C\cap \partial_0Q\cap\partial_1 Q$ 
there exists $H\left(v\right)\in K_1\left(\partial_1 Q\right)$ with
$\partial H\left(v\right)=v-str_{can}\left(v\right)$. 
%(That means we choose $str_{can}\left(v\right)$ in \partial_0Q\cap \partial_1Q$ such that it belongs to the same path-component of $\partial_1Q$ as $v$, and we choose the homotopy in $\partial_1Q$.)
%If $v\in\partial_1Q$ belongs to a path-component $C$ of $\partial Q$ with $C\cap \partial_0Q=\emptyset$, then   
%we define $str_{can}\left(v\right)$ to be the (unique) straight 0-simplex in $C\subset\partial_1Q$ and we fix arbitrarily some $H\left(v\right)\in K_1\left(\partial_1 Q\right)$ with
%$\partial H\left(v\right)=v-str_{can}\left(v\right)$.

If $v\not\in\partial Q$, then we define $str_{can}\left(v\right)$ to be some straight 0-simplex in $\partial Q $ and we fix arbitrarily some $H\left(v\right)\in K_1\left(Q\right)$ with
$\partial H\left(v\right)=v-str_{can}\left(v\right)$.\\

{\em 1-simplices.}\\
For $e\in K_1\left(Q\right)$ we define $$str_{can}\left(e\right):=\left[\overline{H\left(\partial_1e\right)}*e*H\left(\partial_0e\right)\right],$$
where, as always, [.] denotes the unique 1-simplex in $K_1^{str}\left(Q\right)$, which is homotopic rel.\ boundary to the path in the brackets.

$e$ is homotopic to $str_{can}\left(e\right)$ 
by the canonical homotopy which is inverse to the homotopy moving $\overline{H\left(\partial_1e\right)}$ resp.\  
$H\left(\partial_0e\right)$ into constant maps. In particular, the restriction of this homotopy to $\partial_1e,\partial_0e$ gives $H\left(\partial_1e\right),\overline{H\left(\partial_0e\right)}$. Thus, for different
edges with common vertices, the homotopies are compatible. We thus have constructed a homotopy for the 1-skeleton
$\Upsilon_1$. 

We note that, for $v\in \partial_1Q$, the homotopy $H\left(v\right)$ is either constant or
$H\left(v\right)\in K_1\left(\partial_1 Q\right)$,
Thus if $\tau\in K_1\left(\partial_1Q\right)$  then  
$str_{can}\left(\tau\right)\in K_1^{str}\left(\partial_1Q\right)$ and
the homotopy between $\tau$ and $str_{can}\left(\tau\right)$ takes place in $\partial_1 Q$. \\

{\em n-simplices.}\\
We assume inductively, that for some $n\ge 1$, we have 
defined $str_{can}$ on $K_{*\le n}\left(Q\right)$, 
mapping $K_{*\le n}\left(\partial_1 Q\right)$ to $K_{*\le n}^{str}
\left(\partial_1 Q\right)$, and satisfying i),ii),iii).

Let $\tau\in K\left(Q\right)$ be an n+1-simplex. Then we have by ii) a 
homotopy between $\partial\tau$ 
and $str_{can}\left(\partial\tau
\right)$. By  \hyperref[Obs1]{Observation \ref*{Obs1}} this homotopy extends to $\tau$. The resulting simplex $\tau^\prime$ satisfies
$\partial\tau^\prime\in K_n^{str}\left(Q\right)$.
Condition v) from  \hyperref[Lemma11]{Lemma \ref*{Lemma11}} means that
$\tau^\prime$ is homotopic rel.\ boundary to a unique simplex $str_{can}\left(\tau\right)\in 
K_{n+1}^{str}\left(Q\right)$. This proves the inductive step.

If $\tau\in K\left(\partial_1 Q\right)$, then we can inductively assume that the homotopy of $\partial\tau$ has image 
in $\partial_1 Q$. Then condition vii) from  \hyperref[Lemma11]{Lemma \ref*{Lemma11}} implies $str_{can}\left(\tau\right)\in K_{n+1}^{str}\left(\partial_1 Q\right)$. Moreover, since $\partial_1Q$ is
aspherical, the homotopy of $\tau$ can be chosen to have image in $\partial_1 Q$.\\

By construction, for any set of cancellations $\mathcal{C}$,
the induced maps $\tau$ and $\sigma$ 
are homotopic. In particular, if we chose a sufficient set of cancellations in 
the sense of  \hyperref[Def8]{Definition \ref*{Def8}}d), then  \hyperref[Obs6]{Observation \ref*{Obs6}}c) implies that $\sum_{i=1}^r a_i str_{can}\left(\tau_i\right)$ is (relatively) homologous
to $\sum_{i=1}^r a_i\tau_i$.

\end{proof}

However, we want to define a more refined straightening, 
which will be defined only on relative
cycles with some kind of additional information. \\

Before stating the definition of "distinguished 1-simplices" we remark that there is a left and right action 
of the pseudogroup $\Gamma:=\Omega\left(\partial Q\right)$ (as defined in Section \ref{sec:action}) on 
$K_1^{str}\left(Q\right)$: if $e\in K_1^{str}\left(Q\right), \gamma_1\in \pi_1\left(\partial Q,\partial_1e\right), 
\gamma_2\in \pi_1\left(\partial Q,\partial_0e\right)$, 
then let $\gamma_1e\gamma_2$ be the unique straight 1-simplex homotopic rel.\ $\left\{0,1\right\}$ 
to $\gamma_1*e*\gamma_2$. (The left action agrees with the action defined in Section \ref{sec:action}.) The cosets $\Gamma K_1^{str}\left(Q\right)\Gamma$ in  \hyperref[Def9]{Definition \ref*{Def9}} are with respect to this action.

For $x,y\in K_0^{str}\left(Q\right)$ we will denote $K_{1,xy}^{str}:=\left\{e\in K_1^{str}\left(Q\right): \partial_1e=x,\partial_0e=y\right\}$.

\begin{df}\label{Def9} Let $Q,\partial Q,\partial_1 Q,\partial_0Q$ satisfy Assumption I.  

Let $q:Q\rightarrow Q$  and $\left\{x_{E_i^F}\in \partial_0Q\cap\partial_1Q\right\}$ be given by  \hyperref[Cons1]{Construction \ref*{Cons1}}.
% and let $A:=q\left(\partial_0 Q\right)$.

Let $K_*^{str}\left(Q\right)\subset 
S_*\left(Q\right)$ satisfy conditions i)-viii) from  \hyperref[Lemma11]{Lemma \ref*{Lemma11}}.

%Let $$D_0:=\left\{\begin{array}{c}x\in K_0^{str}\left(Q\right): x=x_{E_i^F} \mbox{\ for\ some\ }F, E_i\subset F\cap\partial_0Q,
%\mbox{\ or\ x\ belongs\ to\ some}\\ 
%\mbox{path-component\ F\ of\ } 
%\partial_1 Q\mbox{\ with\ }F\cap\partial_0Q\not=\emptyset,\end{array}\right\}$$

A set $D\subset K_1^{str}\left(Q\right)$ is called a set of distinguished 1-simplices if \\
ix) $\partial_0e,\partial_1e\in K_0^{str}\left(Q\right)$
%K_0^{str}\left(Q\right)$ 
for each $e\in D$,\\
%if $\partial_je$ for $j\in\left\{0,1\right}$, then exists some $i$ with $\partial_je\in U_i$,\\
x) for each $$\left(x,y\right)\in
K_0^{str}\left(Q\right)\times K_0^{str}\left(Q\right)$$ we have that $$D_{xy}:=\left\{e\in D:\partial_1e=x,\partial_0=y\right\}$$ contains exactly one element in
each double coset (w.r.t.\ $\Gamma=\Omega\left(\partial Q\right)$)
$$\Gamma f\Gamma\in \Gamma K_{1,xy}^{str}\left(Q\right)\Gamma, $$ 
%there exists exactly one 1-simplex $e\in D$ with
%$$\partial_1x,\partial_0=y, 
%e\in \Gamma f\Gamma.$$
xi) for all $x\in K_0^{str}\left(Q\right)$, the constant loop $c_x$ belongs to $D$,\\
xii) if $e\in D$, then $\overline{e}\in D$, where $\overline{e}$ denotes the 1-simplex with the opposite orientation,\\
%for all $\left(x,y\right)\in
%K_0^{str}\left( Q\right)\times K_0^{str}\left( Q\right)$ we have $e_{xy}\in D=\overline{e_{yx}}$, i.e.\ $e_{xy}$ and $e_{yx}$ agree up to orientation,\\
xiii) if $F,F^\prime$ are path-components of $\partial_1Q$ and $\left\{x_{E_i^F}\in\partial_0Q\cap F\right\}, 
\left\{x_{E_j^{F^\prime}}\in\partial_0Q\cap F^\prime\right\}$ 
are given\footnotemark\footnotetext[8]{If $F\cap\partial_0Q=\emptyset$ and/or $F^\prime\cap\partial_0Q=\emptyset$, then
there is only one straight 0-simplex $x_{E_0^F}$ resp.\ $x_{E_0^{F^\prime}}$ in $F$ resp.\ $F^\prime$. In particular, if $F\cap\partial_0Q=\emptyset$ 
and $F^\prime\cap\partial_0Q=\emptyset$, then condition xiii) is empty. }
%Similarly, condition xiv) is empty if $C_1\cap\partial_0Q=\emptyset$ and $C_2\cap\partial_0Q=\emptyset$.} 
by  \hyperref[Cons1]{Construction \ref*{Cons1}}, then $q\left(D_{x_{E_i^F}x_{E_j^{F^\prime}}}\right)=D_{x_{E_0^F}x_{E_0^{F^\prime}}}$ for all $x_{E_i^F}, x_{E_j^{F^\prime}}$,\\
xiv) if 
%$e_1,e_2\in D$ with
$x_1,x_2\in C_1, y_1,y_2\in C_2$ for some path-components $C_1,C_2$ of $\partial Q$, then 
for each $e_1\in D_{x_1y_1}$ exists
some $e_2\in D_{x_2y_2}$ with
$q\left(e_2\right)=g q\left(e_1\right)$ for some $g\in H:=q_*\left(\Pi\left(K\left(\partial_0Q\right)\right)\right)$.\\
%xv) if $e_1,e_2\in D$ with
%$\partial_1e_1,\partial_1e_2\in C_1$ for some path-component $C_1$ of $\partial_0 Q$ and 
%$\partial_0e_1=\partial_0e_2\in\partial_1Q$, 
%then $e_2=g e_1$ for some $g\in G:=\Pi\left(K\left(\partial_0Q\right)\right)$.
\end{df}

%Remark: Condition xiii) applies to vertices in the same component of $\partial_1Q$ 
%and will be crucial for the proof of Theorem 1. Condition xiv) applies to vertices 
%that are in the same component of $\partial Q$ but not necessarily in 
%the same component of $\partial_1Q$. Conditon xiv)
%will "only"
%be needed for proving the well-definedness of $q\circ str$ in Corollary 3.

\begin{obs}\label{Obs8} Let the assumptions of \hyperref[Def9]{Definition \ref*{Lemma9}} be satisfied.
Then a set $D$ of distinguished 1-simplices exists.\end{obs}
\begin{proof}

For each path-component $C$ of $\partial Q$ we fix some $x_C\in K_0^{str}\left(C\right)$. 

For each pair $\left\{C_1,C_2\right\}$ of path-components 
we fix one simplex $e$ with $$\partial_1e=x_{C_1},\partial_0e=x_{C_2}$$ 
in each coset of $\Gamma K_{1,x_{C_1}x_{C_2}}^{str}\left(Q\right)\Gamma$ to belong to $D_{x_{C_1}x_{c_2}}$. \\
(For all chosen 1-simplices $e\in D_{x_{C_1}x_{C_2}}$, we 
 choose $\overline{e}$ to belong to $D_{x_{C_2}x_{C_1}}$. If $C_1=C_2$, then in particular for the coset of the constant loop we choose the constant loop to belong to $D_{x_{C_1}x_{C_2}}$.)\\

For each path-component $C$ of $\partial Q$ and each path-component $F$ of $C\cap\partial_1Q$, 
we have that 
$q\left(x_C\right)$ and $q\left(x_{E_0^F}\right)$
belong to the path-connected set $q\left(\partial_0Q\cap C\right)$. Therefore we have
a sequence of 1-simplices $\alpha_1,\ldots,\alpha_m\in K_1\left(\partial_0Q\right)$ with images in distinct path-components
of $\partial_0Q\cap C$,
such that $$\partial_1 q\left(\alpha_1\right)=q\left(x_{C}\right), \partial_0q\left(\alpha_1\right)=\partial_1 q\left(\alpha_2\right),\ldots,\partial_0q\left(\alpha_{m-1}\right)=
\partial_1q\left(\alpha_m\right),\partial_0q\left(\alpha_m\right)=q\left(x_{E_0^F}\right).$$

In order to prepare the definition of the $D_{x,y}$'s, we first describe, for each $x\in C\cap K_0^{str}\left(Q\right)$ a sequence $\left\{\alpha_1,\ldots,\alpha_k\right\}$ of 1-simplices:\\
- if $C\cap \partial_1Q=\emptyset$, then $k=1$ and for each $x\in C$ we choose arbitrarily a 1-simplex $\alpha_1$ in $C$ with $\partial_1\alpha_1=x_C, \partial_0\alpha_1=x$,\\
- if $C\cap\partial_0Q=\emptyset$, then $C\cap K_0^{str}\left(Q\right) =\left\{x_C\right\}$ by  \hyperref[Lemma11]{Lemma \ref*{Lemma11}}, condition ii), and we let $k=0$,\\
- if $C\cap \partial_0Q\cap\partial_1Q\not=\emptyset$, then by condition ii) from  \hyperref[Lemma11]{Lemma \ref*{Lemma11}} we have $x=x_{E_i^F}$ for some path-component $F$ of $\partial_1Q$ and some $i$, thus we have the above-constructed sequence $\alpha_1,\ldots,\alpha_m$ with
$\partial_1 q\left(\alpha_1\right)=q\left(x_{C}\right), \partial_0q\left(\alpha_1\right)=\partial_1 q\left(\alpha_2\right),\ldots,\partial_0q\left(\alpha_{m-1}\right)=
\partial_1q\left(\alpha_m\right),\partial_0q
\left(\alpha_m\right)=q\left(x_{E_i^F}\right)$, where the last equality holds true because $q\left(x_{E_i^F}\right)=x_{E_0^F}=q\left(x_{E_0^F}\right)$.
%and we append $\alpha_{m+1}:=l_{E_0^F,E_i^F}$, oriented such that
%\partial_1\alpha_{m+1}=l_{E_0^F}, \partial_0\alpha_{m+1}=l_{E_i^F}$.\\

Let $x,y\in K^{str}_0\left(Q\right)$. Let $C_1,C_2$ be the path-components of $\partial Q$ with $x\in C_1,y\in C_2$. 

We have constructed
sequences of 1-simplices $\alpha_1,\ldots,\alpha_k\in K_1\left(\partial Q\right)$ resp. $\beta_1,\ldots,\beta_l\in K_1\left(\partial Q\right)$,   
such that $\partial_1 q\left(\alpha_1\right)=q\left(x_{C_1}\right), \partial_0q\left(\alpha_1\right)=\partial_1 q\left(\alpha_2\right),\ldots,\partial_0q\left(\alpha_{k-1}\right)=
\partial_1q\left(\alpha_k\right),\partial_0q\left(\alpha_k\right)=q\left(x\right)$ resp.\ $\partial_1 q\left(\beta_1\right)=q\left(x_{C_2}\right), \partial_0q\left(\beta_1\right)=\partial_1 q\left(\beta_2\right),\ldots,\partial_0q\left(\beta_{k-1}\right)=
\partial_1q\left(\beta_k\right),\partial_0q\left(\beta_k\right)=q\left(y\right)$. Note that all $q\left(\alpha_i\right)$ and $q\left(\beta_i\right)$ 
are either constant or
contained in $q\left(K_1\left(\partial_0Q\right)\right)$.

%{\bf If, for some path-component $F$ of $\partial_1Q$, we consider the points $x_{E_0^F},\ldots,x_{E_s^F}$ 
%given by Construction 1, and if $\alpha_1,\ldots,\alpha_m$ is the sequence of 1-simplices connecting $x_{C_1}$ to $x_{E_0^F}$, 
%then we can and will {\bf choose $\alpha_1,\ldots,\alpha_m,l_{E_0^FE_i^F}$ as the sequence of 1-simplices connecting $x_{C_1}$ to $x_{E_i^F}$}, for $i=1,\ldots,s$. }

%exists a sequence of 1-simplices $\beta_1,\ldots,\beta_n\in K_1\left(\partial_0Q\right)$ in distinct path-components of
%$\partial_0Q\cap C_2$ whose images under $q$ connect $q\left(x_{C_2}\right)$ to $q\left(y\right)$.
%As above, we can make the choice such that, for each path-component $F^\prime$ 
%of $\partial_1Q$, the sequence of 1-simplices connecting $x_{C_2}$ to $x_{E_j}^{F^\prime}$, for $j=1,\ldots,t$, 
%consists of the sequence connecting $x_{C_2}$ to $x_{E^{F^\prime}_0}$ appended by $l_{E^{F^\prime}_0E^{F^\prime}_j}$.

 Let $H:= q_*\left(\Pi\left(K\left(\partial_0Q\right)\right)\right)$. Define 
$$g:=\left\{q\left(\alpha_1\right),q\left(\overline{\alpha}_1\right)\right\}\ldots\left\{q\left(\alpha_k\right),q\left(\overline{\alpha}_k\right)\right\}\left\{q\left(\beta_l\right),q\left(\overline{\beta}_l\right)\right\}\ldots\left\{q\left(\beta_1\right),q\left(\overline{\beta}_1\right)\right\}\in H.$$
(If $k=l=0$, this means just $g=1$.) 

We have that $g=g^{-1}$ and that $$ge\in K_{1,q\left(x\right)q\left(y\right)}^{str}\left(Q\right) \Longleftrightarrow
e\in K_{1,q\left(x_{C_1}\right)q\left(x_{C_2}\right)}^{str}\left(Q\right).$$

%Note: if, for some path-components $F,F^\prime$ of $\partial_1Q$, we consider the points $x_{E_0^F},\ldots,x_{E_s^F}$ resp.\  $x_{E_0}^{F^\prime},\ldots,x_{E_t}^{F^\prime}$ 
%given by 
%Construction 1, then, we have $$q\left(l_{E_0^FE_i^F}\right)\equiv x_{E_0^F}, 
%q\left(l_{E^{F^\prime}_0E^{F^\prime}_j}\right)\equiv x_{E^{F^\prime}_0}$$
%for all $E_i^F,E_j^{F^\prime}$.
%This implies 
%that the $g$ associated
%to $x_{E_0^F},x_{E^{F^\prime}_0}$ is the same as the $g$ associated to $x_{E_i^F},x_{E_j^{F^\prime}}$.
% for $i=1,\ldots,s$ and $j=1,\ldots,t$. 

By construction, the {\bf $g$ associated to $x_{E_i^F}, x_{E_j^{F^\prime}}$ agrees with the $g$ associated to $x_{E_0^F}, x_{E_0^{F^\prime}}$}.\\
%in the setting of condition xiii), we have $$g=\left\{q\left(l_{E_0^FE_i^F}\right),q\left(\overline{l_{E_0^FE_i^F}}\right)\right\},\left\{
%$q\left(l_{E^{F^\prime}_0E^{F^\prime}_j}\right), q\left(\overline{l_{E^{F^\prime}_0E^{F^\prime}_j}}\right)\right\}=1.$$
%Thus condition xiii) is a consequence of condition xiv) and it suffices to prove xiv) for the chosen $g$.\\

We are given $D_{x_{C_1}x_{C_2}}$ and we want to define $D_{xy}$ such that condition xiii) is satisfied. 

First, if $C_1\cap\partial_1Q=\emptyset$ or $C_2\cap\partial_1Q=\emptyset$, then we can fix an arbitrary choice of $D_{x,y}$ 
satisfying conditions x),xi,xii). (Condition xiii) is empty in this case.)
%fix (arbitrarily) some paths $\gamma_1$ from $x$ to $x_{C_1}$ resp.\ $y$ to $x_{C_2}$, let $g

So let us assume $C_1\cap\partial_1Q\not
=\emptyset$ and $C_2\cap\partial_1Q\not=\emptyset$.
%By condition ii) from Lemma 11 this means 
We note that $$q:\left(Q,\partial Q,
\partial_1Q\right)\rightarrow\left(Q,\partial Q,\partial_1Q\right)$$ is homotopic to the identity as a map of triples, by the construction in Section 5.1. This implies that cosets of $\Gamma K^{str}_{1,xy}\left(Q\right)\Gamma$
are in 1-1-correspondence (by applying $q$) to cosets of $\Gamma
K^{str}_{1,q\left(x\right)q\left(y\right)}\Gamma$. It is thus sufficient to describe $q\left(D_{xy}\right)\subset K^{str}_{1,q\left(x\right)q\left(y\right)}$.

Let $$\Gamma f\Gamma\in \Gamma K_{1,q\left(x\right)q\left(y\right)}^{str}\left(Q\right)\Gamma$$ be a double coset. Then the double coset 
$$\Gamma \left(gf\right)\Gamma\in \Gamma  K_{1,q\left(x_{C_1}\right)q\left(x_{C_2}\right)}^{str}\left(Q\right)\Gamma$$
is the image under $q$ of some double coset $$\Gamma
e^\prime\Gamma\in\Gamma
K_{1,x_{C_1}x_{C_2}}^{str}\left(Q\right)\Gamma$$ 
Let $e$ be the unique distinguished simplex in 
$\Gamma
e^\prime\Gamma$. Then we choose $gq\left(e\right)$ to be the 
distinguished simplex in
$\Gamma f\Gamma$.
%we choose the 1-simplex $geg$, where $e$ is the unique distinguished 1-simplex in the double coset 
%$\Gamma \left(gfg\right)\Gamma\in \Gamma  K_{1,q\left(x_{C_1}\right)q\left(x_{C_2}\right)}^{str}\left(Q\right)\Gamma$. 
This is possible because $gq\left(e\right)$ belongs to the double coset $\Gamma f\Gamma$. Indeed 
$$q\left(e\right)\in \Gamma \left(gf\right)\Gamma$$ 
means that $q\left(e\right)=q_*\left(\gamma_1\right)gfq_*\left(\gamma_2\right)$ for some loops $\gamma_1$ and $\gamma_2$ based at $x_{C_1}$ 
resp.\ $x_{C_2}$, and this implies 
$gq\left(e^\prime\right)=q_*\left(\gamma_1^\prime\right) fq_*\left(\gamma_2^\prime\right)$ with 
$$\gamma_1^\prime:=\left[\overline{\alpha}_m*\ldots*\overline{\alpha}_1*\gamma_1*\alpha_1*\ldots*\alpha_m\right], \gamma_2^\prime:=
\left[\overline{\beta}_n*\ldots*\overline{\beta}_1*\gamma_2*\beta_1*\ldots*\beta_n\right].$$

This  defines $D_{xy}$. By construction, condition xiv) is satisfied if $e_1\in D_{x_{C_1}x_{C_2}}$. In general, if $e_1\in D_{x_1y_1}$, then we get $e\in D_{x_{C_1}x_{C_2}}$ and $g_1\in H$ with $q\left(e_1\right)=g_1q\left(e\right)$ and 
$e_2\in D_{x_2y_2}, g_2\in H$ with $q\left(e_2\right)=g_2q\left(e\right)$, thus $q\left(e_2\right)=g_2g_1^{-1}q\left(e_1\right)$.

Condition xiii) is implied because $q\left(x_{E_i^F}\right)=x_{E_0^F}, q\left(x_{E_j^{F^\prime}}\right)=x_{E_0^{F^\prime}}$ and the $g$ associated to $x_{E_i^F}, x_{E_j^{F^\prime}}$ agrees with the $g$ associated to $x_{E_0^F}, x_{E_0^{F^\prime}}$.\\
One checks easily that xi) and xii) are true for $D_{xy}$ since they are true for $D_{x_{C_1}x_{C_2}}$. 
\end{proof}

\begin{df}\label{Def10} 
Let $Q,\partial Q, \partial_0Q, \partial_1Q$ satisfy Assumption I.
Let $z=\sum_{i\in I} a_i\tau_i\in C_n^{inf}\left(Q\right)$
be a singular chain and $\Upsilon$ the associated simplicial set (for some set of cancellations $\mathcal{C}$). 

We say that a labeling of the elements of the 1-skeleton $\Upsilon_1$ by
0's and 1's is admissible, if $\partial e_1\cap\partial e_2=\emptyset$ for all 1-labeled vertices $e_1,e_2$.\end{df}

\begin{lem}\label{Lemma12} Let $Q,\partial Q,\partial_1 Q,\partial_0Q$ satisfy Assumption I.  Let $q:Q\rightarrow Q$ be given by  \hyperref[Cons1]{Construction \ref*{Cons1}}.
% and let $A:=q\left(\partial_0 Q\right)$.

Let $K_*^{str}\left(Q\right)\subset
S_*\left(Q\right)$ satisfy conditions i)-viii) from  \hyperref[Lemma11]{Lemma \ref*{Lemma11}}, and let $D\subset K_1^{str}\left(Q\right)$
be a set of distinguished 1-simplices.

Let $z=\sum_{i\in I}a_i\tau_i\in C_*^{simp,inf}\left(K\left(Q\right)\right)$ 
be a relative cycle
with
%$\partial z=w_0+w_1, w_0\in C_*^{simp,inf}\left(q\left(K\left(\partial_0Q\right)\right)\right),
%w_1
$\partial z\in C_*^{simp,inf}\left(K\left(\partial_1Q\right)\right)$.

Let a set of cancellations $\mathcal{C}$ for $z$ and a minimal
presentation of $\partial z$ be given. Let $\Upsilon,\partial\Upsilon$ be the associated simplicial sets, $\tau:\left(\mid\Upsilon\mid,\mid\partial\Upsilon\mid\right)\rightarrow 
\left(Q,\partial_1 Q\right)$ the associated continuous mapping.

Assume that we have an admissible 0-1-labeling
of $\Upsilon_1$. \\
Then there exists a relative cycle $$z^\prime=\sum_{i\in I}a_i\tau_i^\prime
\in C_*^{simp,inf}\left(K^{str}\left(Q\right),K^{str}\left(\partial_1Q\right)\right)$$ 
%with $\partial
%\sum_{i\in I}a_i\tau_i^\prime
%\in C_*^{simp,inf}\left(K^{str}\left(\partial Q\right)\right)$, 
such that:\\
i) the associated continuous mappings $$\tau,\tau^\prime:
\left(\mid\Upsilon\mid,\mid\partial\Upsilon\mid\right)\rightarrow
\left(Q,\partial_1 Q\right)$$ are homotopic by a homotopy mapping $\mid\partial\Upsilon\mid$ to $\partial Q$,\\ 
ii) if an edge of some $\tau_i$ is labeled by 1, then the corresponding edge of $\tau_i^\prime$ belongs to $D$,\\
%iii) for each $i$, all 0-simplices of $\Upsilon$ are mapped to $D_0$. If two distinct vertices of $\Upsilon$ are mapped to $\partial_0Q$, then they are mapped to distinct points in $\partial_0Q$.

\end{lem}

(Remark: The homotopy in i) does not necessarily map $\mid\partial\Upsilon\mid$ to $\partial_1Q$, but to $\partial Q$.)

\begin{proof}
First we apply the 'canonical straightening' $str_{can}$ from  \hyperref[Obs7]{Observation \ref*{Obs7}}.
The resulting chain $\sum_{i\in I} a_i str_{can}\left(\tau_i\right)$ satisfies i), but not necessarily ii). 

$\sum_{i\in I} a_i str_{can}\left(\tau_i\right)$ inherits the admissible
labeling
from $\sum_{i\in I} a_i \tau_i$. Thus we can w.l.o.g.\ restrict to the case that all $\tau_i$ belong to $K^{str}\left(Q\right)$.

%Moreover, we clearly have a {\em homotopy of pairs} after which all vertices of all $\tau_i$ belong to $K_0^{str}\left(Q\right)$.
%We assume this homotopy to be done.

Let 
$$e\in K_1^{str}\left(Q\right)$$ be a 1-labeled edge, let $x=\partial_1e\in K_0^{str}\left(Q\right), y=\partial_0e\in K_0^{str}\left(Q\right)$.
By  \hyperref[Def9]{Definition \ref*{Def9}}, the coset $\Gamma e\Gamma$ contains a unique distinguished 
1-simplex $str\left(e\right) \in D_{xy}$. (We use the notation from  \hyperref[Def9]{Definition \ref*{Def9}}, in particular $\Gamma:=\Omega\left(\partial Q\right)$.) 

$str\left(e\right)\in \Gamma e\Gamma$ means\footnotemark\footnotetext[9]{If $\partial_0e,\partial_1e\not\in\partial_1Q$, then $str\left(e\right)\in\Gamma e\Gamma$ means, of course, $str\left(e\right)=e$. Similarly, 
if only one vertex of $e$ belongs to $\partial_1Q$, then only that vertex is moved during the 
homotopy.} that there are loops $\gamma_1,\gamma_2\subset \partial Q$ 
based at $x$ resp. $y$ such that $str\left(e\right)\sim \gamma_1*e*\gamma_2$ rel.\ $\left\{0,1\right\}$.
There is an obvious homotopy between $e$ and $\gamma_1*e*\gamma_2$, which moves $\partial_1e$ along $\overline{\gamma}_1$ and $\partial_0e$ along $\gamma_2$. 
(Of course, we change 
the homotopy class relative boundary, so we can not keep
the endpoints fixed during the homotopy.) If $e$ and/or $\partial_0e$ and/or $\partial_1e$ have image in $\partial_1Q$, 
then their images remain in $\partial Q$ (and end up in $\partial_1Q$) during the homotopy.

Using  \hyperref[Obs1]{Observation \ref*{Obs1}}, the so-constructed homotopy between $e$ and $str\left(e\right)$ can be extended to a homotopy from $$\tau:\left(\mid\Upsilon\mid,\mid\partial\Upsilon\mid\right)
\rightarrow \left(Q,\partial_1 Q\right)$$
to some $$\hat{\tau}:\left(\mid\Upsilon\mid,\mid\partial\Upsilon\mid\right)\rightarrow \left(Q,\partial_1 Q\right),$$
such that $\hat{\tau}$ is
a simplicial map from $\Upsilon$ to $S_*\left(Q\right)$. (If a 0-labeled edge has one or both vertices adjacent to 1-labeled edges, then the 0-labeled edge just follows the 
homotopy of the vertices. 0-labeled edges that are not adjacent to 1-labeled edges can remain fixed during the homotopy.) The homotopy maps
$\mid\partial\Upsilon\mid$ to $\partial Q$. 

Next we apply homotopies rel.\ boundary to the (already homotoped images of)
all 0-labeled edges $f\in K_1^{str}\left(Q\right)$, to homotope them to edges 
in $K_1^{str}\left(Q\right)$. 
If $f$ and/or $\partial_0f$ and/or $\partial_1f$ 
have image in $\partial_1Q$, then their images remain in $\partial Q$ (and end up in $\partial_1Q$) during the homotopy.

%Because this will be of importance in the proof of Corollary 3, we explicitly give the homotopy for (already canonically
%straightened) 0-labeled edges. We say that $f$
%is adjacent to a 1-labeled edge if $\partial_1e$ is a vertex of $f$. (We do not consider $\partial_0e$.) \\
%Let $f\in K_1^{str}\left(Q\right)$ be a 0-labeled edge which is at both vertices adjacent
%to 1-labeled edges $e_1,e_2$, i.e.\ $$\partial_0f=\partial_1e_1=x, \partial_1f=\partial_1e_2=v$$ and $$\partial_0e_1=y,\partial_0e_2=w.$$
%Then we have seen that the homotopy from $e_1$ to $e_{xy}$ moves $\partial_1e_1$ along $e_1*\overline{e}_{xy}$, and 
%the homotopy from $e_2$ to $e_{vw}$ moves $\partial_1e_2$ along $e_2*\overline{e}_{vw}$. Hence $$str\left(f\right)=
%\left(e_{vw}*\overline{e}_2*f*e_1*\overline{e}_{xy}\right)_K.$$ Here, $\left(e_{vw}*\overline{e}_2*f*e_1*
%\overline{e}_{xy}\right)_K$ means the unique 1-simplex in $K_1\left(Q\right)$ which is homotopic rel.\ boundary to $e_{vw}*\overline{e}_2*f*e_1*
%\overline{e}_{xy}$.\\
%Similarly if only one vertex of a straight 0-labeled edge $f$ is adjacent to a 1-labeled edge $e$, then we get either
%$str\left(f\right)=
%\left(e_{vw}*\overline{e}_2*f\right)_K$ or
%$str\left(f\right)=
%\left(f*e_1*\overline{e}_{xy}\right)_K$. If $f$ is not adjacent to a 1-labeled edge, then we get $str\left(f\right)=f$.

Now we have a simplicial map $\hat{\tau}:\Upsilon\rightarrow S_*\left(Q\right)$,
such that all 1-simplices are mapped to $K_1^{str}\left(Q\right)$,
and such that $$\hat{\tau}\left(e\right)\in D\subset K_1^{str}\left(Q\right)$$ holds for all 1-labeled edges $e$. Then we can,
as in the proof of Observation 7, by induction on $n$, apply homotopies rel.\ boundary to all  
n-simplices to homotope them into $K_n^{str}\left(Q\right)$. Simplices in $\partial_1Q$ remain in $\partial Q$ (and end up in $\partial_1Q$) during the homotopy.

We obtain a homotopy (of pairs),
which keeps the 1-skeleton fixed, to 
a simplicial map $$\tau^\prime:\Upsilon\rightarrow K^{str}\left(Q\right),$$
mapping $\partial\Upsilon$ to $K^{str}\left(\partial_1 Q\right)$ and satisfying i),ii).
\end{proof}

A somewhat artificial 
formulation of the conclusion of  \hyperref[Lemma12]{Lemma \ref*{Lemma12}} is that we have 
constructed a chain map $$str:C_*^{simp,inf}\left(\Upsilon,
\partial\Upsilon\right)\rightarrow 
C_*^{simp,inf}\left(K^{str}\left(Q\right),K^{str}\left(\partial_1 Q\right)\right).$$
Unfortunately, this somewhat artificial formulation can 
not be simplified because $str$ depends on the chain $\sum_{i\in I}a_i\tau_i$. 
That is, we do not get a chain map $str:C_*^{simp,inf}\left(K\left(Q\right),K\left(\partial_1 Q\right)\right)
\rightarrow C_*^{simp,inf}\left(K^{str}\left(Q\right),K^{str}\left(\partial_1 Q\right)\right)$.

\subsection{Straightening of crushed cycles}\label{sec:crush}

Recall from Section \ref{sec:amenable} that $.\otimes_{{\bf Z}G}{\bf Z}$ means the tensor product with the trivial
${\bf Z}G$-module ${\bf Z}$, that is, the quotient under the $G$-action.
We first state obvious generalizations of Observation 6 to the case of tensor products with a factor with trivial $G$-action.

\begin{obs}\label{Obs9} Let $\left(Q,\partial_1 Q\right)$ be a 
pair of topological spaces. Let $G$ be
a group acting
on a pair $\left(K,\partial K\right)$ with $K\subset S_*\left(Q\right)$ and $\partial K\subset S_*\left(\partial_1 Q\right)$ both closed under face maps.\\
i) If $$z=\sum_{i\in I}a_i\tau_i\otimes 1\in C_*^{simp,inf}
\left(K,\partial K\right)\otimes_{{\bf Z} G}{\bf Z}$$
is a relative cycle, then $$\hat{z}=\sum_{i\in I}\sum_{g\in G} a_i
\left(g\tau_i\right)\in C_*^{simp,inf}\left(K,\partial K\right)$$ is a relative cycle. \\
If $\mathcal{C}$ is a sufficient set of cancellations for $z$, then 
there exists a set of cancellations $\widehat{\mathcal{C}}$ for $\hat{z}$ such
that $\left(\eta_1,\eta_2\right)\in
\widehat{\mathcal{C}}$ implies $\left(\eta_1\otimes 1,\eta_2\otimes 1\right)\in{\mathcal{C}}$.\\
If $\partial z=\sum_{a,i}c_{ai}\partial_a\tau_i\otimes 1$ is a minimal presentation for $\partial z$, then 
$\partial\hat{z}=\sum_{g\in G}\sum_{a,i}c_{ai}\partial_a\left(g\tau_i\right)$ is a 
minimal presentation for $\hat{z}$.\\
ii) Let $\widehat{\Upsilon},\partial\widehat{\Upsilon}$ be the simplicial sets associated to $\hat{z}$,
the sufficient set of cancellations 
$\widehat{\mathcal{C}}$ and the minimal presentation of 
$\partial\hat{z}$. They come with an obvious $G$-action.
Then we have an associated continuous mapping $\hat{\tau}:
\left(\mid \widehat{\Upsilon}\mid,\mid\partial\Upsilon\mid\right)
\rightarrow\left(Q,\partial_1 Q\right)$. \end{obs}

\begin{cor}\label{Cor3}
Let $Q,\partial Q,\partial_1 Q,\partial_0Q$ satisfy Assumption I.  Let $q:Q\rightarrow Q$ be given by  \hyperref[Cons1]{Construction \ref*{Cons1}}.
% and let $A:=q\left(\partial_0 Q\right)$.

Let $K_*^{str}\left(Q\right)\subset
S_*\left(Q\right)$ satisfy conditions i)-viii) from  \hyperref[Lemma11]{Lemma \ref*{Lemma11}}, 
and let $D\subset K_1^{str}\left(Q\right)$ be a set
of distinguished 1-simplices.

Let $G:=\Pi\left(K\left(\partial_0Q\right)
\right)$
%\subset \Pi\left(K\left(A\right)\right)$ 
with its action on $K^{str}\left(Q\right)$ defined in 
 \hyperref[Obs5]{Observation \ref*{Obs5}},
and let $H:=q_*\left(G\right)$ as defined in Section \ref{sec:connected}.
Let $$\sum_{i\in I} a_i\tau_i\otimes 1\in C_n^{simp,inf}\left(K\left(Q\right),GK\left(\partial_1 Q\right)\right)
\otimes_{{\bf Z}G}{\bf Z}$$ 
be a relative cycle. Fix a sufficient set of cancellations
$\mathcal{C}$ and a minimal presentation for $\partial z$. Let $\widehat{\Upsilon},\partial\widehat{\Upsilon}$ be defined by Observation 9.
Assume that we have a $G$-invariant admissible 0-1-labeling of the edges of $\widehat{\Upsilon}$.

Then there is a well-defined chain map $$q\circ str:C_*^{simp,inf}\left(\widehat{\Upsilon}
\right)
\otimes_{{\bf Z}G}{\bf Z}\rightarrow 
C_*^{simp,inf}\left(HK^{str}\left(Q\right)\right)\otimes_{{\bf Z}H}{\bf Z},$$
mapping $C_*^{simp,inf}\left(\partial\widehat{
\Upsilon}\right)
\otimes_{{\bf Z}G}{\bf Z}$ to 
$C_*^{simp,inf}\left(GK^{str}\left(\partial_1 Q\right)\right)\otimes_{{\bf Z}H}{\bf Z}$,
such that:\\
i) if $e\in\widehat{\Upsilon}_1$ is a 1-labeled edge, $str\left(e\otimes 1\right)=f\otimes 1$, then $f\in D$.\\
ii) if $Q$ is an orientable manifold 
with boundary $\partial Q$, 
and
if $$ \sum_{i\in I} a_i\tau_i\otimes 1\in C_*^{simp,inf}\left(
%K\left(Q\right),
K\left(Q\right),GK\left(\partial_1 Q\right)\right)\otimes_{{\bf Z}G}{\bf Z}$$ represents\footnotemark\footnotetext[10]{Cf.\ the footnotes in Section \ref{sec:central}}
the image of $\left[Q,\partial Q\right]\otimes 1$,
% under the canonical homomorphism, 
then 
$$\sum_{i\in I} a_i q\circ str\left(\tau_i\otimes 1\right)\in C_*^{simp,inf}\left(
HK^{str}\left(Q\right),
HK^{str}\left(\partial_1 Q\right)
\right)\otimes_{{\bf Z}H}{\bf Z}$$
represents
the image of $\left[Q,\partial Q\right]\otimes 1$ and 
$$\partial \sum_{i\in I} a_i q\circ str\left(\tau_i\otimes 1\right)\in C_*^{simp,inf}\left( HK^{str}\left(\partial_1 Q\right) \right)\otimes_{{\bf Z}H}{\bf Z}$$ represents the image of $\left[\partial Q\right]\otimes 1$.  % under the canonical homomorphism.  \end{cor} 
\end{cor}

\begin{proof} 

We can apply \hyperref[Lemma12]{Lemma \ref*{Lemma12}} to the infinite chain $\sum_{i\in I,g\in H}a_i \left(g\tau_i\right)$.
Thus \hyperref[Lemma12]{Lemma \ref*{Lemma12}} provides us with a chain map 
$str:C_*^{simp,inf}\left(\widehat{\Upsilon}\right)\rightarrow C_*^{simp,inf}\left(K^{str}\left(Q\right)
\right)$, given by $$str\left(g\tau_i\right):=\left(g\tau_i\right)^\prime.$$
$q:\left(K^{str}\left(Q\right),K^{str}\left(\partial_1Q\right)\right)
\rightarrow \left(K^{str}\left(Q\right),K^{str}\left(\partial_1Q\right)\right)$ is defined by Observation 5.
(Remark: we actually have $q\circ str\left(g\tau_i\right)\in
K^{str}\left(Q\right)$. We need $HK^{str}\left(Q\right)$ in the statement of  \hyperref[Cor3]{Corollary \ref*{Cor3}} just to have the tensor product well-defined.)\\

We are going to 
define $q\circ str\left(\sigma\otimes z\right):=q\left(str\left(\sigma\right)\right)\otimes z$ for each $\sigma\in \widehat{\Upsilon}, z\in{\bf Z}$. For this to be well-defined, we have to check the following claim: \\

{\em for each $\sigma\in
K,g\in G$, there exists $h\in H$ with $q\left(str\left(g\sigma\right)\right)=h q\left(str\left(\sigma\right)\right)$.}\\
\\
By condition viii) from  \hyperref[Lemma11]{Lemma \ref*{Lemma11}} (asphericity of $K^{str}\left(Q\right)$), it suffices to check this for the 1-skeleton. 

It is straightforward to check the claim for the 0-skeleton. 

If $\sigma=v\in S_0\left(\partial_0Q\right)$ then $v$ and $gv$ belong to the same path-component $C$ of $\partial_0Q$, hence $str\left(v\right)$ and $str\left(gv\right)$ belong to the same path-component $C$. 
%Thus $str\left(gv\right)=g^\prime str\left(v\right)$ for some $g^\pr\in S_0\left(\partial_0Q\right)$,
%and $v$ and $gv$
%belong to the same path-component $C$ of $\partial_0Q$. By the proof of Lemma 12, $str\left(gv\right)$ and $str\left(v\right)$ both belong to the path-component
%$C$.
Let $\gamma:\left[0,1\right]\rightarrow\partial_0Q$ be a path with $\gamma\left(0\right)=str\left(v\right), \gamma
\left(1\right)=str\left(gv\right)$. Let $\gamma^\prime$ be the unique 1-simplex in $K\left(\partial_0Q\right)$ which is 
homotopic rel.\ boundary to $\gamma$.
Let $g^\prime:=\left\{\gamma^\prime,\overline{\gamma^\prime}\right\}\in G=\Pi\left(K\left(\partial_0Q\right)\right)$.
Then $g^\prime str\left(v\right)=str\left(gv\right)$, which implies
$q\left( str\left(gv\right)\right)=h q\left( str\left(v\right)\right)$ with $h=q_*\left(g^\prime\right)\in H$.
%,str\left(gv\right)=gv$, hence
%$str\left(gv\right)=g str\left(v\right)$. This implies $q\left(str\left(gv\right)\right)=q_*\left(g\right)q\left(str\left(v\right)\right)$ by the Remark in Section 5.1. 

If $\sigma=v\not\in \partial_0Q$, then $gv=v$, hence $q\left(str\left(gv\right)\right)=q\left(str\left(v\right)\right)$.\\

The proof for 1-simplices consists of two steps. 
In the first step we prove that for $e\in K_1\left(Q\right),g\in G$ we have 
$str_{can}\left(ge\right)=g^\prime str_{can}\left(e\right)$ with $g^\prime\in G$. 
In the second step we 
show that,
if
$e\in K_1^{str}\left(Q\right)$ and $g\in G$, then there exists $h\in H$ with
$q\left(str\left(ge\right)\right)=h q\left(str\left(e\right)\right)$.
Hence altogether we will get $q\left(str\left(ge\right)\right)=q\left(str\left(str_{can}\left(ge\right)\right)\right)
= q\left(str\left(g^\prime\ str_{can}\left(e\right)\right)\right)=h\ q\left(str\left(str_{can}\left(e\right)\right)\right)=h\ q\left(str\left(e\right)\right)$.\\
\\
First step: This is fairly obvious. \\
First case: If both vertices of $e$
do not belong to $\partial_0Q$,
then also both vertices of $str_{can}\left(e\right)$ do not belong to $\partial_0Q$, and we have $ge=e, g str_{can}
\left(e\right)=str\left(e\right)$, which implies the claim.\\
Second Case: If both vertices of $e$
belong to $\partial_0Q$, then $str_{can}\left(e\right)\sim\alpha_1*e*\alpha_2, str_{can}\left(ge\right)\sim\beta_1*ge*\beta_2$ for some paths $\alpha_1,\alpha_2,\beta_1,\beta_2$ in $\partial_0Q$. Moreover, by the definition of the action (Section 3.3) we have $ge\sim\gamma_2*e*\gamma_1$ for some 
$\gamma_1,\gamma_2\in K_1\left(\partial_0Q\right)$. Thus $str_{can}\left(ge\right)\sim
\beta_1*\gamma_1*\alpha_1^{-1}*str_{can}\left(e\right)*\alpha_2^{-1}*\gamma_2*\beta_2$, in particular $str_{can}\left(ge\right)=g^\prime str_{can}\left(e\right)$ for some $g^\prime \in G$. \\
Third case: Finally we consider the case that one vertex, say $\partial_0e$ belongs to 
$\partial_0Q$, but $\partial_1 e$ does not belong. Then we are in the situation of the second case with $\gamma_2=1, \alpha_2=\beta_2$ (except that $\alpha_2$ is not contained in $\partial_0Q$). We get $str_{can}\left(ge\right)\sim \beta_1*\gamma_1*\alpha_1^{-1}*str_{can}\left(e\right)$. Since $ \beta_1*
\gamma_1*\alpha_1^{-1}$ is contained in $\partial_0Q$, this implies that 
$str_{can}\left(ge\right)=g^\prime str_{can}\left(e\right)$ for some $g^\prime \in G$.\\
%If $
%g\partial_0e=\partial_0e$, we can apply the same argument as in the 
%case $\partial_0e,\partial_1e\not\in\partial_0Q$. If $g\partial_0e
%\not=\partial_0e$, then we have (possibly after renumbering) $g=\left\{
%\gamma_1,\ldots\right\}$ with $\partial_1\gamma_1=\partial_0e$. By definition we have that $ge$ is homotopic rel.\ boundary to $e*\gamma_1$, 
%$str\left(ge\right)$ is homotopic rel.\ boundary to $\overline{H\left(\partial_1e\right)}*e*\gamma_1$, $str\left(e\right)$ is homotopic rel.\ boundary to
%$\overline{H\left(\partial_1e\right)}*e$, $g str\left(e\right)$ 
%is homotopic rel.\ boundary to $\overline{H\left(
%\partial_1e\right)}*e*\gamma_1$. Thus $str\left(ge\right)=g str\left(e\right)$.\\
\\
Second step: 
Let $e\in K_1^{str}\left(Q\right)$.

%By the construction of the action in Section 3.3 we have $ge\sim \gamma_1*e*\gamma_2$ where $\gamma_1,\gamma_2$ are either 1-simplices in $\partial_0Q$ or constant. By the proof of Observation 7 
%we have $str_{can}\left(ge\right)\sim \overline{H\left(\partial_1\gamma_1\right)}*\gamma_1*e*\gamma_2*H\left(\partial_0\gamma_2\right)$, where $H\left(\partial_1\gamma_1\right)$ resp.\ $H\left(\partial_0\gamma_2\right)$ are either constant (if $\partial_1e\not\in\partial_0Q$ resp.\ 
%$\partial_0e\not\in\partial_0Q$) or are in $\partial_0Q$ (if $\partial_1e\in\partial_0Q$ resp.\
%$\partial_0e\in\partial_0Q$). 

If $e$ is a 1-labeled edge, with $x=\partial_1e,y=\partial_0e\in K_0^{str}\left(Q\right)$, then we have by condition xiv) from Definition 9 that 
$$q\left(str\left(ge\right)\right)=hq\left(e_2\right)$$
for some $e_2\in D_{xy}$ and some $h\in H$. But $e_2$ belongs to the same coset
in $\Gamma K_1^{str}\left(Q\right)\Gamma$ as $e$, thus $e_2=str\left(e\right)$ which proves the claim for $e$.
%ih then $str\left(e\right)\in D, str\left(ge\right)\in D$. We have that $\partial_0 str\left(e\right)$ and $\partial_0 str\left(ge\right)$ belong to the same path-component of $\partial Q$, and $\partial_1 str\left(e\right)$ and $\partial_1 str\left(ge\right)$ belong %to the same path-component of $\partial Q$. By condition xiii) from Definition 9 this implies $q\left(str\left(ge\right)\right)=h q\left(str\left(e\right)\right)$ for some $h\in H$. Second, let $f$ be a 0-labeled edge with one or two vertices 
%adjacent to 1-labeled edges. 

If $f$ is adjacent to one 1-labeled edge $e$ and $q\left(str\left(ge\right)\right)=h q\left(str\left(e\right)\right)$, then $q\left(str\left(gf\right)\right)=h q\left(str\left(f\right)\right)$ because the homotopy of $f$ resp.\ $gf$ 
just followed the homotopy of $e$ resp.\ $ge$: e.g.\ if 
$\partial_1f=\partial_1e$ and $q\left(str\left(ge\right)\right)
\sim q_*\left(\alpha\right)*q\left(str\left(e\right)\right)*q_*\left(\beta\right)$ with $\alpha,\beta\in K_1\left(\partial_0Q\right)$, then $q\left(str\left(gf\right)\right)\sim q_*\left(\alpha\right)*q\left(str\left(f\right)\right)$. 
%As explained above, 
%we have 
%$str_{can}\left(ge\right)\sim \overline{H\left(\partial_1\gamma_1\right)}*\gamma_1*e*\gamma_2*H\left(\partial_0\gamma_2\right)$. By the proof of Lemma 12 we have
%$str\left(ge\right)\sim \alpha_1*str_{can}\left(ge\right)*\alpha_2$
%for some $\alpha_1,\alpha_2\in\Omega\left(\partial Q\right)$. 
%Then $str\left(gf\right)\sim\alpha_1*\overline{H\left(\partial_1\gamma_1\right)}*
%\gamma_1*f$ On the other hand $str\left(e\right)\sim\beta_1*e*\beta_2$ 
%and $str\left(f\right)\sim\beta_1*f$ for some $\beta_1,\beta_2\in\Omega
%\left(\partial Q\right)$. But $str\left(ge\right)=h str\left(e\right)$ with $h\in H$ means
%$str\left(ge\right)\sim q_*\left(\delta_1\right)*\ldots 
%*q_*\left(\delta_k\right)*str\left(e\right)$ 
%for some 1-simplices $\delta_1,\ldots,\delta_k$ in $\partial_0Q$. 
Similarly if $f$ is adjacent to two 1-labeled edges.
% $e_1$ and $e_2$, and if
%$q\left(str\left(ge_1\right)\right)=h_1 q\left(str\left(e\right)\right), $q\left(str\left(ge\right)\right)=h q\left(str\left(e\right)\right)$$.

Finally, if a 0-labeled straight 1-simplex $f$ is not adjacent to a 1-labeled edge, then $str\left(f\right)=f$ and $str\left(gf\right)=gf$, which implies $str\left(gf\right)=g str\left(f\right)$ and $q\left(str\left(gf\right)\right)=q_*\left(g\right) str\left(f\right)$.

Thus we have proved $q\left(str\left(gf\right)\right)=h q\left(str\left(f\right)\right)$ with some $h\in H$ for any 0-labeled edges $f$.\\
\\
Thus $q\circ str$ is well-defined and satisfies i) by  \hyperref[Lemma12]{Lemma \ref*{Lemma12}}. To prove ii), we first observe that, if $\sum_{i\in I} a_i\tau_i$ 
represents $\left[Q,\partial Q\right]$, then, by  \hyperref[Obs6]{Observation \ref*{Obs6}}c) and
condition i) from  \hyperref[Lemma12]{Lemma \ref*{Lemma12}} (together with $q\sim id$), we have that
$$\sum_{i\in I} a_i q\circ str\left(\tau_i\right)=\sum_{i=1}^ra_iq\left(\tau_i^\prime\right)$$ 
represents $\left[Q,\partial
Q\right]$ and the claim follows. Thus 
it suffices to check: if $\sum_{i\in I} a_i\tau_i\otimes 1$ is (relatively)
homologous to $\sum_{j\in J} b_j\kappa_j\otimes 1$, then
$q\circ str\left(\sum_{i\in I} a_i\tau_i\otimes 1\right)$ is (relatively)
homologous to $q\circ str\left(\sum_{j\in J} b_j\kappa_j\otimes 1\right)$. 

So let $$\sum_{i\in I} a_i\tau_i\otimes 1-
\sum_{j\in J} b_j\kappa_j\otimes 1=\partial \sum_{k\in K} c_k\eta_k\otimes 1\ mod\ C_*^{simp,inf}\left(GK\left(\partial_1 Q\right)\right)\otimes_{{\bf Z}G}{\bf Z}$$
for some chain $ \sum_{k\in K} c_k\eta_k\otimes 1\in
C_*^{simp,inf}\left(K\left(Q\right)\right)\otimes_{{\bf Z}G}{\bf Z}$. 
In complete analogy with  \hyperref[Lemma12]{Lemma \ref*{Lemma12}}, we may extend $ str$ to the simplicial set built by the $g\eta_k$'s, their faces and degenerations, and 
obtain a singular 
chain $q\left(str\left(\sum_{k\in K} c_k\eta_k\right)\right)$
whose boundary is 
$$\partial q\circ str\left(\sum_{k\in K} c_k\eta_k\right)=
q\circ str\left(\sum_{i\in I} a_i\tau_i\otimes 1\right)-
q\circ str\left(\sum_{j\in J} b_j\kappa_j\otimes 1
\right)\ mod\ C_*^{simp,inf}\left(HK^{str}\left(\partial_1 Q\right)\right)\otimes_{{\bf Z}H}{\bf Z}.$$
This gives the first claim of ii). The second claim of ii) follows because $\partial$ 
maps $\left[Q,\partial Q\right]$ to $\left[\partial Q\right]$.

%Since vertices contained in $\partial_0Q$ are fixed during the homotopy, we have $g str_{can}\left(e\right)=str_{can}\left(ge\right)$ for each 1-simplex $e$. 
%If $e$ is a 1-labelled edge, then its vertices belong to 
%$\partial_1Q$. If they belong to a closed component of $\partial_1Q$, then $ge=e$ and $g str\left(e\right)=str
%\left(e\right)$, from which the claim follows.
%If $e$ is adjacent to a 1-labelled edge, then $ge $ is adjacent to the same 1-labelled edge.

%The last claim holds because homotopic relative cycles (for a homotopy preserving $C_*\left(\partial Q\right)$) are relatively homologous.
\end{proof}

\subsection{Removal of 0-homologous chains}\label{sec:remove}

\begin{df}\label{Def11} Let $Q$ be an n-dimensional compact manifold with boundary $\partial Q$.
We define 
$
rmv:
S_*\left(Q\right)\rightarrow S_*\left(Q
\right)$ by $$rmv\left(\sigma\right)=0$$ if $\sigma$ is weakly degenerate ( \hyperref[Def7]{Definition \ref*{Def7}})
and $$rmv\left(\sigma\right)=\sigma$$ else.\end{df}

%It is easily checked (cf.\ the proof of Lemma 2.5.\ in \cite{k1}) that $rmv$ maps cycles to cycles and boundaries to boundaries.
\begin{lem}\label{Lemma13} Assume that $Q$ is a $n$-dimensional compact manifold with boundary $\partial Q$. Let $K_*^{str}\left(Q\right)\subset 
S_*\left(Q\right)$ satisfy the conditions i)-viii) from \hyperref[Lemma11]{Lemma \ref*{Lemma11}}.
Then $$rmv:C_*^{simp}\left(K^{str}\left(Q\right),K^{str}\left(\partial_0 Q\right)\cup K^{str}\left(\partial_1Q\right)
\right)\rightarrow C_*^{simp}\left(
K^{str}\left(Q\right),K^{str}\left(\partial_0 Q\right)\cup K^{str}\left(\partial_1Q\right)\right),$$ 
defined by $$rmv\left(\left[\sigma\right]\right):=\left[rmv\left(\sigma\right)\right],$$
is a well-defined chain
map. Moreover, if $$\sum_{j=1}^r a_j\tau_j\in C_*^{simp}\left(K^{str}\left(Q\right), K^{str}
\left(\partial_0 Q\right)\cup K^{str}\left(\partial_1Q\right)
\right)\subset C_*^{sing}\left(Q,\partial Q\right)$$
represents 
$\left[Q,\partial Q\right]$, then $\sum_{j=1}^r a_j rmv\left(\tau_j\right)$ represents
$\left[Q,\partial Q\right]$.
\end{lem}
\begin{proof}
If $\sigma\in K^{str}\left(\partial_0 Q\right)\cup K^{str}\left(\partial_1 Q\right)$, then $rmv\left(\sigma
\right)\in 
K^{str}\left(\partial_0 Q\right)\cup K^{str}\left(\partial_1 Q\right)$, thus $rmv$ is well-defined.\\
In a first step, we prove that $rmv$ is a chain map.

Assume that $rmv\left(\sigma\right)=0$.\\
If $\sigma$ has image in $\partial Q$, then $rmv\left(\sigma\right)=0$ and $rmv\left(\partial\sigma\right)=0$, thus
$\partial rmv\left(\sigma\right)=rmv\left(\partial\sigma\right)$.\\
If some edge $e$ of $\sigma$, say connecting the $i$-th and $j$-th vertex,
is a constant loop, then all faces of $\sigma$ except possibly $\partial_i\sigma$ and $\partial_j\sigma$
have a constant edge. 
Thus $rmv\left(\partial_k\sigma\right)=0$ if $k\not\in\left\{i,j\right\}$. Moreover, since $e$ is constant, corresponding edges of $\partial_i\sigma$ and $\partial_j\sigma$
are homotopic rel.\ boundary and thus agree (possibly up to orientation) by condition v) from
Lemma 11. By induction on the dimension of subsimplices we get, again 
using condition v) from  \hyperref[Lemma11]{Lemma \ref*{Lemma11}}, that $\partial_i\sigma=\left(-1\right)^{i-j}\partial_j\sigma$. Altogether we get $rmv\left(\partial\sigma\right)=0$, thus
$\partial rmv\left(\sigma\right)=rmv\left(\partial\sigma\right)$.

Assume that $rmv\left(\sigma\right)=\sigma$. Since no edge of $\sigma$ is a constant loop, of course also no edge of
a face $\partial_i\sigma$ is a constant loop. If the image of
$\partial_i\sigma$ is not contained in $\partial Q$, this implies $rmv\left(\partial_i\sigma\right)=\partial_i\sigma=\partial_i rmv\left(\sigma\right)$. If $\partial_i\sigma$ has image in $\partial Q$, then of course $\left[\partial_i\sigma\right]=\left[0\right]=\left[\partial_i rmv\left(\sigma\right)\right]$, which implies $rmv\left(\partial_i\sigma\right)=\partial_i rmv\left(\sigma\right)$.\\
Now we prove that $rmv$ sends relative fundamental cycles to relative fundamental cycles.

Let $\sum_{j=1}^r a_j\tau_j$ be a straight relative cycle, representing the relative homology class $\left[Q,\partial Q\right]$.

We denote by $J_1\subset \left\{1,\ldots,r\right\}$ the indices of those $\tau_j$ which have a constant edge.
The sum $\sum_{j\in J_1} a_j\tau_j$ is a relatively 0-homologous relative cycle. Indeed, each face of $\partial_i\tau_k$ 
not contained in $\partial Q$ has to cancel against some face of some $\tau_l$, because
$\sum_{j=1}^r a_j\tau_j$ is a relative cycle. If $\partial_i\tau_k$
is degenerate, then necessarily $l\in J_1$. 
Moreover, if $\tau_k$ is degenerate and
$\partial_i\tau_k$ is nondegenerate, then we have proved in the first step that $\partial_i\tau_k$ cancels against some $\partial_j\tau_k$. 

Thus $\sum_{j\in J_1} a_j\tau_j$ represents some relative homology class. The isomorphism 
$H_n\left(C_*^{sing}\left(Q,\partial Q\right)
\right)\rightarrow {\bf R}$ is given by pairing with the volume
form of an arbitrary Riemannian metric. After smoothing the relative cycle, we
can apply Sard's lemma, and conclude that degenerate simplices have volume 0. Thus 
$\sum_{j\in J_1} a_j\tau_j$ is 0-homologous.

We denote by $J_2\subset \left\{1,\ldots,r\right\}$ the indices of those $\tau_j$ which are contained in
$\partial Q$. For $j\in J_2$ we have $\left[\tau_j\right]=\left[0\right]\in
C_*^{sing}\left(Q,\partial Q\right)$.
%=C_*^{str}\left(Q\right)/C_*^{str}\left(\partial Q\right)$.

Thus $\sum_{j\not\in J_1\cup J_2}a_j\tau_j$ is another representative of the homology class $\left[Q,\partial Q\right]$.
But, by Definition 11, it also represents $\left(rmv\right)_*\left(\left[Q,\partial Q\right]\right)$.
\end{proof} 

Consider a subgroup $H\subset\Pi\left(K\left(A\right)\right)$ for some $A\subset\partial Q$. (E.g.\ $A=q\left(\partial_0Q\right)$ in the setting of  \hyperref[Cons1]{Construction \ref*{Cons1}}, and $H=q_*\left(\Pi\left(K\left(\partial_0Q\right)\right)\right)\subset \Pi\left(K\left(A\right)\right)$.

A 1-simplex $e$ is a constant loop if and only if $he$ is a constant loop
for all $h\in H$. This implies that a simplex $\sigma$ is degenerate if and only if $h\sigma$ is degenerate for all $h\sigma$. Moreover, $H$ maps simplices in $\partial Q$ to
simplices in $\partial Q$. Thus $rmv\left(\sigma\right)=0$ if and only if $rmv\left(h\sigma\right)=0$ for all $h\in H$, that is, $rmv$ is well defined on $C_*^{simp,inf}
\left(HK^{str}\left(Q\right)\right)\otimes_{{\bf Z}H}{\bf Z}$ for each subgroup $H$.

\begin{lem}\label{Lemma14} Assume that $Q$ is a $n$-dimensional compact manifold with boundary $\partial Q$.
Let the assumptions of  \hyperref[Cor3]{Corollary \ref*{Cor3}} be satisfied. 
Then
%, for each subgroup $H\subset G$, 
we can extend $rmv$ to a well-defined chain map $$rmv:C_*^{simp,inf}
\left(HK^{str}\left(Q\right), HK^{str}\left(\partial_1 Q\right)\right)\otimes_{{\bf Z}H}{\bf Z}
\rightarrow C_*^{simp,inf}\left(HK^{str}\left(Q\right), HK^{str}\left(\partial_1 Q\right)
\right)\otimes_{{\bf Z}H}{\bf Z}$$ by defining $$rmv\left(\sigma
\otimes z\right)=\left\{\begin{array}{cc} 0:& rmv\left(\sigma\right)=0 \\
\sigma\otimes z:&\mbox{\ else\ }\end{array}\right\}.$$
Moreover, if $\sum_{j\in J} a_j\tau_j\otimes 1\in C_*^{simp,inf}\left(HK^{str}\left(Q\right), HK^{str}\left(\partial_1 Q\right)
\right)\otimes_{{\bf Z}H}{\bf Z}
$ represents the image 
of $\left[Q,\partial Q\right]\otimes 1$, then $\sum_{\in J} a_j rmv\left(\tau_j\otimes 1\right)$ represents the image of $\left[Q,\partial Q\right]\otimes 1$.
\end{lem}
\begin{proof} Well-definedness of $rmv$ follows from the remark before  \hyperref[Lemma14]{Lemma \ref*{Lemma14}}. The same proof as for  \hyperref[Lemma13]{Lemma \ref*{Lemma13}} shows that
$rmv$ is a chain map. 
%(The only non-obvious point is the case that 
%$rmv\left( g\sigma \right)=0$ and $g\sigma$ has a constant edge, with 
%vertices $v_i$ and $v_j$. In this case all $g\partial_k\sigma$ with $k\not\in\left\{
%i,j\right\}$ have a constant edge. By the definition of the $G$-action in Observation 9, 
%we have that the action of $G$, and thus of $H\subset G$,
%fixes constant edges. Thus all $\partial_k\sigma$ with $k\not\in\left\{
%i,j\right\}$ have a constant edge. Moreover, $g\partial_i\sigma=\left(-1\right)^{i-j}g\partial_j\sigma$, which implies
%$\partial_i\sigma=\left(-1\right)^{i-j}\partial_j\sigma$,
%thus $rmv\left(\partial\sigma\right)=0$.)

If $\sum_{j=1}^r a_j\tau_j$ represents $\left[Q,\partial Q\right]$, then the second claim follows from Lemma 13.
If $\sum_{j\in J} a_j\tau_j\otimes 1$ is homologous to $\sum_{i=1}^s b_i\kappa_i\otimes 1$ and 
$\sum_{i=1}^s b_i\kappa_i$ represents $\left[Q,\partial Q\right]$, then (because $rmv$ is a
chain map) 
$rmv\left(\sum_{j\in J} a_j\tau_j\otimes 1\right)$ is homologous to 
$rmv\left(\sum_{i=1}^s b_i\kappa_i\otimes 1\right)$, which implies the second claim.\end{proof}

The proof of  \hyperref[Thm1]{Theorem \ref*{Thm1}} will pursue the idea of straightening a given cycle such that many simplices either become 
weakly degenerate or will have an edge in $\partial_0Q$. In the first case, they will disappear after application of $rmv$. In the second case, they disappear in view of the following observation, which is a variant of an argument used in \cite{gro}.

\begin{lem}\label{Lemma15} a) Let Assumption I be satisfied for a manifold $Q$ and 
consider the action of $G=\Pi\left(K\left(\partial_0Q\right)\right)$ on $K\left(Q\right)$.
Let $\sigma\in K\left(Q\right)$ be a simplex.
%a relative cycle $z=\sum_{i\in I} a_i\tau_i\otimes 1\in C_n^{simp,inf}\left(K\left(Q
%\right),K\left(\partial Q\right)\right)\otimes_{{\bf Z}G}
%{\bf Z}$. Let $\sigma\subset\tau_i$ be a subsimplex of some $\tau_i, i\in I$.

If $str\left(\sigma\right)$ has an edge in $\partial_0Q$, then $$str\left(\sigma\otimes 1\right)=0\in
C_*^{simp,inf}\left(K\left(Q
\right)\right)\otimes_{{\bf Z}G}
{\bf Z}.$$ 
b) If $q:Q\rightarrow Q$ is given by  \hyperref[Cons1]{Construction \ref*{Cons1}}, $H=q_*\left(G\right)$,
and $\sigma\in K\left(Q\right)$ a simplex such that
$q\left(str\left(\sigma\right)\right)$
has an edge in $q\left(\partial_0Q\right)$, then 
$$q\left(str\left(\sigma\otimes 1\right)\right)=0\in
C_*^{simp,inf}\left(K\left(Q
\right)\right)\otimes_{{\bf Z}H}
{\bf Z}.$$

\end{lem}
\begin{proof} a) Let $\gamma$ be the edge of $str\left(\sigma\right)$ with image in $\partial_0Q$,
then $g=\left\{\gamma,\overline{\gamma}\right\}$ is an element of 
$G=\Pi\left(K\left(\partial_0Q\right)\right)$ and $g str\left(\sigma\right)=\overline{str\left(\sigma\right)}$.
In the simplicial chain complex 
%$C_*^{simp,inf}
%\left(K^{str}\left(Q\right),K^{str}\left(\partial Q\right)\right)\subset
$C_*^{simp,inf}\left(K\left(Q\right)\right)$, one has $\overline{str\left(\sigma\right)}=-str\left(\sigma\right)$. Thus $g str\left(\sigma\right)=-str\left(\sigma\right)$, which implies $str\left(\sigma\otimes 1\right)=str\left(\sigma\right)\otimes 1=0$. \\
b) Let $\gamma$ be the edge of $q\left(str\left(\sigma\right)\right)$ with image in $q\left(\partial_0Q\right)$. Let $\gamma^\prime$ be the corresponding edge of $str\left(\sigma\right)$. Let $g=\left\{\gamma^\prime,\overline{\gamma}^\prime\right\}\in G$ and $h=q_*\left(g\right)=\left\{\gamma,\overline{\gamma}\right\}\in H$. The same argument as in a) shows $h q\left(str\left(\sigma\right)\right)=-q\left(str\left(\sigma\right)\right)$. \end{proof}

\section{Proof of Main Theorem}

\subsection{Motivating examples}\label{sec:motivating}

{\bf Example 1}: {\em Let $M$ be a connected, orientable, hyperbolic $n$-manifold, $F$ an orientable, geodesic n-1-submanifold, $Q=\overline{M-F}$. For simplicity we assume
that $M$ and $F$ are closed, thus $Q$ is a hyperbolic manifold with geodesic boundary $\partial_1Q\not=\emptyset$,
and $\partial_0Q=\emptyset$.}

{\bf Outline of proof of $\parallel M\parallel_F^{norm}\ge\frac{1}{n+1}\parallel \partial Q\parallel$}:
Start with a fundamental cycle 
$\sum_{i=1}^r a_i\sigma_i$ of $M$, 
such that $\sigma_1,\ldots, \sigma_r$ are
normal to $F$. 
Since we want to consider laminations without isolated leaves, we replace $F$ by a trivially foliated product neighborhood $\mathcal{F}$.
We can assume after a suitable homotopy that each component of $\sigma_i^{-1}\left(\partial Q\right)$ either contains no vertex of $\Delta^n$ or 
consists of exactly one vertex, and that each vertex of $\Delta^n$ belongs to $\sigma_i^{-1}\left({\mathcal{F}}\right)$, for $i=1,\ldots,r$.

Each $\sigma_i^{-1}\left(Q\right)$
consists of polytopes, which can each be further triangulated (without introducing new vertices)
in a coherent way (i.e., such that boundary cancellations between different $\sigma_i$'s will remain)
into $\tau_{i1},\ldots,\tau_{is\left(i\right)}$. 

$\sum_{i=1}^r a_i\left(\tau_{i1}
+\ldots+\tau_{is\left(i\right)}\right)$ is a relative fundamental cycle for $Q$.

For each $\sigma_i$, preimages of
the boundary leaves of $\mathcal{F}$ cut $\Delta^n$ into regions which we colour 
with black (components of $\sigma_i^{-1}\left(\mathcal{F}\right)$) and white 
(components of $\sigma_i^{-1}\left(Q\right)$). Moreover, if $\sigma_i^{-1}\left(
\partial Q\right)$ contains vertices, these vertices are coloured black. This is a canonical colouring (\hyperref[Def4]{Definition \ref*{Def4}}). 

The edges of the simplices $\tau_{i,j}$ fall into two classes:
'old edges', i.e.\ subarcs of edges of 
$\sigma_i$, and 'new edges', which are contained in the interior of some subsimplex of $\sigma_i$ of dimension $\ge 2
$.

We define the labeling of the edges of $\tau_{ij}$ such that 'old edges' 
are labelled 1 and 'new edges' are labelled 0. This is an admissible labeling (\hyperref[Def10]{Definition \ref*{Def10}}). With this labeling, we
apply the straightening procedure\footnotemark\footnotetext[11]{Under the assumptions of Example 1, straight simplices can be chosen to be the totally geodesic simplices with vertices in $S_0^{str}\left(Q\right)$. Distinguished simplices are chosen according to \hyperref[Obs8]{Observation \ref*{Obs8}}.} from Section 5 to get a straight cycle $\sum_{i=1}^r a_i\left(str\left(\tau_{i1}\right)
+\ldots+str\left(\tau_{is\left(i\right)}\right)\right)$. 
(Thus 'old edges' are straightened to distinguished 1-simplices.)

After straightening we remove all weakly degenerate simplices (simplices contained in $\partial Q$ or having a constant edge), i.e.\ we apply 
the map $rmv$ from Section \ref{sec:crush}. By \hyperref[Lemma13]{Lemma \ref*{Lemma13}}, this does not change the homology 
class.
In particular, the boundary of the relative cycle,
$\partial \sum_{i,j}a_i rmv\left(str\left(\tau_{ij}\right)\right)$ still represents the fundamental class $\left[\partial Q\right]$ of $\partial Q$.

{\em Claim: for each $\sigma_i$, after straightening there remain at most $n+1$ faces of
nondegenerate simplices $str\left(\tau_{ij}\right)$ contributing to $\partial 
\sum_{i,j} a_i rmv\left(str\left(\tau_{ij}\right)\right)$.}

In view of  \hyperref[Lemma10]{Lemma \ref*{Lemma10}}, it suffices 
to show the following subclaim: if, for a fixed $i$, $T_1=\partial_{k_1}\tau_{ij_1},T_2=\partial_{k_2}\tau_{ij_2}$ are faces of some $\tau_{ij_1}$ resp.\ $\tau_{ij_2}$ such 
that $T_1,T_2$ have a white-parallel arc (\hyperref[Def6]{Definition \ref*{Def6}}), then $rmv\left(str\left(\tau_{ij_1}\right)\right)=0, 
rmv\left(str\left(\tau_{ij_2}\right)\right)=0$.
 and in particular the corresponding straightened faces $str\left(T_1\right),
str\left(T_2\right)$ do not occur (with nonzero coefficient)
in $\partial \sum_{i,j} rmv\left(str\left(\tau_{ij}\right)\right)$. 
(Notational remark: for a subsimplex $T$ of an affine subset $S\subset \Delta^n$ we get a singular simplex $\sigma_i\mid_T$ by restricting $\sigma_i$ to $T$. We denote by $str\left(T\right)$ the straightening of $\sigma_i\mid_T$.)

To prove the subclaim, 
%assume that $T_1,T_2$ have a white-parallel arc. 
let $W$ be the white region of $\Delta^n$
containing $T_1$ and $T_2$ in its boundary. 
By assumption of the subclaim, there is a white square bounded by two arcs $e_1\subset T_1, e_2\subset T_2$ and two arcs
$f_1,f_2$ which are subarcs of edges of $\Delta^n$. (The square is a formal sum of two triangles, 
$U_1+U_2$, which are 2-dimensional faces of some $\tau_{ij}$'s.) 

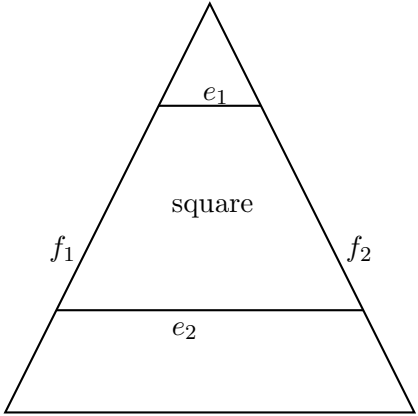
\psset{unit=0.1\hsize}
$$\pspicture(0,-1)(15,4)
\pspolygon[](0,0)(4,0)(2,4)(0,0)
\psline(0.5,1)(3.5,1)

\psline(1.5,3)(2.5,3)

\uput[0](1.5,0.8){$e_2$}
\uput[0](0.3,1.6){$f_1$}
\uput[0](1.8,3.1){$e_1$}
\uput[0](3.2,1.6){$f_2$}
\uput[0](1.5,2){square}

\endpspicture$$

We want to show that all edges of $str\left(\tau_{ij_1}\right)$ belong to $S_1^{str}\left(\partial
Q\right)$.
Note that $T_1,T_2\subset \partial W$ are mapped to $\partial Q$. Let $x_1\in S_0^{str}\left(Q\right)$
resp.\ $x_2\in S_0^{str}\left(Q\right)$ be the unique elements of
$S_0^{str}\left(Q\right)$ in the same connected component 
$C_1$ resp.\ $C_2$ of $\partial Q$ as $\sigma_i\left(
T_1\right)$ resp.\ $\sigma_i\left(T_2\right)$. 
In particular $\partial_0str\left(e_1\right)=x_1=\partial_1str\left(e_1\right)$ and $\partial_0str\left(e_2\right)=x_2=
\partial_1str\left(e_2\right)$.
Thus $e_1$ and $e_2$ are straightened to loops $str\left(e_1\right)$ resp.\ $str\left(e_2\right)$ based at $x_1$ resp.\ $x_2$. 
The straightenings of
the other two arcs, $str\left(f_1\right),str\left(f_2\right)$
connect $x_1$ to $x_2$, and they are distinguished 1-simplices because they arise as straightenings of 'old edges'.
Thus $str\left(f_1\right)=str\left(f_2\right)$, by uniqueness of distinguished 1-simplices in each coset $\Gamma K_1^{str}\left(Q\right)\Gamma$ of $\Gamma=\Omega\left(\partial Q\right)$. 
This is why we have performed the straightening construction
in Section 5 such that 
there should be only one distinguished 1-simplex, in each coset,
for any given pair of connected components.\\

This means that the square is straightened to a cylinder. 

But $\left(Q,\partial Q\right)$ is acylindrical, thus either both $str\left(e_1\right)$ and $str\left(e_2\right)$ are
constant (in which case $rmv\left(str\left(\tau_{ij_1}\right)\right)=rmv\left(str\left(\tau_{ij_2}\right)\right)=0$),
or the cylinder 
must be homotopic into $\partial Q$. In the latter case, $str\left(f_1\right)$ must (be
homotopic into and therefore) be contained in $\partial Q$. In particular, $\partial_0 str\left(f_1\right)$ and $\partial_1 str\left(f_1\right)$ belong to the same component of $\partial Q$. This implies
$\partial_0 str\left(f_1\right)=\partial_1 str\left(f_1\right)$. Since $str\left(f_1\right)$ is a distinguished 1-simplex, this implies that $str\left(f_1\right)$ is constant.

Let $P_1,P_2$ be the affine planes whose intersections with $\Delta^n$ contain $T_1$ resp.\ $T_2$.
We have now that there is an arc $f_1$ connecting 
$P_1\cap\Delta^n$ to $P_2\cap\Delta^n$ such that $str\left(f_1\right)$ is contained in $\partial Q$. This
implies that for each other arc $f$ connecting
$P_1\cap\Delta^n$ to $P_2\cap\Delta^n$ its straightening $str\left(f\right)$
must (be homotopic into and therefore) be contained in $\partial Q$.

If $P_1$ and $P_2$ are of the same type, then all edges of $str\left(\tau_{ij_1}\right)$ connect
$P_1\cap\Delta^n$ to $P_2\cap\Delta^n$, hence all
edges of $str\left(\tau_{ij_1}\right)$ belong to $S_1^{str}\left(\partial Q\right)$.
If $P_1$ and $P_2$ are not of the same type, then
existence of a parallel arc implies that
at least one of them, say $P_1$, must be of type $\left\{0a_1\ldots a_k\right\}$ with $k\not\in\left\{0,n-1\right\}$. Then, if $P_3$ 
is any other plane bounding $W$, it follows from \hyperref[Cor2]{Corollary \ref*{Cor2}} that $P_3$ has a white-parallel arc with $P_1$. 
Thus, repeating the argument in the last paragraph with $P_1$ and $P_3$ in place of $P_1$ and $P_2$, we 
obtain that for each arc
$f$ connecting
$P_1\cap\Delta^n$ to $P_3\cap\Delta^n$ its straightening $str\left(f\right)$
must (be homotopic into and therefore) be 
contained in $\partial Q$. Hence, for each $\tau_{ij_1}$ in the chosen triangulation of $W$, its 1-skeleton is straightened into $\partial Q$. 

Since straight simplices $\sigma$ (of dimension $\ge 2$) with $\partial\sigma$ 
in the geodesic boundary $\partial Q$, must be in $\partial Q$,
this implies by induction that the $k$-skeleton of $str\left(\tau_{ij_1}\right)$ is in $\partial Q$, for each $k$.
In particular, $str\left(\tau_{ij_1}\right)\in S_n^{str}\left(\partial Q\right)$. Hence $rmv\left(
str\left(\tau_{ij_1}\right)\right)=0$.
We have proved the subclaim.\\
\\
By  \hyperref[Lemma10]{Lemma \ref*{Lemma10}}, the {\em subclaim} implies the {\em claim}. Since $\sum_{i=1}^r a_i \partial \sum_j rmv\left(str\left(\tau_{ij}\right)\right)$ represents the fundamental class $\left[\partial Q\right]$,
we conclude 
$\parallel \partial Q\parallel\le \left(n+1\right)\sum_{i=1}^r\mid a_i\mid$.\\
\\
\\
The simplifications of Example 1 in comparison to the proof in Section \ref{sec:proof} are essentially all due to the fact that $\partial_0Q=\emptyset$. We remark that in Example 2, if $F$ is not geodesic, then $Q\not=N$ and thus $\partial_0Q\not=\emptyset$ (even though 
$\partial M=\emptyset, \partial F=\emptyset$). Thus the generalization to $\partial_0Q\not=\emptyset$ would be necessary even if one only wanted to consider closed manfifolds $M$ and $F$. \\
\\
{\bf Example 2}: {\em Let $M$ be a connected, closed, hyperbolic 3-manifold, $F\subset M$ a closed, incompressible surface, $N=\overline{M-F}, Q=Guts\left(N\right)$. }

{\bf Outline of proof of $\parallel M\parallel_F^{norm}\ge\frac{1}{4}\parallel \partial Q\parallel$}:
Start with a fundamental cycle
$\sum_{i=1}^r a_i\sigma_i$ of $M$,
such that $\sigma_1,\ldots, \sigma_r$ are
normal to $F$.
As in Example 1 we get a relative fundamental cycle $\sum_{i=1}^r a_i\left(\tau_{i1}                                                                             +\ldots+\tau_{is\left(i\right)}\right)$
of $N$. We can not apply the argument from Example 1 to $N$ because $N$ is not acylindrical. Therefore we 
would like to work with
a relative fundamental cycle for the acylindrical manifold $Q$.

$N$ is aspherical. Using  \hyperref[Lemma2]{Lemma \ref*{Lemma2}}, we can assume that all $\tau_{ij}$ belong to $K\left(N\right)$. Then we can apply the retraction $r$ from  \hyperref[Lemma5]{Lemma \ref*{Lemma5}}. Since $r$ is only defined after tensoring $.\otimes_{{\bf Z}G}{\bf Z}$,
we get $r\left(\tau_{ij}\otimes 1\right)=\kappa_{ij}\otimes 1$ with $\kappa_{ij}\in K\left(Q\right)$ only determined up to choosing one $\kappa_{ij}$ in its $G$-orbit. 

Since $Q$ is aspherical, we have $K\left(Q\right)=\widehat{K}\left(Q\right)$, that is, the $\kappa_{ij}$ can be considered as simplices in $Q$  and we can apply  \hyperref[Lemma6]{Lemma \ref*{Lemma6}}b) to obtain a fundamental cycle for $\partial Q$.

The rest of the proof then basically boils down to
copying the proof of Example 1 (with $\tau_{ij}$ replaced by $\kappa_{ij}$), but taking care of the ambiguity in the choice of $\kappa_{ij}$. The details can be found in the next section.

\subsection{Proof}\label{sec:proof}

We refer to the introduction for the statement of  \hyperref[Thm1]{Theorem \ref*{Thm1}} and the relevant definitions. In this section we are going to prove  \hyperref[Thm1]{Theorem \ref*{Thm1}}.\\

{\em Proof:}
%We prove the claimed inequality for the normal Gromov norm. Afterwards, at the end, we will give the necessary modifications for the proof of the stronger inequality for the (transverse) Gromov norm.\\

If $n=1$, then  \hyperref[Thm1]{Theorem \ref*{Thm1}} is trivially true. Hence we can restrict to the case $n\ge 2$.

If $\partial_1 Q$ were empty, then $\partial
Q=\partial_0Q$ and amenability of $\pi_1\partial_0Q$ would imply
$\parallel\partial Q\parallel=0$, in particular  \hyperref[Thm1]{Theorem \ref*{Thm1}} would be trivially true. Hence we can restrict to the case
$\partial_1Q\not=\emptyset$. In particular, $Q$ satisfies Assumption I from Section 5.\\

Consider a relative cycle $\sum_{i=1}^r a_i\sigma_i$, representing $\left[M,\partial M\right]$,
such that $\sigma_1,\ldots, \sigma_r$ are
normal to $\mathcal{F}$. 
%(This is in particular the case if
%$\sigma_1,\ldots, \sigma_r$ are
%transverse to $\mathcal{F}$.)
Our aim is to show: 
$\sum_{i=1}^r\mid a_i\mid\ge \frac{1}{n+1}\parallel\partial Q\parallel$.\\
% (resp.\ $\sum_{i=1}^r\mid a_i\mid\ge \frac{1}{\left[\frac{n}{2}\right]+1}\parallel\partial Q\parallel$ if $\sigma_1,\ldots,\sigma_r$ are even transverse to $\mathcal{F}$).\\

Denote $$N=\overline{M-{\mathcal{F}}}.$$ 
Since each $\sigma_i$ is normal to $\mathcal{F}$,
we have for each $i=1,\ldots,r$ that, 
after application of a simplicial homeomorphism $h_i:\Delta^n
\rightarrow\Delta^n$,
the image of
$\sigma_i^{-1}\left(N\right)$
consists of polytopes, which can each be further triangulated in a coherent way (i.e., such that boundary cancellations between different $\sigma_i$'s will remain)
into simplices $\theta_{ij},j\in \hat{J}_i$. (It is possible that
$\mid \hat{J}_i\mid=\infty$, because $N$ may be noncompact.) 
We choose these triangulations of the $\sigma_i^{-1}\left(N\right)$ 
to be minimal (\hyperref[Def6]{Definition \ref*{Def6}}), that is, we do not introduce new vertices. (Indeed, 
compatible minimal triangulations of the $\sigma_i^{-1}\left(N\right)$ do exist: one starts with common minimal triangulations of 
the common faces and extends them to minimal triangulations of each polytope.)
 
Because boundary cancellations are preserved, we have that $$\sum_{i=1}^r a_i\sum_{j\in \hat{J}_i}\theta_{ij}
$$ is a countable (possibly infinite) relative cycle 
representing the fundamental class
$\left[N,\partial N\right]$ in the sense of section 3.2.\\
We fix a sufficient set of cancellations ${\mathcal{C}}^M$ for the 
relative cycle $\sum_{i=1}^r a_i\sigma_i$, in the sense of Definition 8. 
This induces a sufficient set of cancellations ${\mathcal{C}}^N$ for the relative cycle $\sum_{i=1}^r \sum_{j\in \hat{J}_i}a_i\theta_{ij}$.\\

%Recall that we can assume that $\mathcal{F}$ has no isolated leaves, upon replacing isolated leaves with canonically foliated $I$-bundles.
If $\partial M$ is a leaf of $\mathcal{F}$, then all faces of $z$ contributing to $\partial z$ are contained in $\partial N$. We call these faces exterior faces. We can assume that, for each $i$, \\
- each component of $\sigma_i^{-1}\left(\partial N\right)$ either contains no vertex of $\Delta^n$, or
consists of exactly one vertex, or consists of an exterior face,\\
- and each vertex of $\Delta^n$ belongs to $\sigma_i^{-1}\left(\mathcal{F}\right)$.\\
Indeed, by a small homotopy of the relative 
fundamental cycle $\sum_{i=1}^r a_i\sigma_i$, preserving normality, we can obtain that no component of $\sigma_i^{-1}\left(\partial N\right)$ contains 
a vertex of $\Delta^n$, except for exterior faces. Afterwards, if some vertices of $\sum_{i=1}^r a_i\sigma_i$ do not belong to $\mathcal{F}$, 
we may homotope a small neighborhood of the vertex,
until the vertex (and no other point of the neighborhood) meets $\partial N$. This, of course, 
preserves normality to
%of $M$. The homotopy can be chosen to preserve normality of the $\sigma_i$ to 
$\mathcal{F}$. 
%Namely, if some vertices belong to $N$, then we may homotope these vertices until
%$\partial N$, without affecting the intersection of the $\sigma_i$ with $\mathcal{F}$, and afterwards we may slightly homotope the 
%vertices into the interior of the 
%laminated part, such that near to the vertices the intersections of leaves 
%$with the simplex are normal.

Since each $\sigma_i$ is normal to $\mathcal{F}$, in particular each $\sigma_i$ is normal to the union of boundary leaves 
$$\partial_1 N:=\overline{\partial N - 
\left(\partial M\cap \partial N\right)}.$$
Thus for each $\sigma_i$, after application of a simplicial
homeomorphism $h_i:\Delta^n\rightarrow\Delta^n$,
the image of $\sigma_i^{-1}\left(\partial_1 N\right)$ consists of a (possibly infinite) 
set $$Q_1,Q_2,\ldots\subset\Delta^n,$$
such that $$Q_i=P_i\cap\Delta^n$$ for some
affine hyperplanes $P_1,P_2,\ldots$. 
We define a colouring by declaring that (images under $h_i$ of)
components of $$\sigma_i^{-1}
\left(int\left(N\right)\right):=\sigma_i^{-1}\left(N-\partial_1N\right)$$ are coloured
white and (images under $h_i$ of)
components of $\sigma_i^{-1}\left({\mathcal{F}}\right)$ 
are coloured black. (In particular, all $Q_i$ are coloured black.)
Since we assume that all vertices of $\Delta^n$ belong to 
$\sigma_i^{-1}\left({\mathcal{F}}\right)$, and since each boundary leaf is adjacent to at least one component of $\sigma_i^{-1}\left(int\left(N\right)\right)$, this is a canonical colouring (\hyperref[Def4]{Definition \ref*{Def4}}).\\
\\

%All vertices of the relative cycle $$z=\sum_{i=1}^r a_i\sum_{j\in \hat{J}_i}\theta_{ij}$$ belong to $\partial N$.
%We can homotope the relative cycle such that the following condition holds for its vertices: if a vertex $v$
%of some $\theta_{ij}$ belongs to some path-component $C$ of $\partial N$ with $C\cap\partial_0Q\not=\emptyset$, then $v\in\partial_0Q$.

By \hyperref[Lemma2]{Lemma \ref*{Lemma2}}a), 
we can homotope the relative cycle $\sum_{i=1}^r\sum_{j\in\hat{J}_i}\theta_{ij}\in C_n^{inf}\left(N,\partial N\right)$
to a relative cycle $$\sum_{i=1}^r a_i\sum_{j\in \hat{J}_i}\hat{\theta}_{ij}$$
such that each $\hat{\theta}_{ij}$ is a simplex of $\widehat{K}\left(N\right)$, as defined
in Section \ref{sec:construction}, and such that the boundary $\partial \sum_{i=1}^r\sum_{j\in\hat{J}_i}\theta_{ij}$
is homotoped into $\widehat{K}\left(\partial N\right)$. 
%Then we apply Lemma 2c) and obtain that
Then consider $$\sum_{i=1}^r\sum_{j\in\hat{J}_i}a_i\tau_{ij}:=
\sum_{i=1}^r\sum_{j\in\hat{J}_i}a_ip\left(\hat{\theta}_{ij}\right)\in
C_n^{simp,inf}\left(K\left(N\right)\right),$$
%represents $\left[N,\partial N\right]$, 
where $p:\widehat{K}\left(N\right)\rightarrow K\left(N\right)$ is the projection defined at the end of Section \ref{sec:construction}, and $\tau_{ij}:=p\left(\hat{\theta}_{ij}\right)$ for all $i,j$. 
%(Remark: to apply Lemma 2c), 
%we need $\parallel N,\partial N\parallel\not=0$. If we had $\parallel N,\partial N\parallel=0$, 
%then Lemma 6 would easily imply $\parallel
%\partial Q\parallel=0$ and the conclusion of Theorem 1 would be trivially true. Hence we can restrict to the case that
%the assumption of 2c) holds true.) \\

Consider $Q\subset N$ as in the assumptions of  \hyperref[Thm1]{Theorem \ref*{Thm1}}.
We denote $$G:=\Pi\left(K\left(\partial_{0} Q\right)\right).$$

We have by assumption that $N=Q\cup R$ is an essential decomposition (as defined in the introduction), which means exactly that
the assumptions of \hyperref[Lemma5]{Lemma \ref*{Lemma5}} are satisfied.
%with $R^\prime=N^\prime\cap R, Q^\prime=N^\prime\cap Q=\partial Q, Q_0^\prime=\partial_{0}Q$. 
Thus, according to \hyperref[Lemma5]{Lemma \ref*{Lemma5}}, there exists a retraction 
$$r:C_n^{simp,inf}\left(K\left(N\right)
\right)\otimes_{{\bf Z}G}{\bf Z}
\rightarrow C_n^{simp,inf}\left(K\left(Q\right)\right)\otimes_{{\bf Z}G}{\bf Z}$$ 
for $n\ge 2$, mapping $C_n^{simp,inf}\left(GK\left(\partial N\right)
\right)\otimes_{{\bf Z}G}{\bf Z}$ to $C_n^{simp,inf}\left( GK\left(\partial_1 Q
\right)\right)\otimes_{{\bf Z}G}{\bf Z}$, 
such that, for each simplex $\tau_{ij}\in K\left(N\right)$, we either have $r\left(\tau_{ij}\otimes 1\right)=0$ or
$$r\left(\tau_{ij}\otimes 1\right)=\kappa_{ij}\otimes 1$$ for some simplex $\kappa_{ij}\in K\left(Q\right).$
(Recall that we assume from the beginning $n\ge 2$.)

%(According to Lemma 5, the image of $r$ is actually
%$C^{simp,lf}_*\left(K\left(Q\right),K\left(\partial Q\right)\right)\otimes_{{\bf Z}G}{\bf Z}$. 
%Compactness of $Q$ implies that $C^{simp,lf}_*\left(K\left(Q\right),K\left(\partial Q\right)\right)=C_*^{simp}
%\left(K\left(Q\right),K\left(\partial Q\right)\right)$.)

Thus $$r\left(\sum_{i=1}^r a_i\sum_{j\in \hat{J}_i}\tau_{ij}\otimes 1\right)=
\sum_{i=1}^r a_i\sum_{j\in J_i}\kappa_{ij}\otimes 1$$ with $J_i\subset\hat{J}_i$ for all $i$. (It may still be possible that $\mid J_i\mid=\infty$.)\\
We remark that $\kappa_{ij}$ is only determined up to choosing one $\kappa_{ij}$ in its G-orbit.

Since $r$ is a chain map, we get a sufficient set of cancellations for \\
$\sum_{i=1}^r a_i\sum_{j\in J_i}\kappa_{ij}\otimes 1$
by $${\mathcal{C}}^Q:=\left\{\left(\partial_{k}\kappa_{i_1j_1}\otimes 1, \partial_{l}\kappa_{i_2j_2}\otimes 1\right): \left(\partial_{k}\tau_{i_1j_1},
\partial_{l}\tau_{i_2j_2}\right)\in{\mathcal{C}}^N\right\}.$$

By assumption, $Q$ is aspherical. We can therefore apply 
%part b) of 
\hyperref[Lemma6]{Lemma \ref*{Lemma6}} and have that 
$$\partial\left(\sum_{i=1}^r a_i \sum_{j\in J_i} \kappa_{ij}\otimes 1\right)\in C_*^{simp,inf}
\left(GK\left(\partial_1 Q\right)\right)\otimes_{{\bf Z}G}
{\bf Z}$$
represents (the image of) $\left[\partial Q\right]\otimes 1$.

%Amenability of $\pi_1\partial_0Q$ implies that $\parallel\partial Q\parallel=\parallel\partial_1Q\parallel$, see \cite{gro} or Theorem 1 in \cite{k2}.

\hyperref[Lemma4]{Lemma \ref*{Lemma4}}a) gives that $G$ is amenable. Together with  \hyperref[Lemma7]{Lemma \ref*{Lemma7}} this implies $$\parallel\partial Q\parallel
\le \sum_{i=1}^r \mid a_i\mid
\left(n+1\right)\mid J_i\mid.$$

In the remainder of the proof, we will use  \hyperref[Lemma14]{Lemma \ref*{Lemma14}} to improve this inequality and, in 
particular, get rid of the unspecified (possibly infinite)
numbers $\mid J_i\mid$.\\

%We assumed that all vertices of $\theta_{ij}$ either belonged to $\partial_0Q$ or to some path-component of $\partial N$ that is disjoint from $\partial_0Q$. By construction of $r$, this implies that all 
%vertices of $\kappa_{ij}$ either belong to $\partial_0Q$ or belong to some path-component $C$
%of $\partial Q$ with $C\cap\partial_0Q=\emptyset$. The same will, by construction of $str$, be true for all vertices 
%of $str\left(\kappa_{ij}\right)$.\\

$Q,\partial Q, \partial_0Q, \partial_1Q$
satisfy Assumption I from Section 5. Thus there exists a simplicial
set $$K_*^{str}\left(Q\right)\subset S_*\left(Q\right)$$ 
satisfying conditions i)-viii) from  \hyperref[Lemma11]{Lemma \ref*{Lemma11}}, and a set $$D\subset K_1^{str}\left(Q\right)$$
of distinguished 1-simplices (\hyperref[Def9]{Definition \ref*{Def9}}).

Recall that, for each $i$, $$\sum_{j\in\hat{J}_i} \theta_{i,j}$$ was defined by choosing a triangulation
of $\sigma_i^{-1}\left(N\right)$.
%, where the triangulations for different $\sigma_i$'s had to be chosen 
%compatibly at common boundary faces. 
The simplices $\theta_{i,j}$ thus have 'old edges', i.e.\ subarcs of edges of $\sigma_i$, and 'new edges', whose interior
is contained in the interior of some subsimplex of $\sigma_i$ of dimension $\ge 2$. 

Associated to $z=\sum_{i=1}^r a_i\sum_{j\in\hat{J}_i} \theta_{ij}$ and ${\mathcal{C}}^N$ (and 
an arbitrary minimal presentation of $\partial z$) are, by \hyperref[Def8]{Definition \ref*{Def8}}, simplicial sets $\Upsilon^N,\partial\Upsilon^N$. 
%Let $\Upsilon^r:=\left\{\tau\in\Upsilon^N:r\left(\tau\otimes 1\right)\not=0\right\}\subset\Upsilon$. 

The only possibility that two 'old edges' have a vertex in $\Upsilon^N$ in common is that this vertex is a vertex of $\sigma_i$. 
%Hence there does not exist a cycle of 'old edges' in $\Upsilon^N$. 

So the labeling of edges of $$\sum_{i=1}^r a_i\sum_{j\in\hat{J}_i}\theta_{ij}$$ by labeling 'old edges' not containig a vertex of any $\sigma_i$ with
label 1 and all other edges with label 0 is an admissible labeling (Definition 10).\\

Associated to $$w=\sum_{i=1}^ra_i \sum_{j\in J_i} \kappa_{ij}\otimes 1$$ and ${\mathcal{C}}^Q$ 
(and
an arbitrary minimal presentation of $\partial w$) there are simplicial sets $\Upsilon,\partial\Upsilon$.
By our definition of ${\mathcal{C}}^Q$, $\Upsilon$ is isomorphic to a 
simplicial subset of $\Upsilon^N$, namely to the subset generated by the set 
$$\left\{\tau\in\Upsilon^N:r\left(\tau\otimes 1\right)\not=0\right\}$$ together with all iterated faces and degenerations.
In particular, the admissible 0-1-labeling of $\Upsilon^N$ induces an admissible 
0-1-labeling of $\Upsilon$.
% (even though possibly the $q\left(\kappa_{ij}\right)$'s may yield a 1-labeled cycle in $Q$).

%For each $\left(i,j\right)$, we give the edges of $\kappa_{i,j}$ the same labelling as the 
%corresponding edges of $\tau_{i,j}$. 
%Then $\sum_{i=1}^r a_i\sum_{j=1}^{s\left(i\right)}\kappa_{ij}\otimes 1\in 
%C_*\left(K,\partial K\right)$ is a relative cycle, with an admissible 0-1-labelling
%of its edges.

%$Q$ satisfies Assumption I. 
By  \hyperref[Cons1]{Construction \ref*{Cons1}}, there is a map of triples $q:\left(Q,\partial Q,\partial_1Q\right)\rightarrow \left(Q,\partial Q,\partial_1Q\right)$ which is
(as a map of triples) homotopic to the identity, and such that $q\left(\partial_0Q\cap C\right)$ is path-connected for each path-component $C$ of $\partial Q$.

We denote $$A:=q\left(\partial_0Q\right), H:=q_*\left(G\right)=q_*\left(\Pi\left(K\left(\partial_0Q\right)\right)\right)\subset \Pi\left(K\left(A\right)\right).$$

%The remark after Construction 1 (Section 5.2) means that $q$ induces a well-defined homomorphism
%$$q: C_*^{simp,inf}
%\left(K\left(Q,\partial Q\right)\right)\otimes_{{\bf Z}G}
%{\bf Z}\rightarrow 
%C_*^{simp,inf}
%\left(K\left(Q,\partial Q\right)\right)\otimes_{{\bf Z}H}
%{\bf Z}.$$

We observe that $H$ is a quotient of $G$, hence amenable, even though $\Pi\left(K\left(A\right)\right)$ need not be amenable.

Let $\widehat{\Upsilon},\partial\widehat{\Upsilon}$ be defined by  \hyperref[Obs9]{Observation \ref*{Obs9}}.
By \hyperref[Cor3]{Corollary \ref*{Cor3}}, there is a chain map 
$$q\circ str:C_*^{simp,inf}\left(\widehat{\Upsilon}
\right)\otimes_{{\bf Z}G}{\bf Z}
\rightarrow
C_*^{simp,inf}\left(HK^{str}\left(Q\right)
\right)\otimes_{{\bf Z}H}{\bf Z},$$
mapping $C_*^{simp,inf}\left(\partial\widehat{\Upsilon}
\right)\otimes_{{\bf Z}G}{\bf Z}$ to
$C_*^{simp,inf}\left(
HK^{str}\left(\partial_1 Q\right)
\right)\otimes_{{\bf Z}H}{\bf Z}$, such that 
$$\partial \sum_{i=1}^ra_{i} \sum_{j\in J_i}q\left(str\left(\kappa_{ij}\right)\right)\otimes 1$$
represents (the image of) $\left[\partial Q\right]\otimes 1$
and such that 1-labeled edges 
are mapped to distinguished 1-simplices. (We keep in mind that $\kappa_{ij}$ is only determined up to G-action, thus $q\left(str\left(\kappa_{ij}\right)\right)$
is determined only up to choosing one simplex in its $H$-orbit.)

We then apply \hyperref[Lemma14]{Lemma \ref*{Lemma14}} to get the cycle
%$$\sum_{i=1}^r a_i \sum_{j\in J_i}
%rmv\left(str\left(q\left(\kappa_{ij}\right)\right)\otimes 1\right)\in C_*^{simp,inf}\left(K^{str}\left(Q\right),K^{str}\left(\partial Q\right)\right)\otimes_{{\bf Z}H}{\bf Z}$$
%representing (the image of) $\left[Q,\partial Q\right]\otimes 1$. Thus
$$\partial \sum_{i=1}^r a_i \sum_{j\in J_i}
rmv\left(q\left(str\left(\kappa_{ij}\right)\right)\otimes 1\right)\in C_*^{simp,inf}\left(HK^{str}\left(\partial_1 Q\right)\right)
\otimes_{{\bf Z}H}{\bf Z}$$ representing (the image of) $\left[\partial Q\right]\otimes 1$. 
We want to show that it is actually 
a finite chain of $l^1$-norm at most $\left(n+1\right)\sum_{i=1}^r\mid a_i\mid$.\\
\\
{\bf Claim}: {\em For each $i$, $$\partial\sum_{j\in J_i}rmv\left( 
q\left(str\left(\kappa_{ij}\right)\right)\otimes 1\right)$$
is the formal sum of at most $n+1$ n-1-simplices $L\otimes 1$ with coefficient 1.} \\

%This will imply the wanted inequality
%$$\parallel\partial Q\parallel\le\sum_{i=1}^r\mid a_i\mid \left(n+1\right).$$

%To prove the claim, our aim is to show that, under the assumption of pared acylindricity of $\left(Q,\partial_1Q\right)$, for 
%each $i$, the set of n-1-simplices $L$ such that $L\otimes 1$
%occurs with coefficient 1 in 
%$\partial\sum_{j\in J_i}rmv\left(
%str\left(\kappa_{i,j}\right)\otimes 1\right)$
%can 
%not contain any n-1-simplices $L_1\otimes 1=str\left(r\left(T_1\otimes 
%1\right)\right)$ and
%$L_2\otimes 1=str\left(r\left(T_2\otimes
%1\right)\right)$
%such that $T_1$ 
%does have a 
%white-parallel arc with with some $T_2$ (Definition 5).
The claim will be a consequence of the following subclaim and Lemma 10.\\
%This will imply, with Lemma 10,
%that there can be at most $n+1$ n-1-simplices 
%occuring with coefficient 1, yielding the wanted inequality.\\
\\
{\bf Subclaim}: 
{\em Assume that for some fixed $i\in I$, for the 
chosen triangulation 
$$\sigma_i^{-1}\left(N\right)=\bigcup_{j\in \hat{J}_i}\theta_{ij},$$ 
and the associated canonical colouring, there exist 
$j_1,j_2\in \hat{J}_i, k_1,k_2\in \left\{0,\ldots,n\right\}$ such that the
faces $$T_1=\partial_{k_1}\theta_{ij_1}\in S_{n-1}\left(\partial N\right), T_2=\partial_{k_2}\theta_{ij_2}\in S_{n-1}\left(\partial N\right)$$ 
have a white-parallel arc (Definition 6). Then $$rmv\left(q\left(str\left(\kappa_{ij_1}\right)\right)\otimes 1\right)=0,
rmv\left(q\left(str\left(\kappa_{ij_2}\right)\right)\otimes 1\right)=0.$$}

{\bf We are going to prove the subclaim.}  $$\partial_{k_l}\theta_{ij_l}\in S_{n-1}\left(\partial N\right)$$
%S_{n-1}\left(N^\prime\right)$$ 
implies (by  \hyperref[Lemma5]{Lemma \ref*{Lemma5}} and  \hyperref[Cons1]{Construction \ref*{Cons1}})
$$\partial_{k_l} q\left(str\left(\kappa_{ij_l}\right)\right)\in HK_*^{str}\left(\partial_1 Q\right)$$ for $l=1,2$.
Argueing by contradiction, we assume that $$rmv\left(q\left(str\left(\kappa_{ij_1}\right)\right)\otimes 1\right)\not=0.$$

By assumption of the subclaim, there are white-parallel arcs $e_1,e_2$ of $T_1$ resp.\ $T_2$.
This means that there are arcs
$e_1,e_2$ in a 2-dimensional subsimplex $\tau^2  \subset\Delta^n$ of the standard simplex, and that there are arcs $f_1,f_2$, which are subarcs of some edge of $\tau^2$, such that
$$\partial_0e_1=\partial_1f_2,\partial_0f_2=\partial_0e_2,\partial_1e_2=\partial_0f_1,\partial_1f_1=\partial_1e_1$$ 
and such that $$e_1,f_2,e_2,f_1$$ bound a square in the boundary of 
a white component. (Cf.\ the picture in section 6.1. We will use the same letter for an affine subset of $\Delta^n$ and 
for the singular simplex obtained by restricting $\sigma_i$ to this subset.) The square is of the form $U_1+U_2$, where $U_1,U_2$ are n-2-fold iterated faces of some $\theta_{ij}$'s.
Hence $$\partial U_1=e_1+f_2+\partial_2U_1$$ and $$\partial U_2=-e_2-f_1-\partial_2U_1,$$ i.e. 
$$\partial \left(U_1+U_2\right)=e_1+f_2-e_2-f_1$$ and
$$\partial_2U_1=-\partial_2U_2.$$ 
We emphasize that we assume $e_1$ resp.\ $e_2$ to be edges of
$\theta_{ij_1}$ resp.\ $\theta_{ij_2}$ but that $f_1,f_2$ need not be edges of $\theta_{ij_1}$ or $\theta_{ij_2}$.\\
\\
Notational convention: {\em for each iterated face $f=\partial_{k_1}\ldots\partial_{k_l}\theta_{ij}$ with $i\in I,j\in J_i$, we
will denote $f^\prime$ the $n-l$-simplex with
$$f^\prime\otimes 1=\partial_{k_1}\ldots\partial_{k_l}\kappa_{ij}\otimes 1=
r\left(\partial_{k_1}\ldots\partial_{k_l}\tau_{ij}\otimes 1\right)=r\left(
\partial_{k_1}\ldots\partial_{k_l} p\left(\hat{\theta}_{ij}\right)\otimes 1\right).$$}
(The last two equations are 
true because $r,p$ and the homotopy from $\sum_{i,j}a_i\theta_{ij}$ to $\sum_{i,j}\hat{\theta}_{ij}$ are chain maps.) In other words, if $f$ is an iterated face of some $\tau_{ij}$, then $f^\prime$ is, up to the ambiguity by the $H$-action, the corresponding iterated face of $\kappa_{ij}$.\\

By \hyperref[Lemma5]{Lemma \ref*{Lemma5}} we have $e_1^\prime,e_2^\prime\in GK\left(\partial_1Q\right)$.
Thus we can (and will) choose $\kappa_{i{j_1}},\kappa_{ij_2}$ in their $G$-orbits such that we have $e_1^\prime,e_2^\prime\in K\left(\partial_1Q\right)$, hence $str\left(e_1^\prime\right), str\left(e_2^\prime\right)\in K^{str}\left(\partial_1Q\right)$.

Since $r,p$ and the homotopy are chain maps, we have 
%$$\partial_0e_1^\prime\otimes 1=\partial_1f_2^\prime\otimes 1,\partial_0f_2^\prime\otimes 1=\partial_0e_2^\prime\otimes 1,
%\partial_1e_2^\prime\otimes 1=\partial_0f_1^\prime\otimes 1,\partial_1f_1^\prime\otimes 1=\partial_1e_1^\prime\otimes 1$$
%and 
%$$\partial \left(U_1^\prime\otimes 1\right)= e_1^\prime\otimes 1+f_2^\prime\otimes 1 +
$$\partial_2U_1^\prime\otimes 1=-
\partial_2 U_2^\prime\otimes 1.$$
%=-e_2^\prime\otimes 1-f_1^\prime\otimes 1-\partial_2U_1^\prime\otimes 1.$$
%(One should keep in mind that $x\otimes 1=y\otimes 1$ means that $y=gx$ for some $g\in G$.)
That is, $$\partial_2U_1^\prime=g\overline{\partial_2U_2^\prime}$$ for some $g\in G$.

Since $U_1^\prime$ and $U_2^\prime$ belong to different $\kappa_{ij}$'s, say $\kappa_{ij_1}$ and $\kappa_{ij_2}$, we 
can, upon replacing $\kappa_{ij_2}$ by $g\kappa_{ij_2}$, assume that $\partial_2U_1^\prime=\overline{\partial_2U_2^\prime}$, that is, 
$U_1^\prime+U_2^\prime$ is a square. (Since $g$ 
maps $\partial e_2^\prime$ to $\partial e_1^\prime$, this second choice of $\kappa_{ij_2}$ in its $G$-orbit preserves the condition that
$e_2^\prime\in K^{str}\left(\partial_1Q\right)$.) 

%Finally, if vertices of $\kappa_ij_{1}$ or $\kappa_{ij_2}$ belong to $\partial_0Q\cap\partial_1Q$, 
%then we can, upon replacing 
%$\kappa_ij_{1}$ and $\kappa_{ij_2}$ with $g_2\kappa_ij_{1}$ and $g_2\kappa_{ij_2}$ for the same $g_2\in G$, assume that all vertices in $\partial_0Q\cap\partial_1Q$ belong to one of the sets $U_i$ given by construction 1 in Section 5.1.

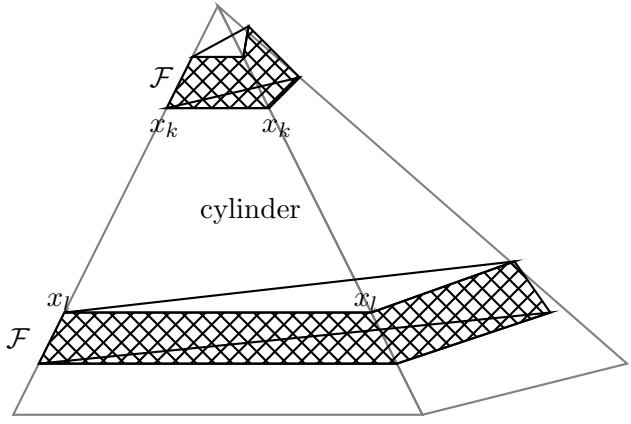
\psset{unit=0.1\hsize}
$$\pspicture(0,-1)(15,4)
\pspolygon[linecolor=gray](0,0)(4,0)(2,4)(0,0)
\pspolygon[linecolor=gray](4,0)(2,4)(6,0.5)(4,0)
\pspolygon(0.5,1)(3.5,1)(4.9,1.5)
\pspolygon(1.5,3)(2.5,3)(2.8,3.3)
\pspolygon(0.25,0.5)(3.75,0.5)(5.25,1)
\pspolygon(1.75,3.5)(2.25,3.5)(2.3,3.8)
\pspolygon[fillstyle=crosshatch](1.5,3)(2.5,3)(2.8,3.3)(2.3,3.8)(2.25,3.5)(1.75,3.5)
\pspolygon[fillstyle=crosshatch](0.25,0.5)(3.75,0.5)(5.25,1)(4.9,1.5)(3.5,1)(0.5,1)
\uput[0](-0.2,0.75){${\mathcal{F}}$}

\uput[0](1.2,3.3){${\mathcal{F}}$}
\uput[0](0.2,1.1){$x_l$}
\uput[0](3.2,1.1){$x_l$}
\uput[0](1.2,2.8){$x_k$}
\uput[0](2.3,2.8){$x_k$}
\uput[0](1.7,2){cylinder}

\endpspicture$$

Let $F$ resp.\ $F^\prime$ be the path-components of $\partial_1Q$ with $e_1^\prime\subset F$ resp.\ $e_2^\prime\subset F^\prime$. 
Then we have $\partial_1str\left(f_1^\prime\right),\partial_0str\left(f_2^\prime\right)\in F, \partial_0str\left(f_1^\prime\right),\partial_1str\left(f_2^\prime\right)\in F^\prime$. 

We note that $f_1^\prime$ and $f_2^\prime$
are edges with label 1. By condition (i) of Corollary 3,
this implies that $str\left(f_1^\prime\right)$ and $str\left(f_2^\prime\right)$ are distinguished
1-simplices.

By Condition ix) and Condition xiii) of Definition 9 we have that 
$$ \partial_1q\left(str\left(f_1^\prime\right)\right)=x_{E_0^F}=\partial_0q\left(str\left(f_2^\prime\right)\right), \partial_0
q\left(str\left(f_1^\prime\right)\right)=x_{E_0^{F^\prime}}=\partial_1q\left(str\left(f_2^\prime\right)\right).$$
%$q\left(str\left(f_1^\prime\right)\right)$ and 
%$q\left(str\left(f_2^\prime\right)\right)$ are distinguished 1-simplices with vertices in 
%$x_{E_0}$. 
That is, $q\left(str\left(e_1^\prime\right))\right)$ and $q\left(str
\left(e_2^\prime\right)\right)$ are loops in $\partial_1Q$, based at $x_{E_0^F}$ resp.\ $x_{E_0^{F^\prime}}$. 

%Moreover, $e_1^\prime,e_2^\prime\in K\left(\partial_1Q\right)$ implies, by Lemma 12 and Construction 1 that $q\left(str\left(e_1^\prime\right)\right),
%q\left(str\left(e_2^\prime\right)\right)\in K\left(\partial_1Q\right)$. 

Since the 
square $q\left(str\left(U_1^\prime+U_2^\prime\right)\right)$ realizes a homotopy between 
$q\left(str\left(f_1^\prime\right)\right)$ and $q\left(str\left(f_2^\prime\right)\right)$, we have that $$q\left(str\left(f_1^\prime\right)\right)=\gamma_1q\left(str\left(f_2^\prime\right)\right)\gamma_2$$
with $$\gamma_1=q\left(str\left(e_1^\prime\right))\right), \gamma_2=q\left(str\left(e_2^\prime\right))\right)\in \Omega\left(\partial_1Q\right)\subset \Gamma=\Omega\left(\partial Q\right).$$
By condition x) from  \hyperref[Def9]{Definition \ref*{Def9}} this implies  $$q\left(str\left(f_1^\prime\right)\right)=
q\left(str\left(f_2^\prime\right)\right).$$

%Moreover, $str\left(f_1^\prime\right)$ and 
%$str\left(f_2^\prime\right)$ have their respective endpoints in the same path-components 
%of $\partial Q$.
%By Definition 9, this implies that there exists $h\in H$ with $$q\left(str\left(f_1^\prime\right)\right)=h q\left(str\left(f_2^\prime\right)\right).$$

This means that $ q\left(str\left(U_1^\prime\right)\right)+ q\left(str\left(U_2^\prime\right)\right)$ is a cylinder with the boundary circles $q\left(str\left(e_1^\prime\right)\right)$ and $q\left(str\left(e_2^\prime\right)\right)$ in $\partial_1Q$.

(This is why we have performed the straightening construction 
in Section 5 such that 
there should be only one distinguished 1-simplex in each coset.) 
\\
%We have $e_1,e_2\in K_1\left(\partial N\right)$. Lemma 5, Lemma 11 vii) and the Remark in Section 5.1.\ imply $$q\left(str\left(e_i^\prime\otimes 1\right)\right)=q\left(str\left(r\left(e_i\otimes 1\right)\right)\right)\in HK_1\left(\partial_1Q\right)\otimes_{{\bf Z}H}{\bf Z}$$
%for $i=1,2$.  That is, $hU_1^\prime+hU_2^\prime$ is in the $H$-orbit of a cylinder with boundary in $\partial_1
%Q$. \\

The assumption $rmv\left(q\circ str\left(\kappa_{ij_1}\right)\otimes 1\right)\not=0$ implies that
the loops $ q\left(str\left(e_1^\prime\right)\right)$ and $ q\left(str\left(e_2^\prime
\right)\right)$ are not 0-homotopic.\\
Indeed,
if one of them, say $ q\left( str\left(e_1^\prime\right)\right)$, were a 0-homotopic (thus constant) loop, then also $ q\left( str\left(e_2^\prime\right)\right)$ would be a
0-homotopic (thus constant) loop, because they are homotopic through the cylinder.
%$str^\prime\left(f_1^\prime\right)=str^\prime\left(
%f_2^\prime\right)$ and $str^\prime\left(e_1^\prime\right)*str^\prime\left(f_1^\prime\right)*str^\prime
%\left(e_2^\prime\right)^{-1}* str^\prime\left(f_2^\prime\right)^{-1}$ 
%is 0-homotopic. \\
But  
$q\left(str\left(e_1^\prime\right)\right), q\left( str\left(e_2^\prime\right)\right)$ are edges of $q\left(str\left(\kappa_{ij_1}\right)\right)$ resp.\ $q\left(str\left(\kappa_{ij_2}\right)\right)$. 
In particular, $q\left(str\left(\kappa_{ij_1}\right)\right)$ and $q\left(str\left(\kappa_{ij_2}\right)\right)$ would have
a constant loop as an edge. 
By  \hyperref[Lemma14]{Lemma \ref*{Lemma14}} and  \hyperref[Def7]{Definition \ref*{Def7}}, this would prove the wanted equalities
$rmv\left(q\circ str\left(\kappa_{ij_1}\right)\otimes 1\right)=0, rmv\left(q\circ str\left(\kappa_{ij_2}\right)\otimes 1\right)=0$.

Thus we can assume that
$ q\left(str\left(e_1^\prime\right)\right)$ and $ q\left(str\left(e_2^\prime\right)\right)$ 
are not 0-homotopic, that is,
the cylinder $$ q\left(str\left(U_1^\prime\right)\right)+ q\left(str\left(U_2^\prime\right)\right)$$ is $\pi_1$-injective as a 
map of pairs. Since $\left(Q,\partial_1Q\right)$ is a 
pared acylindrical manifold, the cylinder must then 
be homotopic into $\partial Q$, as a map of
pairs 
$$\left({\bf S}^1\times \left[0,1\right],{\bf S}^1\times\left\{0,1\right\}\right)                              \rightarrow \left(Q,\partial_1Q\right).$$

Since $\partial_1Q$ is acylindrical, the cylinder must then either 
degenerate (${\bf S}^1\times\left[0,1\right]\rightarrow\partial Q$ homotopes to a map that factors over the projection ${\bf S}^1\times\left[0,1\right]\rightarrow{\bf S}^1$, in particular $q\left(str\left(e_1^\prime\right)\right)=q\left(str\left(e_2^\prime\right)
\right)$) or be homotopic into $\partial_0Q$ (and hence into $q\left(\partial_0Q\right)$, since $q\sim id$). In the second case 
the vertices $x_{E_0^F},x_{E_0^{F^\prime}}$ must belong to
$\partial_0Q$ and
%since $q\sim id$, 
we get by condition vii) from \hyperref[Lemma11]{Lemma \ref*{Lemma11}} that 
%$str\left(U_1^\prime+U_2^\prime\right)$ is freely homotopic into $\partial_0Q$. In particular 
$q\left(str\left(e_1^\prime\right)\right),q\left(str\left(e_2^\prime\right)\right)\in K_1^{str}\left(\partial_0Q\right)$. By \hyperref[Lemma15]{Lemma \ref*{Lemma15}} this 
implies that $q\left(str\left(\kappa_{ij_1}\right)\right)\otimes 1=0, q\left(str\left(\kappa_{ij_2}\right)\right)\otimes 1=0$. \\
Thus we can assume that the cylinder degenerates. In particular
$ q\left(str\left(f_1^\prime\right)\right), q\left(str\left(f_2^\prime\right)\right)\in K_1^{str}\left(\partial_1 Q\right)$.\\
%, and 
%thus also $str\left(f_1^\prime\right), str\left(f_2^\prime\right)$, are (homotopic into and therefore contained) in $\partial_1 Q$.
%By Lemma 11 vii) this implies
%$$str\left(f_1^\prime\right), str\left(f_2^\prime\right)\in K_1^{str}\left(\partial_1 Q\right).$$
\\
Let $P_1,P_2$ be the affine planes whose intersections with $\Delta^n$ contain $T_1$ resp.\ $T_2$. Let $W$ be the white component whose boundary contains the 
white-parallel arcs of $T_1,T_2$.
We have seen that there is are arcs $f_1,f_2$ connecting
$P_1\cap\Delta^n$ to $P_2\cap\Delta^n$ such that $$q\left(str\left(f_1^\prime\right)\right), q\left(str\left(f_2^\prime\right)\right)\in K_1^{str}\left(\partial_1 Q\right).$$ 
This
implies that for each other arc $f$ connecting
$P_1\cap\Delta^n$ to $P_2\cap\Delta^n$ the straightening $q\left(str\left(f^\prime\right)\right)$
must be (homotopic into and therefore by condition vii) from \hyperref[Lemma11]{Lemma \ref*{Lemma11}}) contained in $\partial_1 Q$. 
%Since $\partial_i str\left(f^\prime\right)$ belongs to the same component of $\partial Q$ as $\partial_istr\left(f_1^\prime\right)$, for $i=1,2$,
%the same argument as for $str\left(f_1^\prime\right)$ implies that $str\left(f^\prime\right)$ must be constant.

If $P_1$ and $P_2$
are of the same type (\hyperref[Def2]{Definition \ref*{Def2}}), then this shows that for all arcs $f\subset W$:
$$q\left(str\left(f^\prime\right)\right)\in K_1^{str}\left(\partial_1 Q\right)$$

If $P_1$ and $P_2$ are not of the same type, then the existence of a parallel arc
implies that at least one 
of them, say $P_1$, must be of type $\left\{0a_1\ldots
a_k\right\}$ with $k\not\in\left\{0,n-1\right\}$. Then, for each plane $P_3\not=P_1$ with $P_3\cap\Delta^n\subset \partial W$, it follows
from \hyperref[Cor2]{Corollary \ref*{Cor2}} that $P_3\cap\Delta^n$ has a white-parallel arc with $P_1\cap\Delta^n$.
Thus, repeating the
argument with $P_1$ and $P_3$ in place of $P_1$ and $P_2$, we prove that there are arcs in $\partial_1 Q$
connecting
$P_1\cap\Delta^n$ to $P_3\cap\Delta^n$, and consequently
for each arc
$f\subset W$ connecting
$P_1\cap\Delta^n$ to $P_3\cap\Delta^n$, the straightening $str\left(f^\prime\right)$
must be (homotopic into and therefore) contained in $\partial_1 Q$.

Consequently, also for arcs 
connecting
$P_2\cap\Delta^n$ to $P_3\cap\Delta^n$, we have that
$q\left(str\left(f^\prime\right)\right)$
must be (homotopic into and therefore) contained in $\partial_1 Q$. This finally shows that the 
1-skeleta of 
$q\left(str\left(\kappa_{ij_1}\right)\right)$ and $q\left(str\left(\kappa_{ij_2}\right)\right)$ belong to $K_1^{str}\left(\partial_1 Q\right)$.
By $\pi_1$-injectivity of $\partial_1Q\rightarrow Q$, asphericity of $K\left(\partial_1Q\right)$, and condition vii) from  \hyperref[Lemma11]{Lemma \ref*{Lemma11}}, this implies
that the
2-skeleta of
$q\left(str\left(\kappa_{ij_1}\right)\right)$ and $q\left(str\left(\kappa_{ij_2}\right)\right)$ belong to $K_1^{str}\left(\partial_1 Q\right)
$. Inductively, if the $k$-skeleta of
$q\left(str\left(\kappa_{ij_1}\right)\right)$ and $q\left(str\left(\kappa_{ij_2}\right)\right)$ belong to 
$K_k^{str}\left(\partial_1 Q\right)
$, then by asphericity of $K\left(Q\right)$, asphericity of $K\left(\partial_1Q\right)$, and 
condition vii) from \hyperref[Lemma11]{Lemma \ref*{Lemma11}} we obtain
that the $k+1$-skeleta of
$q\left(str\left(\kappa_{ij_1}\right)\right)$ and $q\left(str\left(\kappa_{ij_2}\right)\right)$ belong to 
$K_{k+1}^{str}\left(\partial_1 Q\right)
$. This provides the inductive step and thus our inductive proof shows that
$q\left(str\left(\kappa_{ij_1}\right)\right)$ and 
$q\left(str\left(\kappa_{ij_2}\right)\right)$
belong to $K^{str}\left(\partial_1Q\right)$. 

By \hyperref[Def7]{Definition \ref*{Def7}}, \hyperref[Def11]{Definition \ref*{Def11}} and \hyperref[Lemma14]{Lemma \ref*{Lemma14}} this implies $$rmv\left(
q\left(str\left(\kappa_{ij_1}\right)\right)\otimes 1\right)=0,
rmv\left(
q\left(str\left(\kappa_{ij_2}\right)\right)\otimes 1\right)=0.$$

So we have shown the subclaim: if $T_1=\partial_{k_1}\theta_{ij_1}, T_2=\partial_{k_2}\theta_{ij_2}$
have a white-parallel arc, then $rmv\left(
q\left(str\left(\kappa_{ij_1}\right)\right)\otimes 1\right)=0,
rmv\left(
q\left(str\left(\kappa_{ij_2}\right)\right)\otimes 1\right)=0$. In particular, $$q\left(str\left(T_1^\prime\right)\right),
q\left(str\left(T_2^\prime\right)\right)$$ do not occur (with non-zero coefficient)
in $$\partial \sum_{j\in J_i} rmv\left(q\left(str\left(\kappa_{ij}\right)\right)\otimes 1
\right).$$
\\
By \hyperref[Lemma10]{Lemma \ref*{Lemma10}}, for a canonical colouring associated to a
set of affine planes $P_1,P_2,\ldots$, and a fixed triangulation of each $Q_i=P_i\cap\Delta^n$,
we have at most $n+1$ n-1-simplices whose 1-skeleton does not contain a white-parallel arc.
Therefore the {\em subclaim} implies the {\em claim}. \\

Thus we have presented $\left[\partial Q\right]\otimes 1 $ as a finite chain of
$l^1$-norm at most\\
$ \left(n+1\right)\sum_{i=1}^r\mid a_i\mid$. By \hyperref[Lemma4]{Lemma \ref*{Lemma4}}a) we know that $G=\Pi\left(K\left(\partial_{0}Q\right)\right)$ is amenable. Hence $H=q_*\left(G\right)$ is amenable.
Thus \hyperref[Lemma7]{Lemma \ref*{Lemma7}}, applied to $X=\partial Q$ and $K=HK^{str}\left(\partial_1 Q\right)$ with its $H$-action, implies
$$\parallel \partial Q\parallel\le \left(n+1\right)\sum_{i=1}^r\mid a_i\mid.$$
QED\\
\\
\\
We remark that \hyperref[Thm1]{Theorem \ref*{Thm1}} is not true without assuming amenability of $\pi_1\partial_{0}Q$.
Counterexamples can be found, for example,
using \cite{jun} or \cite{k1}, Theorem 6.3.\\
\\
In \cite{ag}, \hyperref[Thm1]{Theorem \ref*{Thm1}}  has been proven for incompressible surfaces in hyperbolic 3-manifolds. We compare the steps of the proof in \cite{ag} with the arguments in our paper:\\
Step 1 in \cite{ag} is the normalization procedure, which we have restated in \hyperref[Lemma1]{Lemma \ref*{Lemma1}}. \\
Step 2 in \cite{ag} consists in choosing compatible triangulations of the
polytopes $\sigma_i^{-1}\left(N\right)$.\\
Step 3 in \cite{ag} boils down to the statement that, for each component $Q_i$ of $Q$,
there exists a retraction $r:\hat{N}\rightarrow p^{-1}\left(Q_i\right)$, for the covering $p:\hat{N}\rightarrow N$ corresponding to $\pi_1Q_i$.
Such a statement can not be correct because it would (together with step 7 from \cite{ag}) imply $\parallel N\parallel
\ge \parallel Q\parallel$ whenever $Q$ is a $\pi_1$-injective submanifold
of $N$. This inequality is true for
submanifolds with amenable boundary, but not in general. In fact, one only has the more complicated retraction
$r:C_*\left(K\left(N\right),
K\left(N^\prime\right)\right)\otimes_{{\bf Z}G}{\bf Z}
\rightarrow C_*\left(K\left(Q\right),K\left(\partial Q\right)
\right)\otimes_{{\bf Z}G}{\bf Z}$, with $G=\Pi\left(K\left(\partial_{0}Q\right)\right)$.
This more complicated retraction is the reason that much of the latter arguments
become notationally awkward, although conceptually not much is changing. Moreover, the action of
the group $G$ is basically the reason that  \hyperref[Thm1]{Theorem \ref*{Thm1}} is true only for amenable $G$.\\
Basically, the reason why the retraction $r:\hat{N}\rightarrow Q$ does not exist, is as follows. Let $R_j$ be the connected components of $\hat{N}-p^{-1}\left(Q_i\right)$.
Then $R_j$ is homotopy equivalent to each connected component of $\partial R_j$. If $\partial R_j$ were
connected for each $j$, this homotopy equivalence could be extended to a homotopy equivalence $r:\hat{N}\rightarrow p^{-1}\left(Q_i\right)$. However, in most cases $\partial R_j$
will be disconnected, and then such an $r$ can not exist.\\
We note that also the weaker construction of cutting off simplices does not work. A simplex may intersect $Q_i$ in many components and it is not clear which component to choose.\\
Step 5 in \cite{ag} is the straightening procedure, it corresponds to sections 5.2-5.4 in this paper.
We remark that the straightening procedure must be slightly more complicated than in \cite{ag} because it
is not possible, as suggested in \cite{ag}, to homotope
all edges between boundary components of $\partial Q$ into shortest geodesics.
This is the reason why we can only straighten chains with an admissible 0-1-labeling of their edges
(and why
our straightening homomorphism in Section \ref{sec:straighten} is only defined on $C_*^{simp}\left(\mid\Upsilon\mid\right)$ and not on all of
$C_*^{sing}\left(Q\right)$).\\
Step 6 in \cite{ag} consists in removing degenerate simplices. This corresponds to Section \ref{sec:remove} in this paper.\\
Step 7 in \cite{ag} proves that each triangle in $\sigma_i^{-1}\left(\partial N\right)$ contributes only once 
to the constructed fundamental cycle of $\partial Q$. Since, in our argument, we do not work with the covering $p:\hat{N}\rightarrow N$, we have no need
for this justification.\\
Step 8 in \cite{ag} counts the remaining triangles per simplex (after removing degenerate simplices).
It seems to have used the combinatorial arguments which we work out for arbitrary dimensions in Section 4.\\
We mention that the arguments of Section 4 are the only part of the proof which gets easier
if one restricts to
3-manifolds rather than arbitrary dimensions. Moreover, the proof for laminations
is the same as for hypersurfaces except for \hyperref[Lemma1]{Lemma \ref*{Lemma1}}. Thus, upon these two points it seems that even in the
case of incompressible surfaces in 3-manifolds the proof of \hyperref[Thm1]{Theorem \ref*{Thm1}} can not be further simplified.

\section{Specialization to 3-manifolds}

{\bf Guts of essential laminations.} We start with recalling the guts-terminology.
 Let $M$ be a compact 3-manifold with (possibly empty) boundary consisting of incompressible tori, and $\F$ an essential lamination
transverse or tangential to the boundary. $N=\ol$ is a, possibly noncompact,
irreducible 3-manifold with incompressible, aspherical boundary $\partial N$. We denote $\partial_0N=\partial
N\cap\partial M$ and $\partial_1N=\overline{\partial N - \partial_0N}$.
(Thus $\partial_1N$ is the union of boundary leaves of the lamination.)
By the proof of \cite{gk}, Lemma 1.3., the noncompact ends of $N$ are essential
$I$-bundles
over noncompact subsurfaces of $\partial_1N$. After cutting off
each of these ends along an essential, properly fibered annulus, one obtains a compact 3-manifold to which one can 
apply the JSJ-decomposition of \cite{js}, \cite{joh}. 
Hence we have a decomposition of $N$ into the characteristic submanifold 
$Char\left(N\right)$ (which consists of $I$-bundles 
and Seifert fibered 
solid tori, where the fibrations have to respect boundary patterns 
as defined in \cite{joh}, p.83) and the guts of $N$, $Guts\left(N\right)$. The $I$-fibered ends of $N$ will be added to the characteristic submanifold, which thus may become noncompact, while $Guts\left(N\right)$ is compact.
(We mention that there are
different notions of guts in the literature. Our notion is compatible with
\cite{ag}, \cite{ag2}, but differs from the definition in \cite{gk} or
\cite{cd} by taking the Seifert fibered solid tori into the 
characteristic submanifold and not into the guts. Thus, solid torus guts in the
paper of Calegari-Dunfield is the same as empty guts in our setting.) If $\partial_0N\cap \partial Q\not=\emptyset$ consists of annuli $A_1,\ldots,A_k$, then, to be consistent with the setting of Theorem 1, we add components $A_i\times\left[0,1\right]$ to $Char\left(N\right)$ (without changing the homeomorphism type of $N$), which implies $\partial_0N\cap\partial Q=\emptyset$.

For $Q=Guts\left(N\right)$ we denote $\partial_1 Q=\partial_1N\cap\partial Q=\partial
N\cap\partial Q=Q\cap\partial N$ and $\partial_0Q=\overline{\partial Q-\partial_1 Q}$.
For $R=Char\left(N\right)$ we denote $\partial_1 R=\partial
N\cap\partial R$ 
and $\partial_0R=\overline{\partial R-\partial_1 R}$. $\partial_0N\cap\partial Q=\emptyset$ implies then
$\partial_0Q=Q\cap R$.

$\partial_0Q$ consists of essential tori and annuli, in particular $\pi_1\partial_0Q$ is amenable.
The guts of $N$ has the following properties: the pair $\left(Q,\partial_1Q\right)$ is a pared acylindrical manifold as defined in Definition 3, 
$Q,\partial_1Q, \partial_1R$ are aspherical, 
and the inclusions $\partial_0Q\rightarrow Q, \partial_1Q\rightarrow Q, Q\rightarrow N$,
$\partial_0R\rightarrow R,\partial_1R\rightarrow R,R\rightarrow N$
are $\pi_1$-injective (see \cite{js},\cite{joh}).
It follows from Thurston's hyperbolization theorem for Haken manifolds that $Q$ admits a hyperbolic metric with geodesic boundary $\partial_1Q$ and 
cusps 
corresponding to $\partial_0Q$. (In particular, $\chi\left(\partial Q\right)\le 0$, thus $\partial Q$ is aspherical, and $\partial_1Q$ is a hyperbolic surface, thus acylindrical.)\\
\\
\\
{\bf Theorem 2} : {\em Let $M$ be a compact 3-manifold with (possibly empty) boundary consisting
of incompressible tori, and let
$\mathcal{F}$ be an essential lamination of
$M$. Then $$\parallel
M,\partial M\parallel^{norm}_{\mathcal{F}}\ge-\chi\left(Guts\left({\mathcal{F}}\right)
\right).$$
%, \parallel M\parallel_{\mathcal{F}} \ge -2\chi\left(Guts\left({\mathcal{F}}\right)\right).$$
More generally, if $P$ is a polyhedron with $f$ faces, then
 $$ \parallel
M,\partial M\parallel^{norm}_{{\mathcal{F}},P}\ge-\frac{2}{f-2}\chi\left(Guts\left({\mathcal
{F}}\right)\right).$$}

\begin{proof} 

Let $N=\overline{M-{\mathcal{F}}}$. Since $\mathcal{F}$ is essential, $N$ is irreducible (hence aspherical, since $\partial N\not=\emptyset$)
and has incompressible, aspherical boundary.
Let $R=Char\left(N\right)$ be the characteristic submanifold and
$Q=Guts\left(N\right)$ be the complement of the characteristic submanifold of $N$. The discussion before Theorem 2 shows that the decomposition $N=Q\cup R$
satisfies the assumptions of  \hyperref[Thm1]{Theorem \ref*{Thm1}}.
%We denote $\partial_1 Q=\partial
%N\cap\partial Q,
%\partial_0Q=\overline{\partial Q-\partial_1 Q}=Q\cap R$,
%and
%$\partial_1 R=\partial N\cap\partial R.$
%It follows from the construction of the
%characteristic submanifold (see \cite{js},\cite{joh}) that\\
%- $Q$ is compact and aspherical,\\
%- $\left(Q,\partial_1Q\right)$ is pared acylindrical, $\partial_1Q$  and $\partial_1R$ are aspherical, $\partial_1Q\rightarrow Q\rightarrow N$ and $\partial_1R\rightarrow R\rightarrow N$ are
%$\pi_1$-injective, and\\
%- $\partial_0Q$ consists of essential tori and annuli, in particular $\pi_1\partial_0Q$ is amenable and injects into $\pi_1 Q$ and $\pi_1R$.

%(Even stronger, Thurston has proven that $Q$ is a hyperbolic manifold with geodesic boundary $\partial_1Q$ and cusps corresponding to $\partial_0Q$. In particular, $\chi\left(\partial Q\right)\le 0$, thus $\partial Q$ is aspherical, and $\partial_1Q$ is a hyperbolic surface, thus acylindrical.)

%In conclusion, the decomposition $N=Q\cup R$ satisfies the assumptions of Theorem 1.

From the computation of the simplicial volume for surfaces
(\cite{gro}, section 0.2.) and
$\chi\left(Q\right)=\frac{1}{2}\chi\left(\partial Q\right)$ (which is a consequence of
Poincare duality for the closed 3-manifold $Q\cup_{\partial Q}Q$), it follows that
$$-\chi\left(Guts\left({\mathcal{F}}\right)\right)=
-\frac{1}{2}\chi\left(\partial Guts\left({\mathcal{F}}\right)\right)=\frac{1}{4}\parallel \partial
Guts\left({\mathcal{F}}\right)\parallel.$$
Thus, the first claim is obtained as application of Theorem 1 to
$Q=Guts\left({\mathcal{F}}\right)$. 

The second claim, that is the generalisation to arbitrary polyhedra, is
obtained as in \cite{ag}. 
Namely, one
uses the same straightening as above, and asks again how many nondegenerate 
2-simplices may, after straightening, occur in the intersection of 
$\partial Q$ with some polyhedron $P_i$. In \cite{ag}, p.\ 11, it is shown that this
number is at most $2f-4$, where $f$ is the number of faces of $P_i$. The same argument as above shows then $\sum_{i=1}^r\mid a_i\mid\ge\frac{1}{2f-4}\parallel
\partial Guts\left({\mathcal{F}}\right)\parallel$, giving the wanted inequality.
\end{proof}

\noindent The following corollary applies, for example, to all hyperbolic manifolds
obtained by Dehn-filling the complement of the figure-eight knot in ${\bf S}^3
$. (Note that Hatcher has proved in \cite{hat} that each hyperbolic manifold
obtained by Dehn-filling the complement of the figure-eight knot in ${\bf S}^3
$ carries essential laminations.)

\begin{cor}\label{Cor4} If $M$ is a finite-volume hyperbolic manifold with $Vol\left(M\right)< 2V
_3 =2.02...$, then $M$ carries no essential lamination $\mathcal{F}$
with $\parallel M,\partial M\parallel^{norm}_{{\mathcal{F}},P}=\parallel M,\partial M\parallel_P$ for
 all polyhedra P, and nonempty guts. In particular, there is no tight essential
lamination with nonempty guts.\end{cor}

\begin{proof} The derivation of Corollary 4 from \hyperref[Thm2]{Theorem \ref*{Thm2}} is exactly
the same as in \cite{ag} for the usual (non-laminated)
Gromovnorm. Namely, by \cite{stt} (or \cite{ag}, end of Section 6) there exists a sequence $P_n$ of straight polyhedra in ${\bf H}^3$ with $\lim_{n\rightarrow\infty} \frac{Vol\left(P_n\right)}{f_n-2}=V_3$, with $f_n$ denoting
the number of faces of $P_n$. Assuming
that $M$ carries a lamination $\mathcal{F}$ with $ \parallel M,\partial M\parallel^{norm}_{{\mathcal{F}},P_n}=\parallel M,\partial M\parallel_{P_n}$
for all $n$, one gets $$-\chi\left(Guts\left({\mathcal{F}}\right)\right)\le
\frac{f_n-2}{2}\parallel  M,\partial M\parallel_{{\mathcal{F}},P_n}=
\frac{f_n-2}{2}\parallel  M,\partial M\parallel_{P_n}$$
$$\le \frac{f_n-2}{2}\frac{Vol\left(M\right)}{Vol\left(P_n\right)} 
\rightarrow \frac{Vol\left(M\right)}{2V_3} < 1.$$
On the other hand, if $Guts\left({\mathcal{F}}\right)$ is not empty, then it is
a hyperbolic manifold with nonempty geodesic boundary, hence $$
\chi\left(Guts\left({\mathcal{F}}\right)\right)\le -1,$$
giving a contradiction.\end{proof}

\begin{df}\label{Def12} The Weeks manifold is the closed 3-manifold obtained by $\left(-\frac{5}{1},-\frac{5}{2}\right)$-surgery at the Whitehead link (\cite{ro}, p.68).\end{df}

It is known that the Weeks manifold is hyperbolic and that its hyperbolic volume is
approximately $0.94..$. (It is actually the hyperbolic 3-manifold of smallest volume.)

\begin{cor}\label{Cor5} (\cite{cd}, Conjecture 9.7.): The Weeks manifold admits no tight
lamination $\mathcal{F}$.\end{cor}
\begin{proof} According to \cite{cd}, the Weeks manifold can not carry
a tight lamination with empty guts. Since tight laminations satisfy
$\parallel M\parallel^{norm}_{{\mathcal{F}},P}=\parallel M\parallel$ for each 
polyhedron (see \hyperref[Lemma1]{Lemma \ref*{Lemma1}}),
and since the
Weeks manifold has volume smaller than $2V_3$, it follows from  \hyperref[Cor4]{Corollary \ref*{Cor4}} that
it can not carry a tight lamination with nonempty guts neither.\end{proof}

The same argument shows that a hyperbolic 3-manifold $M$ with \\
- $Vol\left(M\right) < 2V_3$, and\\
- no injective homomorphism $\pi_1M\rightarrow Homeo^+\left({\bf S}^1\right)$\\
can not carry a tight lamination, because it was shown by Calegari-Dunfield
in \cite{cd} that the existence of a tight lamination with empty guts implies the existence of an injective homomorphism $\pi_1M\rightarrow Homeo^+\left({\bf S}^1\right)$. Some methods for excluding the existence 
of injective homomorphisms $\pi_1M\rightarrow Homeo^+\left({\bf S}^1\right)$
have been developed in \cite{cd} (which yielded in particular the nonexistence
of such homomorphisms for the Weeks manifold, used in the corollary above), but
in general it is still hard to apply this criterion to other hyperbolic 3-manifolds of volume $< 2V_3$.

As indicated in \cite{cal2}, an approach to
a generalization of some of the above arguments 
to essential, non-tight laminations, yielding possibly a proof for nonexistence of essential laminations on the Weeks manifold, could consist in trying to define a straightening of cycles (as in the proof of \hyperref[Lemma1]{Lemma \ref*{Lemma1}}) upon possibly changing the essential lamination. 

As a consequence of a recent paper of Tao Li, one can at least exclude the 
existence 
of transversely orientable essential laminations on the Weeks manifold.  
\begin{cor}\label{Cor6} The Weeks manifold admits no transversely orientable essential lamination $\mathcal{F}$.\end{cor}
\begin{proof} According to \cite{li}, Theorem 1.1, the following statement is true: if a closed, orientable, atoroidal 3-manifold $M$ contains a transversely orientable essential lamination, then it contains a transversely orientable tight essential lamination. Hence Corollary 6 is a direct consequence of \hyperref[Cor5]{Corollary \ref*{Cor5}}.\end{proof}

\section{Higher dimensions}

We want to finish this paper with showing that \hyperref[Thm1]{Theorem \ref*{Thm1}} is interesting also in higher dimensions. While
in dimension 3 the assumptions of \hyperref[Thm1]{Theorem \ref*{Thm1}} hold for each essential lamination, it is likely that this will not be the case for many laminations 
in higher dimensions. However, the most straightforward, but already interesting  
application of the inequality is  \hyperref[Cor7]{Corollary \ref*{Cor7}} which
means that, for a given negatively curved manifold $M$, we can give an explicit bound on the topological complexity of geodesic hypersurfaces. Such a bound seems to be new except, of course, in the 3-dimensional case where it is due to Agol (\cite{ag}) and (with nonexplicit constants) to Hass (\cite{ha}).

\begin{cor}\label{Cor7} Let $M$ be a compact Riemannian $n$-manifold of negative sectional curvature and finite volume. Let $F\subset M$ be a geodesic $n-1$-dimensional
hypersurface of finite volume. Then $\parallel
F\parallel\le \frac{n+1}{2}\parallel M\parallel$.\end{cor}

\begin{proof} Consider $N=\overline{M-F}$. $\left(N,\partial N\right)$
is acylindrical. This is well-known and
can be seen as follows: assume that $N$ contained an essential
cylinder, then the double $DN=N\cup_{\partial_1 N}N$ would contain an essential
2-torus.
% and, in particular, its fundamental group contained
%a non-peripheral subgroup isomorphic to ${\bf Z}^2$.
But, since $N$ is a negatively curved manifold with geodesic boundary,
we can glue the Riemannian metrics to get a complete negatively curved Riemannian metric
on $DN$. In particular, $DN$ contains no essential 2-torus, giving a contradiction.
%$\pi_1\left(DN\right)$ must be (relatively) word-hyperbolic (relative to $\pi_1\left(\partial DN\right)$)
%and can thus not contain any non-peripheral subgroup isomorphic to ${\bf Z}^2$.

Moreover, the geodesic boundary $\partial N$ is $\pi_1$-injective and negatively curved, thus aspherical. 
Therefore we can choose $Q=N$, in which case the other assumptions of Theorem 1
are
trivially satisfied. From \hyperref[Thm1]{Theorem \ref*{Thm1}} we conclude \\
$\parallel M\parallel_{F}^{norm}\ge\frac{1}{n+1}\parallel \partial N\parallel$. The boundary of $N$ consists of two copies of $F$, hence $\parallel \partial N\parallel=2\parallel F\parallel$. The leaf space of $\widetilde{F}\subset\widetilde{M}$ is a Hausdorff tree, thus \hyperref[Lemma1]{Lemma \ref*{Lemma1}}b) implies
$\parallel M\parallel_{F}^{norm}=\parallel M\parallel$. The claim follows.\end{proof}

This statement should be read as follows: for a given manifold $M$ (with given volume) one has an upper bound on the topological complexity of compact geodesic hypersurfaces.

%We finish this section by a reformulation of Corollary 7 for hyperbolic manifolds.
For hyperbolic manifolds one can use the proportionality principle and the Chern-Gau\ss -Bonnet Theorem 
to reformulate \hyperref[Cor7]{Corollary \ref*{Cor7}} as follows: If $M$ is a closed hyperbolic n-manifold 
and $F$ a closed n-1-dimensional
geodesic hypersurface, then $Vol\left(M\right)\ge
C_n\chi\left(F\right)$ for a constant $C_n$ depending only on $n$.

%\begin{cor}: Let $M$ be a closed hyperbolic n-manifold and $F$ a closed n-1-dimensional
%geodesic hypersurface. Then $$Vol\left(M\right)\ge
%C_n\chi\left(F\right)$$ for a constant $C_n$ depending only on $n$.\end{cor}

%\begin{pf}
%Assume $n$ is odd. (If $n$ is even, then $\chi\left(F\right)=0$ and the claim is trivially true.) Then, $F$ is an even-dimensional hyperbolic manifold and, by the Chern-Gau\ss -Bonnet-Theorem we have $Vol\left(F\right)=G_{n-1}\chi\left(F\right)$
%for a constant $G_{n-1}$. Together with $Vol\left(M\right)=V_n\parallel M\parallel,
%Vol\left(F\right)=V_{n-1}\parallel F\parallel$ and Corollary 7,
%we conclude $Vol\left(M\right)\ge \frac{2}{n+1}\frac{V_n}{V_{n-1}}G_{n-1}\chi\left(F\right)$.\end{pf}

Thilo Kuessner\\
Mathematisches Institut, Universit\"at M\"unster\\
Einsteinstra\ss e 62\\
D-48149 M\"unster\\
Germany\\
e-mail: kuessner@math.uni-muenster.de


\begin{thebibliography}{29}
\bibitem{ag}I.\ Agol, {\em Lower bounds on volumes of hyperbolic Haken 3-manifolds}, Preprint.  http://front.math.ucdavis.edu/math.GT/9906182.
\bibitem{ag2}I.\ Agol, P.\ Storm, W.\ Thurston, with an appendix by N.\ Dunfield, {\em Lower bounds on volumes of 
hyperbolic Haken 3-manifolds}, J.\ Amer.\ Math.\ Soc.\ {\bf 20}, pp.1053-1077 (2007).
\bibitem{bp}R.\ Benedetti, C.\ Petronio, {\em Lectures on
Hyperbolic Geometry}. Universitext, Springer-Verlag, Berlin (1992).
\bibitem{bri}M.\ Brittenham, {\em Essential laminations and Haken normal form}, 
Pac.\ J.\ Math.\ {\bf 168}, 217-234 (1995).
\bibitem{cal}D.\ Calegari, {\em The Gromov norm and foliations}.
GAFA 11, pp.1423-1447 (2001).
\bibitem{cal2}D.\ Calegari, {\em Problems in foliations and laminations on 3-manifolds}.
 Proc.\ Sympos.\ Pure Math.\ {\bf 71}, 297-335 (2003).
\bibitem{cd}D.\ Calegari, N.\ Dunfield, {\em Laminations and groups of homeomorphisms
of the circle}, Invent.\ Math.\ {\bf 152}, 149-207 (2003).
\bibitem{cs}M.\ Culler, P.\ Shalen, {\em Volumes of hyperbolic Haken manifolds}, 
Invent.\ Math.\ {\bf 118}, pp.285-329 (1994).
%\bibitem{dt}A.\ Dold, R.\ Thom, {\em Quasifaserungen und unendliche symmetrische Produkte}, Ann.\ Math.\ (2) {\bf 67}, pp.239-281 (1958).
\bibitem{fe}S.\ Fenley, {\em Laminarfree hyperbolic manifolds}, Comm.\ Math.\ Helv.\ {\bf 82}, pp.247-321 (2007).
\bibitem{fr}S.\ Francaviglia, {\em Hyperbolic volume of representations of 
fundamental groups of cusped 3-manifolds}, Int.\ Math.\ Res.\ Not.\ {\bf 9}, 425-459 (2004).
\bibitem{gab}D.\ Gabai, {\em Essential laminations and Kneser normal form}, J.\ 
Diff.\ Geom.\ {\bf 53}, 517-574 (1999).
\bibitem{gk}D.\ Gabai, W.\ Kazez, {\em Group negative curvature for 3-manifolds with genuine laminations}, Geom.\ Topol.\ {\bf 2}, 65-77 (1998).
\bibitem{goe}D.\ Gabai, U.\ Oertel, {\em Essential laminations in 3-manifolds}, 
Ann.\ of Math.\ {\bf 130}, 41-73 (1989).
\bibitem{gro}M.\ Gromov, {\em Volume and Bounded
Cohomology}, Public.\ Math.\ IHES {\bf 56}, 5-100 (1982).
\bibitem{ha}J.\ Hass, {\em Acylindrical surfaces in 3-manifolds}, Michigan Math.\ J.\ 42, pp.357-365 (1995).
\bibitem{hat}A.\ Hatcher, {\em Some examples of essential laminations in 
3-manifolds}, Ann.\ Inst.\ Fourier {\bf 42}, 313-325 (1992).
\bibitem{iva}N.\ Ivanov, {\em Foundations of the theory of
bounded cohomology}, J.\ Sov.\ Math.\ {\bf 37}, 1090-1114 (1987).
\bibitem{js}W.\ Jaco, P.\ Shalen, {\em Seifert fibered spaces in 3-manifolds}, 
Mem.\ Amer.\ Math.\ Soc., {\bf 21}, (1979).
\bibitem{joh}K.\ Johannson, {\em Homotopy Equivalences of 3-Manifolds with Boundaries},
Lecture Notes in Mathematics 761, Springer Verlag, (1979).
\bibitem{jun}D.\ Jungreis, {\em Chains that realize the Gromov
invariant of hyperbolic manifolds},  Ergodic Theory and Dynamical Systems {\bf 17},
643-648 (1997).
\bibitem{k1}T.\ Kuessner, {\em Efficient fundamental cycles for
cusped hyperbolic manifolds}, Pac.\ J.\ Math.\ {\bf 211}, 283-314 (2003).
h\bibitem{k2}T.\ Kuessner, {\em Multicomplexes, 
bounded cohomology and additivity 
of simplicial volume}, Preprint, 
http://www.math.uni-muenster.de/reine/u/kuessner/preprints/bc.pdf
\bibitem{k3}T.\ Kuessner, {\em Foliated norms on fundamental group and homology}, Top.\ Appl.\ {\bf 144}, 229-254 (2004).
\bibitem{li}T.\ Li, {\em Compression branched surfaces and tight essential laminations}, Preprint.
\bibitem{may}J.\ P.\ May, {\em Simplicial objects in algebraic topology} (Chicago Lectures in mathematics, UCP, 1992).
\bibitem{rs}R.\ Roberts, M.\ Stein, {\em Group actions on order trees}, Top.\ Appl.\ {\bf 115}, 175-201 (2001).
\bibitem{rss}R.\ Roberts, J.\ Shareshian, M.\ Stein, {\em Infinitely many hyperbolic 3-manifolds which contain no Reebless foliation}, J.\ Amer.\ Math.\ Soc.\ {\bf 16}, 639-679 (2003).
\bibitem{ro}D.\ Rolfsen, {\em Knots and Links}, AMS Chelsea Publishing (2003).
\bibitem{stt}D.\ Sleator, R.\ Tarjan, W.\ Thurston, {\em Rotation distance, triangulations, and hyperbolic geometry}, J.\ Amer.\ Math.\ Soc.\ {\bf 1}, 647-681 (1988).
\bibitem{thu}W.\ Thurston, {\em The Geometry and Topology of 3-Manifolds
}. Lecture Notes, Princeton.
\bibitem{thua}W.\ Thurston, {\em Hyperbolic geometry and 3-manifolds}.
Low-dimensional topology, Proc.Conf.Bangor, Lond.Math.Soc.Lect.NoteSer.48, pp.9-
25 (1982).
\end{thebibliography}
\end{document}